\documentclass[10pt]{amsart}%{scrartcl,scrreprt}
\usepackage[T1]{fontenc}
\usepackage{geometry}
\usepackage[latin1] {inputenc}
\usepackage{amsmath}
\usepackage{amsfonts, amssymb, textcomp}
\usepackage[colorlinks=flase, linkcolor=red,urlcolor=green, citecolor=blue]{hyperref}
\usepackage{subeqnarray}

\usepackage{latexsym}
\usepackage{fancyhdr}
\usepackage{longtable}
\usepackage{amsmath, amssymb}
\usepackage{graphicx}
\usepackage{enumerate}

\numberwithin{equation}{section}

\theoremstyle{plain}
\newtheorem{proposition}{Proposition}[section]
\newtheorem{theorem}[proposition]{Theorem}
\newtheorem{lemma}[proposition]{Lemma}
\newtheorem{corollary}[proposition]{Corollary}
\newtheorem{definition}[proposition]{Definition}

\newtheorem{remark}[proposition]{Remark}

\newcommand{\RR}{\mathbb{R}}
\newcommand{\CC}{\mathbb{C}}
\newcommand{\NN}{\mathbb{N}}

\newcommand{\id}{\operatorname{id}}

\let\on=\operatorname
\newenvironment{demo}[1]{\par\smallskip\noindent{\bf #1.}}{\par\smallskip}

\makeatletter
\@namedef{subjclassname@2020}{%
  \textup{2020} Mathematics Subject Classification}
\makeatother

\newsavebox{\fmbox}
\newenvironment{fmpage}[1]
 {\begin{lrbox}{\fmbox}\begin{minipage}{#1}}
 {\end{minipage}\end{lrbox}\fbox{\usebox{\fmbox}}}

\title[Generalized upper and lower Legendre conjugates]
{Generalized upper and lower Legendre conjugates for Braun-Meise-Taylor weight functions}

\author[G.~Schindl]{Gerhard Schindl}

\address{G.~Schindl: Fakult\"at f\"ur Mathematik, Universit\"at Wien, Oskar-Morgenstern-Platz~1, A-1090 Wien, Austria.}
\email{gerhard.schindl@univie.ac.at}

\begin{document}

\begin{abstract}
We apply recent knowledge and techniques of the new generalized upper and lower Legendre conjugates to the theory of weight functions in the sense of Braun-Meise-Taylor and study in detail the effects on the corresponding associated weight matrices. An immediate and concrete application of the main statements is also provided. More precisely, we generalize a very recent result concerning the continuity and the range of the resolvent operator when being considered on weighted spaces of globally defined functions of Gelfand-Shilov type.
\end{abstract}

\thanks{This research was funded in whole or in part by the Austrian Science Fund (FWF) 10.55776/PAT9445424}
\date{\today}

\maketitle

\section{Introduction}
In this work we continue the study about the new conjugates $\sigma\check{\star}\tau$ and $\sigma\widehat{\star}\tau$ for weight functions $\sigma$, $\tau$ from our recent paper \cite{genLegendreconj}; see Section \ref{generalizedLegendresection} for the precise definitions. $\sigma\check{\star}\tau$ denotes the \emph{generalized lower Legendre conjugate} whereas $\sigma\widehat{\star}\tau$ is the \emph{generalized upper Legendre conjugate.} Recall that the known Legendre conjugates correspond, up to an inversion of the variable, to the case $\tau=\id: t\mapsto t$. In order to proceed, $\sigma$ and $\tau$ are required to satisfy certain growth conditions and in \cite{genLegendreconj} we have considered a general notion for being a weight function, more precisely it turns out that it is sufficient to assume that $\sigma,\tau$ are weight functions in the general setting of \cite{index}; see Definition \ref{weightfctdef}. Then $\sigma\check{\star}\tau$ is again a weight function but, however, in order to ensure that $\sigma\widehat{\star}\tau$ is a weight function, $\sigma$ and $\tau$ have to satisfy a certain growth relation and technical conditions have to be investigated; we refer to \cite[Sect. 4.2]{genLegendreconj}.\vspace{6pt}

In \cite[Sect. 5]{genLegendreconj} we have then considered a more special case for weight functions, namely $\sigma=\omega_{\mathbf{M}}$, $\tau=\omega_{\mathbf{N}}$ with $\mathbf{M},\mathbf{N}\in\RR_{>0}^{\NN}$ being weight sequences satisfying certain growth properties and $\omega_{\mathbf{M}}$, $\omega_{\mathbf{N}}$ are denoting the corresponding \emph{associated weight functions;} see Section \ref{assofunctionsection} for details. It turns out that within this setting more detailed information can be obtained which is due to the fact that associated weight functions satisfy more growth and regularity properties and the underlying weight sequences admit the possibility to express desired properties in terms of $\mathbf{M}$, $\mathbf{N}$ directly. We refer to the main results \cite[Thm. 5.4 \& 5.11]{genLegendreconj} where it is shown that the operations $\sigma\check{\star}\tau$ and $\sigma\widehat{\star}\tau$ are precisely corresponding to the point-wise product and point-wise quotient of sequences, respectively. Also in the weight sequence setting for introducing $\omega_{\mathbf{M}}\widehat{\star}\omega_{\mathbf{N}}$ the sequences $\mathbf{M}$, $\mathbf{N}$ are required to satisfy a certain growth relation since in general this function is not well-defined anymore; see \cite[Sect. 5.4]{genLegendreconj}.\vspace{6pt}

The aim of this recent work is now to study $\sigma\check{\star}\tau$ and $\sigma\widehat{\star}\tau$ for more specific \emph{weight functions in the sense of Braun-Meise-Taylor;} see \cite{BraunMeiseTaylor90}, \cite{BonetMeiseMelikhov07} and \cite{compositionpaper} for instance and more details can be found in Section \ref{BMTdefsection}. In this situation the weights satisfy additional growth and regularity properties and, indeed, associated weight functions satisfy some of the basic requirements in this setting automatically; see Lemma \ref{assoweightomega0}.

In \cite{compositionpaper} and \cite{dissertation} we have introduced and investigated the new notion of \emph{associated weight matrices:} To each Braun-Meise-Taylor weight function $\omega$ one can associate a weight matrix $\mathcal{M}_{\omega}$, see Section \ref{assomatrixsection} for details, and via using $\mathcal{M}_{\omega}$ one is able to study and define corresponding weighted spaces equivalently. In \cite{compositionpaper} this has been exclusively done for \emph{ultradifferentiable function classes of Roumieu- and Beurling type $\mathcal{E}_{\{\omega\}}$ and $\mathcal{E}_{(\omega)}$.} However, the crucial result \cite[Lemma 5.9]{compositionpaper} only involves growth properties for the weight, see estimate \eqref{newexpabsorb}, and no information concerning the underlying weighted space is required. Since in other weighted settings analogous defining expressions are used, the corresponding statement immediately holds for more weighted spaces: The symbol/functor $\mathcal{B}$ denotes \emph{globally defined ultradifferentiable classes,} $\mathcal{A}$ is the category of \emph{ultraholomorphic functions,} $\mathcal{S}$ are weighted globally defined function classes of \emph{Gelfand-Shilov type,} $\Lambda$ is standing for \emph{weighted spaces of sequences of complex numbers} and $\mathcal{F}$ denotes \emph{weighted formal power series with complex coefficients.} In this work, except for $\mathcal{S}$ in Section \ref{applicationlowersection}, we are not defining the weighted spaces explicitly since we are mainly focusing on weights and their growth properties and hence our new results admit applications in different directions and settings. For the precise definition of the category $\mathcal{E}$ we refer (again) to \cite{BraunMeiseTaylor90}, \cite{BonetMeiseMelikhov07} and \cite{compositionpaper}, for $\mathcal{A}$ see \cite{sectorialextensions1}, for $\mathcal{S}$ see \cite{nuclearglobal2}.

Note that weighted spaces appear pair-wise: One considers the \emph{Roumieu type} $\mathfrak{F}_{\{\omega\}}$ and the \emph{Beurling type} $\mathfrak{F}_{(\omega)}$ with $\mathfrak{F}\in\{\mathcal{E}, \mathcal{B}, \mathcal{A}, \mathcal{S}, \Lambda, \mathcal{F}\}$.\vspace{6pt}

We are interested in and concerned with the study of the effects of the operations $\check{\star}$ and $\widehat{\star}$ for the corresponding associated weight matrices and ask:

\begin{itemize}
\item[$(*)$] \emph{How do the mappings $(\sigma,\tau)\mapsto\sigma\check{\star}\tau$ and $(\sigma,\tau)\mapsto\sigma\widehat{\star}\tau$ modify the matrices $\mathcal{M}_{\sigma}$, $\mathcal{M}_{\tau}$?}

\item[$(*)$]\emph{How can the weighted spaces defined in terms of $\sigma\check{\star}\tau$ and $\sigma\widehat{\star}\tau$ be expressed via $\mathcal{M}_{\sigma}$, $\mathcal{M}_{\tau}$?}

\item[$(*)$] \emph{Are these induced modifications allowing for some applications?}
\end{itemize}

Note that, when $\sigma$ and $\tau$ are given, then $\sigma\check{\star}\tau$ and $\sigma\widehat{\star}\tau$ define related weight structures and hence it is natural to investigate the ``controlled transformation of regularity'' w.r.t. these conjugates.

As expected and suggested by the statements \cite[Thm. 5.4 \& 5.11]{genLegendreconj}, also for Braun-Meise-Taylor weight functions the point-wise product(s) resp. point-wise quotient(s) of the sequences belonging to the associated matrices is becoming relevant and give rise to the notion of products resp. quotients of (general) weight matrices. Since the associated weight matrices provide a convenient alternative method to describe weighted classes expressed by Braun-Meise-Taylor weight functions, the study of these questions is of particular interest for explicit applications. Due to the general setting and flavor of our main statements one can expect numerous applications in different weighted settings. Indeed, already in this recent article we provide a generalization of the very recent main statement \cite[Thm. 4.6 $2.$]{ArizaFernandezGalbis25} dealing with the continuity of the resolvent operator considered on weighted spaces of globally defined functions of Gelfand-Shilov type. This question is related to the study of \emph{dynamics of certain composition operators in the Gelfand-Shilov setting} (being defined in terms of a polynomial). More precisely, we are able to replace the special Gevrey-weights treated in \cite{ArizaFernandezGalbis25} by more general Braun-Meise-Taylor weight functions $\sigma$, $\tau$ and when involving the conjugates $\check{\star}$, $\widehat{\star}$.

Finally, we mention that the non-standard situations in the weight sequence setting treated and studied in detail in \cite{genLegendreconj}, more precisely dealing with $\mathbf{M}$ such that
$\mathbf{M}_{\iota}:=\liminf_{p\rightarrow+\infty}\left(\frac{M_p}{M_0}\right)^{1/p}=\liminf_{p\rightarrow+\infty}(M_p)^{1/p}<+\infty$, are not becoming relevant in this recent work. This is due to the fact that the sequences under consideration are frequently belonging to associated weight matrices and hence to the ``regular set of sequences'' $\mathcal{LC}$; see Section \ref{weightsequencesection} and $(i)$ in Section \ref{assomatrixsection}.\vspace{6pt}

The paper is structures as follows: First, in Section \ref{weightsection} we gather and list all required technical growth properties and notation for weight sequences, (associated) weight functions and (associated) weight matrices, and revisit also the generalized lower and upper Legendre conjugate and crucial growth indices for weight functions. Then, in Section \ref{lowersection} we study the effect of $\check{\star}$ and prove the main results Theorem \ref{mainweighfctthm} and Corollary \ref{mainweighfctcor}. In Section \ref{applicationlowersection} we give an application of these statements in the Gelfand-Shilov setting and study the behavior of the resolvent operator; see Theorem \ref{ArizaFernandezGalbisthm}. In Sections \ref{uppersection} and \ref{applicationuppersection} the analogous statements are shown for $\widehat{\star}$; see the main results Theorems \ref{mainweighfctcorinv}, \ref{mainweighfctinv}, and \ref{ArizaFernandezGalbisthminv}. However, again the upper Legendre conjugate is much more technical and involved and first, in Section \ref{welldefinedsection}, we have to provide a detailed study in terms of weight functions and their associated weight matrices in order to ensure that $\sigma\widehat{\star}\tau$ is well-defined.

\section{Weights, growth conditions and definitions}\label{weightsection}
We write $\NN=\{0,1,2,\dots\}$ and $\NN_{>0}:=\{1,2,\dots\}$ and set $\RR_{>0}:=(0,+\infty)$.

\subsection{Weight sequences}\label{weightsequencesection}
Given a sequence $\mathbf{M}=(M_p)_p\in\RR_{>0}^{\NN}$ we also use the notation $\mu=(\mu_p)_p$ with $\mu_p:=\frac{M_p}{M_{p-1}}$, $p\ge 1$, $\mu_0:=1$, and analogously for all other arising sequences. $\mathbf{M}$ is called \emph{normalized,} if $1=M_0\le M_1$ is valid. For the \emph{point-wise product sequence} we simply write $\mathbf{M}\cdot\mathbf{N}$ and $\frac{\mathbf{M}}{\mathbf{N}}$ for the \emph{point-wise quotient sequence.}\vspace{6pt}

$\mathbf{M}$ is called \emph{log-convex} if
$$\forall\;p\in\NN_{>0}:\;M_p^2\le M_{p-1} M_{p+1},$$
equivalently if $\mu$ is nondecreasing. If $\mathbf{M}$ is log-convex and normalized, then $\mu_p\ge 1$ for all $p\in\NN$, $p\mapsto M_p$ and $p\mapsto(M_p)^{1/p}$ are nondecreasing, $(M_p)^{1/p}\le\mu_p$ for all $p\in\NN_{>0}$ and finally $M_pM_q\le M_{p+q}$ for all $p,q\in\NN$; see e.g. \cite[Lemmas 2.0.4 \& 2.0.6]{diploma}.\vspace{6pt}

For our purpose it is convenient to consider the following set of sequences
$$\hypertarget{LCset}{\mathcal{LC}}:=\{\mathbf{M}\in\RR_{>0}^{\NN}:\;\mathbf{M}\;\text{is normalized, log-convex},\;\lim_{p\rightarrow+\infty}(M_p)^{1/p}=+\infty\}.$$
We see that $\mathbf{M}\in\hyperlink{LCset}{\mathcal{LC}}$ if and only if $1=\mu_0\le\mu_1\le\dots$ and $\lim_{p\rightarrow+\infty}\mu_p=+\infty$ (see e.g. \cite[p. 104]{compositionpaper}) and there is a one-to-one correspondence between $\mathbf{M}$ and $\mu=(\mu_p)_p$ by taking $M_p:=\prod_{i=0}^p\mu_i$.\vspace{6pt}

Let $\mathbf{M},\mathbf{N}\in\RR_{>0}^{\NN}$ be given, then  write
\begin{itemize}
\item[$(*)$] $\mathbf{M}\le\mathbf{N}$ if $M_p\le N_p$ for all $p\in\NN$,

\item[$(*)$] $\mathbf{M}\hypertarget{preceq}{\preceq}\mathbf{N}$ if
$$\sup_{p\in\NN_{>0}}\left(\frac{M_p}{N_p}\right)^{1/p}<+\infty,$$

\item[$(*)$] $\mathbf{M}\hypertarget{mtriangle}{\vartriangleleft}\mathbf{N}$ if $\lim_{p\rightarrow+\infty}\left(\frac{M_p}{N_p}\right)^{1/p}=0$; i.e. if
\begin{equation*}\label{triangleestim}
\forall\;h>0\;\exists\;C_h\ge 1\;\forall\;p\in\NN:\;\;\;M_p\le C_hh^pN_p.
\end{equation*}
\item[$(*)$] $\mathbf{M}$ and $\mathbf{N}$ are called \emph{equivalent,} formally denoted by $\mathbf{M}\hypertarget{approx}{\approx}\mathbf{N}$, if $\mathbf{M}\hyperlink{preceq}{\preceq}\mathbf{N}$ and $\mathbf{N}\hyperlink{preceq}{\preceq}\mathbf{M}$.

\item[$(*)$] Obviously, $\mathbf{M}\hyperlink{mtriangle}{\vartriangleleft}\mathbf{N}$ implies $\mathbf{M}\hyperlink{preceq}{\preceq}\mathbf{N}$, but \hyperlink{mtriangle}{$\vartriangleleft$} neither is reflexive nor symmetric.
\end{itemize}

$\mathbf{M}$ satisfies \emph{moderate growth,} abbreviated by \hypertarget{mg}{$(\on{mg})$}, if
$$\exists\;C\ge 1\;\forall\;p,q\in\NN:\;M_{p+q}\le C^{p+q+1} M_p M_q.$$
In view of \cite{Komatsu73} this condition is also known under \emph{(M.2) or stability under ultradifferential operators} and by definition moderate growth is obviously preserved under equivalence. In case $M_0=1$, then it suffices to consider $C^{p+q}$ in the estimate.\vspace{6pt}

An important example, also for the considerations in this work, are the \emph{Gevrey sequences} $\mathbf{G}^s:=(p!^s)_{p\in\NN}$, $s>0$. Note that $\mathbf{G}^s$ is equivalent to $\overline{\mathbf{G}}^s$ with $\overline{\mathbf{G}}^s:=(p^{ps})_{p\in\NN}$ by \emph{Stirling's formula.} Each $\mathbf{G}^s$ satisfies \hyperlink{mg}{$(\on{mg})$}.

\subsection{Associated weight functions}\label{assofunctionsection}
For the following definition we refer to \cite[Chapitre I]{mandelbrojtbook} and the more recent work \cite{regularnew} where more properties are summarized. Let $\mathbf{M}=(M_p)_{p\in\NN}\in\RR_{>0}^{\NN}$ be given with $M_0=1$, then the \emph{associated function} $\omega_{\mathbf{M}}: [0,+\infty)\rightarrow[0,+\infty)\cup\{+\infty\}$ is defined as follows:
\begin{equation}\label{assofunc}
\omega_{\mathbf{M}}(t):=\sup_{p\in\NN}\log\left(\frac{t^p}{M_p}\right),\qquad t\ge 0,
\end{equation}
with the convention that $0^0:=1$. This ensures $\omega_{\mathbf{M}}(0)=0$ and $\omega_{\mathbf{M}}(t)\ge 0$ for any $t\ge 0$ since $\frac{t^0M_0}{M_0}=1$ for all $t\ge 0$. \eqref{assofunc} corresponds to \cite[$(3.1)$]{Komatsu73} and we immediately have that $\omega_{\mathbf{M}}$ is non-decreasing and satisfying $\lim_{t\rightarrow+\infty}\omega_{\mathbf{M}}(t)=+\infty$.

\subsection{Weight matrices}\label{matrixsection}
Let $\mathcal{I}=\RR_{>0}$ be the index set. A \emph{weight matrix} $\mathcal{M}$ associated with $\mathcal{I}$ is the set $\mathcal{M}=\{\mathbf{M}^{(x)}: x>0\}$ such that $\mathbf{M}^{(x)}\le\mathbf{M}^{(y)}$ for all $0<x\le y$. $\mathcal{M}$ is called
\begin{itemize}
\item[$(*)$]\emph{constant} if $\mathbf{M}^{(x)}\hyperlink{approx}{\approx}\mathbf{M}^{(y)}$ for all $x,y>0$,

\item[$(*)$] \emph{log-convex} if $\mathbf{M}^{(x)}$ is log-convex for all $x>0$,

\item[$(*)$] \emph{standard log-convex} if $\mathbf{M}^{(x)}\in\hyperlink{LCset}{\mathcal{LC}}$ for all $x>0$.
\end{itemize}

Let us consider some growth conditions on $\mathcal{M}$ and relations between matrices $\mathcal{M}$ and $\mathcal{N}$; see \cite[Sect. 4.1 \& 4.2]{compositionpaper}:

\hypertarget{R-mg}{$(\mathcal{M}_{\{\on{mg}\}})$} \hskip1cm $\forall\;\alpha\in\mathcal{I}\;\exists\;C>0\;\exists\;\beta\in\mathcal{I}\;\forall\;p,q\in\NN:\;\;\;M^{(\alpha)}_{p+q}\le C^{p+q+1} M^{(\beta)}_p M^{(\beta)}_q$,\par\vskip.3cm
\hypertarget{B-mg}{$(\mathcal{M}_{(\on{mg})})$} \hskip1cm $\forall\;\alpha\in\mathcal{I}\;\exists\;C>0\;\exists\;\beta\in\mathcal{I}\;\forall\;p,q\in\NN:\;\;\;M^{(\beta)}_{p+q}\le C^{p+q+1} M^{(\alpha)}_p M^{(\alpha)}_q$,\par\vskip.3cm
\hypertarget{R-L}{$(\mathcal{M}_{\{\on{L}\}})$} \hskip1cm $\forall\;\alpha\in\mathcal{I}\;\forall\;h>0\;\exists\;\beta\in\mathcal{I}\;\exists\;C>0\;\forall\;p\in\NN:\;\;\;h^pM^{(\alpha)}_{p}\le CM^{(\beta)}_p$,\par\vskip.3cm
\hypertarget{B-L}{$(\mathcal{M}_{(\on{L})})$} \hskip1cm $\forall\;\alpha\in\mathcal{I}\;\forall\;h>0\;\exists\;\beta\in\mathcal{I}\;\exists\;C>0\;\forall\;p\in\NN:\;\;\;h^pM^{(\beta)}_{p}\le CM^{(\alpha)}_p$,\par\vskip.3cm

Concerning the growth relations between weight matrices write:
$$\mathcal{M}\{\preceq\}\mathcal{N}:\;\;\;\forall\;\alpha>0\;\exists\;\beta>0:\;\;\;\mathbf{M}^{(\alpha)}\hyperlink{preceq}{\preceq}\mathbf{N}^{(\beta)},$$
$$\mathcal{M}(\preceq)\mathcal{N}:\;\;\;\forall\;\alpha>0\;\exists\;\beta>0:\;\;\;\mathbf{M}^{(\beta)}\hyperlink{preceq}{\preceq}\mathbf{N}^{(\alpha)},$$
$$\mathcal{M}\vartriangleleft\mathcal{N}:\;\;\;\forall\;\alpha,\beta>0:\;\;\;\mathbf{M}^{(\alpha)}\hyperlink{mtriangle}{\vartriangleleft}\mathbf{N}^{(\beta)}.$$
We call $\mathcal{M}$ and $\mathcal{N}$ \emph{R-equivalent,} written $\mathcal{M}\{\approx\}\mathcal{N}$, if $\mathcal{M}\{\preceq\}\mathcal{N}$ and $\mathcal{N}\{\preceq\}\mathcal{M}$ and \emph{B-equivalent,} denoted by $\mathcal{M}(\approx)\mathcal{N}$, if $\mathcal{M}(\preceq)\mathcal{N}$ and $\mathcal{N}(\preceq)\mathcal{M}$. Relation $\vartriangleleft$ is neither reflexive nor symmetric.\vspace{6pt}

\emph{Convention:} We write the symbol $[\cdot]$ as a uniform notation if mean either $\{\cdot\}$ or $(\cdot)$ but not mixing the cases. A similar convention is used for the corresponding weighted spaces. Recall that in the literature $\{\cdot\}$ refers to the \emph{Roumieu-type spaces} and $(\cdot)$ to the \emph{Beurling-type.}

\subsection{Weight functions and the growth indices $\gamma(\omega)$ and $\overline{\gamma}(\omega)$}\label{growthindexsect}
We briefly recall the definitions of the growth indices $\gamma(\omega)$ and $\overline{\gamma}(\omega)$; see \cite[Sect. 2.3 $(6)$, Sect. 2.4 $(7)$]{index} and the references in these sections. First, recall the notion of being a weight function in the sense of \cite{index} and \cite{genLegendreconj}:

\begin{definition}\label{weightfctdef}
$\omega:[0,+\infty)\rightarrow[0,+\infty)$ is called a \emph{weight function} if
\begin{itemize}
\item[$(*)$] $\omega$ is non-decreasing and

\item[$(*)$] $\lim_{t\rightarrow+\infty}\omega(t)=+\infty$.
\end{itemize}
\end{definition}
We are going to use the following notation for any (weight) function $\omega:[0,+\infty)\rightarrow[0,+\infty)$ and arbitrary $\alpha>0$: set $\omega^{\iota}(t):=\omega(\frac{1}{t})$, $t>0$, and $\omega^{1/\alpha}(t):=\omega(t^{1/\alpha})$, $t\ge 0$. Note that $\omega^{1/\alpha}$ is obtained by a power-substitution and it is again a weight function. Moreover, $(\omega^{\iota})^{1/\alpha}(t)=(\omega^{1/\alpha})^{\iota}(t)$ for all $t>0$ and finally write $\id^{1/\alpha}$ for $t\mapsto t^{1/\alpha}$.\vspace{6pt}

Weight functions in the sense of Definition \ref{weightfctdef} are sufficient in order to introduce and work with the crucial indices $\gamma(\omega)$ and $\overline{\gamma}(\omega)$. Let $\omega$ be a weight function and $\gamma>0$. We say that $\omega$ has property $(P_{\omega,\gamma})$ if
\begin{equation}\label{newindex1}
\exists\;K>1:\;\;\;\limsup_{t\rightarrow+\infty}\frac{\omega(K^{\gamma}t)}{\omega(t)}<K.
\end{equation}
If $(P_{\omega,\gamma})$ holds for some $K>1$, then $(P_{\omega,\gamma'})$ is satisfied for all $\gamma'\le\gamma$ with the same $K$ since $\omega$ is non-decreasing. Moreover we can restrict to $\gamma>0$, because for $\gamma\le 0$ condition $(P_{\omega,\gamma})$ is satisfied for any weight function $\omega$, again because $\omega$ is non-decreasing and $K>1$. Then put
\begin{equation}\label{newindex2}
\gamma(\omega):=\sup\{\gamma>0:\;\;(P_{\omega,\gamma})\;\;\text{is satisfied}\},
\end{equation}
and if none condition $(P_{\omega,\gamma})$ (with $\gamma>0$) holds then set $\gamma(\omega):=0$.

Analogously, for $\gamma>0$ we say that $\omega$ has property $(\overline{P}_{\omega,\gamma})$ if
\begin{equation}\label{newindex3}
\exists\;A>1:\;\;\;\liminf_{t\rightarrow+\infty}\frac{\omega(A^{\gamma}t)}{\omega(t)}>A.
\end{equation}
If $(\overline{P}_{\omega,\gamma})$ holds for some $A>1$, then $(\overline{P}_{\omega,\gamma'})$ is satisfied for all $\gamma'\ge\gamma$ with the same $A$ since $\omega$ is non-decreasing. Moreover, we can restrict to $\gamma>0$ because for $\gamma\le 0$ condition $(\overline{P}_{\omega,\gamma})$ is never satisfied for any weight function ($\omega$ is assumed to be non-decreasing and $A>1$). Then set
\begin{equation}\label{newindex4}
\overline{\gamma}(\omega):=\inf\{\gamma>0: \;\;(\overline{P}_{\omega,\gamma})\;\;\text{is satisfied}\}.
\end{equation}
We obtain $(0\le)\gamma(\omega)\le\overline{\gamma}(\omega)$, see \cite[Sect. 2.3 \& 2.4]{index}, and by definition the following identities are valid:
\begin{equation}\label{newindex5}
\forall\;\alpha>0:\;\;\;\gamma(\omega^{1/\alpha})=\alpha\gamma(\omega),\hspace{15pt}\overline{\gamma}(\omega^{1/\alpha})=\alpha\overline{\gamma}(\omega).
\end{equation}

Finally, let us recall some crucial growth relations for weight functions: We write
\begin{itemize}
\item[$(*)$] $\sigma\hypertarget{ompreceq}{\preceq}\tau$ if
\begin{equation}\label{bigOrelation}
\tau(t)=O(\sigma(t))\;\text{as}\;t\rightarrow+\infty,	
\end{equation}
\item[$(*)$] $\sigma\hypertarget{omvartriangle}{\vartriangleleft}\tau$ if
\begin{equation}\label{smallOrelation}
	\tau(t)=o(\sigma(t))\;\text{as}\;t\rightarrow+\infty.	
\end{equation}
\item[$(*)$] $\sigma$ and $\tau$ are called \emph{equivalent,} written $\sigma\hypertarget{sim}{\sim}\tau$, if $\sigma\hyperlink{ompreceq}{\preceq}\tau$ and $\tau\hyperlink{ompreceq}{\preceq}\sigma$.

\item[$(*)$] Note that relation \hyperlink{omvartriangle}{$\vartriangleleft$} neither is reflexive nor symmetric.
\end{itemize}

The indices $\gamma(\cdot)$ and $\overline{\gamma}(\cdot)$ are preserved under \hyperlink{sim}{$\sim$}; this follows from the characterizations in the main results \cite[Thm. 2.11 \& 2.16]{index}, see also \cite[Rem. 2.12]{index} and the comments after \cite[Thm. 2.16]{index}.

\subsection{Generalized lower and upper Legendre conjugates}\label{generalizedLegendresection}
Let $\sigma,\tau$ be weight functions according to Definition \ref{weightfctdef}. We introduce the \emph{generalized lower Legendre conjugate} by
\begin{equation}\label{wedgeformula}
\sigma\check{\star}\tau(t):=\inf_{s>0}\{\sigma(s)+\tau(t/s)\},\;\;\;t\in[0,+\infty),
\end{equation}
and the \emph{generalized upper Legendre conjugate} by
\begin{equation}\label{widehatformula}
\sigma\widehat{\star}\tau(t):=\sup_{s\ge 0}\{\sigma(s)-\tau(s/t)\},\;\;\;t\in(0,+\infty).
\end{equation}
We also set $\sigma\widehat{\star}\tau(0):=\sigma(0)-\tau(0)$ which is justified by \cite[Lemma 4.1 $(b)$]{genLegendreconj}. These operations and their effects on the growth indices $\gamma(\cdot)$, $\overline{\gamma}(\cdot)$ are studied in detail in \cite[Sect. 3 \& 4]{genLegendreconj} for (general) weight functions and in \cite[Sect. 5]{genLegendreconj} when $\sigma=\omega_{\mathbf{M}}$, $\tau=\omega_{\mathbf{N}}$.

$\check{\star}$ yields again a weight function but, in general, for $\widehat{\star}$ this is not clear; see requirements $(A)$ and $(B)$ in \cite[Sect. 4.2]{genLegendreconj}: First, in order to ensure $\sigma\widehat{\star}\tau(t)\ge 0$ for all $t$ it suffices to assume that $\tau(0)=0$, see \cite[Lemma 4.2]{genLegendreconj}, and second, in order to guarantee that $\sigma\widehat{\star}\tau$ is well-defined, i.e. that $\sigma\widehat{\star}\tau(t)<+\infty$ for all $t$, it is required to have
\begin{equation}\label{equ39}
\forall\;t\in(0,+\infty)\;\exists\;D_t>0\;\forall\;s\ge 0:\;\;\;\sigma(s)-\tau(s/t)\le D_t.
\end{equation}
This condition corresponds to $(B)$ in \cite[Sect. 4.2]{genLegendreconj} resp. to \cite[$(4.4)$]{genLegendreconj} with $t_0=+\infty$.\vspace{6pt}

Special but important cases are provided by considering $\tau:=\id^{1/\alpha}$, $\alpha>0$, which gives
\begin{equation}\label{lowerenvelope}
\sigma\check{\star}\id(t)=\inf_{s>0}\{\sigma(s)+t/s\}=\inf_{u>0}\{\sigma(1/u)+tu\}=:(\sigma^{\iota})_{\star}(t),
\end{equation}
$h_{\star}(t):=\inf_{u>0}\{h(u)+tu\}$ denoting the known \emph{lower Legendre conjugate (or envelope),} see \cite[Sect. 2.5]{index}, and
\begin{equation}\label{lowerLegendregeneral}
\forall\;\alpha>0\;\forall\;t\ge 0:\;\;\;\sigma\check{\star}\id^{1/\alpha}(t)=(((\sigma^{\iota})^{\alpha})_{\star})^{1/\alpha}(t);
\end{equation}
see \cite[$(3.3)$ \& $(3.4)$]{genLegendreconj}. Moreover, recall \cite[$(4.9)$ \& $(4.10)$]{genLegendreconj}:
\begin{equation}\label{upperLegendre}
\sigma\widehat{\star}\id(t)=\sup_{s\ge 0}\{\sigma(s)-s/t\}=:\sigma^{\star}(1/t)=(\sigma^{\star})^{\iota}(t),\;\;\;t\in(0,+\infty),
\end{equation}
with $\sigma^{\star}(t):=\sup_{s\ge 0}\{\sigma(s)-ts\}$ denoting the known \emph{upper Legendre conjugate (or envelope),} see again \cite[Sect. 2.5]{index}, and
\begin{equation}\label{upperLegendregeneral}
\forall\;\alpha>0\;\forall\;t\ge 0:\;\;\;\sigma\widehat{\star}\id^{1/\alpha}(t)=(((\sigma^{\alpha})^{\star})^{\iota})^{1/\alpha}(t).
\end{equation}

\subsection{Weight functions in the sense of Braun-Meise-Taylor}\label{BMTdefsection}
We consider now a crucial subclass of weight functions; i.e. functions in the sense of Braun-Meise-Taylor (BMT-weights for short). Let $\omega:[0,+\infty)\rightarrow[0,+\infty)$ be continuous, non-decreasing, $\omega(0)=0$ and $\lim_{t\rightarrow+\infty}\omega(t)=+\infty$. If $\omega$ satisfies in addition $\omega(t)=0$ for all $t\in[0,1]$, then $\omega$ is called \emph{normalized.} For convenience we write that $\omega$ has $\hypertarget{om0}{(\omega_0)}$ if it satisfies all these assumptions; see e.g. \cite{BraunMeiseTaylor90}, \cite[Sect. 2.1]{index} and \cite[Sect. 2.2]{sectorialextensions}.\vspace{6pt}

Moreover we consider the following conditions; this list of properties has already been used in ~\cite{dissertation}.

\begin{itemize}
\item[\hypertarget{om1}{$(\omega_1)}$] $\omega(2t)=O(\omega(t))$ as $t\rightarrow+\infty$, i.e. $\exists\;L\ge 1\;\forall\;t\ge 0:\;\;\;\omega(2t)\le L(\omega(t)+1)$.

%\item[\hypertarget{om2}{$(\omega_2)$}] $\omega(t)=O(t)$ as $t\rightarrow+\infty$.

\item[\hypertarget{om3}{$(\omega_3)$}] $\log(t)=o(\omega(t))$ as $t\rightarrow+\infty$.

\item[\hypertarget{om4}{$(\omega_4)$}] $\varphi_{\omega}:t\mapsto\omega(e^t)$ is a convex function on $\RR$.

\item[\hypertarget{om5}{$(\omega_5)$}] $\omega(t)=o(t)$ as $t\rightarrow+\infty$.

\item[\hypertarget{om6}{$(\omega_6)$}] $\exists\;H\ge 1\;\forall\;t\ge 0:\;2\omega(t)\le\omega(H t)+H$.
\end{itemize}

Finally, we recall the \emph{strong nonquasianalyticity} condition for weight functions
\begin{equation}\label{assostrongnq}
\exists\;C\ge 1\;\forall\;y\ge 0:\;\;\;\int_1^{+\infty}\frac{\omega(yt)}{t^2}dt\le C\omega(y)+C.
\end{equation}

In the literature one can find different assumptions for BMT-weights. However, some of them are basic and for convenience we introduce
$$\hypertarget{omset0}{\mathcal{W}_0}:=\{\omega:[0,\infty)\rightarrow[0,\infty): \omega\;\text{has}\;\hyperlink{om0}{(\omega_0)},\hyperlink{om3}{(\omega_3)},\hyperlink{om4}{(\omega_4)}\}.$$
In the forthcoming \hyperlink{omset0}{$\mathcal{W}_0$} will be understood to be the set of (normalized) weight functions in the sense of Braun-Meise-Taylor. Then, for any $\omega\in\hyperlink{omset0}{\mathcal{W}_0}$ we define the \emph{Legendre-Fenchel-Young-conjugate} of $\varphi_{\omega}$ by
\begin{equation}\label{legendreconjugate}
\varphi^{*}_{\omega}(x):=\sup\{x y-\varphi_{\omega}(y): y\ge 0\},\;\;\;x\ge 0,
\end{equation}
with the following properties, see e.g. \cite[Rem. 1.3, Lemma 1.5]{BraunMeiseTaylor90}: It is convex and non-decreasing, $\varphi^{*}_{\omega}(0)=0$, $\varphi^{**}_{\omega}=\varphi_{\omega}$, $\lim_{x\rightarrow+\infty}\frac{x}{\varphi^{*}_{\omega}(x)}=0$ and finally $x\mapsto\frac{\varphi_{\omega}(x)}{x}$ and $x\mapsto\frac{\varphi^{*}_{\omega}(x)}{x}$ are non-decreasing on $[0,+\infty)$. Note that by normalization we can extend the supremum in \eqref{legendreconjugate} from $y\ge 0$ to $y\in\RR$ without changing the value of $\varphi^{*}_{\omega}(x)$ for given $x\ge 0$.\vspace{6pt}

We mention the following known result, see e.g. \cite[Lemma 2.8]{testfunctioncharacterization} resp. \cite[Lemma 2.4]{sectorialextensions} and the references mentioned in the proofs there.

\begin{lemma}\label{assoweightomega0}
Let $\mathbf{M}\in\hyperlink{LCset}{\mathcal{LC}}$, then $\omega_{\mathbf{M}}\in\hyperlink{omset0}{\mathcal{W}_0}$ holds and \hyperlink{om6}{$(\omega_6)$} for $\omega_{\mathbf{M}}$ if and only if $\mathbf{M}$ has \hyperlink{mg}{$(\on{mg})$}.
\end{lemma}

Finally, we recall that the growth indices from Section \ref{growthindexsect} can equivalently be expressed by standard growth conditions for BMT-weight functions:

\begin{remark}\label{om1om6indexrem}
\emph{Let $\omega$ be a weight function. Then $\omega$ has \hyperlink{om1}{$(\omega_1)$} if and only if $\gamma(\omega)>0$, see \cite[Thm. 2.11, Rem. 2.12, Cor. 2.14]{index}, and by \cite[Thm. 2.16, $(7)$, Cor. 2.17]{index} condition \hyperlink{om6}{$(\omega_6)$} holds if and only if $\overline{\gamma}(\omega)<+\infty$.}
\end{remark}

\subsection{Associated weight matrix}\label{assomatrixsection}
We summarize some facts which are shown in \cite[Sect. 5]{compositionpaper} and are needed in this work; all properties listed below are valid for $\omega\in\hyperlink{omset0}{\mathcal{W}_0}$ except \eqref{newexpabsorb} for which \hyperlink{om1}{$(\omega_1)$} is crucial.

\begin{itemize}
\item[$(i)$] The idea was that to each $\omega\in\hyperlink{omset0}{\mathcal{W}_0}$ one can associate a weight matrix $\mathcal{M}_{\omega}:=\{\mathbf{W}^{(\ell)}=(W^{(\ell)}_p)_{p\in\NN}: \ell>0\}$ by\vspace{6pt}

    \centerline{$W^{(\ell)}_p:=\exp\left(\frac{1}{\ell}\varphi^{*}_{\omega}(\ell p)\right)$.}\vspace{6pt}

We have that $\mathbf{W}^{(\ell)}\in\hyperlink{LCset}{\mathcal{LC}}$ for each $\ell>0$, i.e. $\mathcal{M}_{\omega}$ is \emph{standard log-convex,} and we even have the order relation
\begin{equation}\label{quotientorderequ}
\forall\;\ell_2\ge\ell_1>0:\;\;\;\vartheta^{(\ell_1)}\le\vartheta^{(\ell_2)},
\end{equation}
with $\vartheta^{(\ell)}$ denoting the corresponding sequence of quotients; see \cite[Sect. 2.5]{whitneyextensionweightmatrix}.

\item[$(ii)$] $\mathcal{M}_{\omega}$ satisfies
    \begin{equation}\label{newmoderategrowth}
    \forall\;\ell>0\;\forall\;p,q\in\NN:\;\;\;W^{(\ell)}_{p+q}\le W^{(2\ell)}_pW^{(2\ell)}_q.
    \end{equation}
    so both \hyperlink{R-mg}{$(\mathcal{M}_{\{\on{mg}\}})$} and \hyperlink{B-mg}{$(\mathcal{M}_{(\on{mg})})$} are valid.

\item[$(iii)$] \hyperlink{om6}{$(\omega_6)$} holds if and only if some/each $\mathbf{W}^{(\ell)}$ satisfies \hyperlink{mg}{$(\on{mg})$} if and only if $\mathbf{W}^{(\ell)}\hyperlink{approx}{\approx}\mathbf{W}^{(n)}$ for each $\ell,n>0$. Thus \hyperlink{om6}{$(\omega_6)$} is characterizing the case when $\mathcal{M}_{\omega}$ is \emph{constant.}

\item[$(iv)$] In case $\omega$ has in addition \hyperlink{om1}{$(\omega_1)$}, then $\mathcal{M}_{\omega}$ also satisfies
     \begin{equation}\label{newexpabsorb}
     \forall\;h\ge 1\;\exists\;d\ge 1\;\forall\;\ell>0\;\exists\;D\ge 1\;\forall\;p\in\NN:\;\;\;h^jW^{(\ell)}_p\le D W^{(d\ell)}_p,
     \end{equation}
hence both \hyperlink{R-L}{$(\mathcal{M}_{\{\on{L}\}})$} and \hyperlink{B-L}{$(\mathcal{M}_{(\on{L})})$}. This estimate is crucial for verifying $\mathcal{E}_{[\mathcal{M}_{\omega}]}=\mathcal{E}_{[\omega]}$ (as locally convex vector spaces).

\item[$(v)$] We have $\omega\hyperlink{sim}{\sim}\omega_{\mathbf{W}^{(\ell)}}$ for each $\ell>0$, more precisely
\begin{equation}\label{goodequivalenceclassic}
\forall\;\ell>0\,\,\exists\,D_{\ell}>0\;\forall\;t\ge 0:\;\;\;\ell\omega_{\mathbf{W}^{(\ell)}}(t)\le\omega(t)\le 2\ell\omega_{\mathbf{W}^{(\ell)}}(t)+D_{\ell};
\end{equation}
for a proof see \cite[Theorem 4.0.3, Lemma 5.1.3]{dissertation}, \cite[Lemma 5.7]{compositionpaper} and also \cite[Lemma 2.5]{sectorialextensions}. Note that, on the one hand, for proving \eqref{goodequivalenceclassic} the convexity condition \hyperlink{om4}{$(\omega_4)$} is indispensable but, on the other hand, \hyperlink{om4}{$(\omega_4)$} is only required for the second estimate in \eqref{goodequivalenceclassic}.
\end{itemize}

Additionally, by inspecting the proofs in \cite{compositionpaper} (and \cite{dissertation}) it turns out:

\begin{itemize}
\item[$(a)$] For all listed statements above, except the second estimate in \eqref{goodequivalenceclassic} and for the ``only if'' in $(iii)$, the convexity condition \hyperlink{om4}{$(\omega_4)$} is not required necessarily: For transferring properties from $\omega$ to $\mathcal{M}_{\omega}$ only the (automatically) convexity for $\varphi^{*}_{\omega}$ is crucial and \hyperlink{om4}{$(\omega_4)$} is required for showing characterizations; see also the proofs in \cite[Sect. 4 \& 5]{dissertation} and the comments in \cite[Sect. 3.1]{dissertation}.

\item[$(b)$] Let $\omega\in\hyperlink{omset0}{\mathcal{W}_0}$ be given and $\sigma$ be a weight function in the sense of Definition \ref{weightfctdef}. Assume that $\sigma\hyperlink{sim}{\sim}\omega$ holds, then automatically \hyperlink{om3}{$(\omega_3)$} is valid for $\sigma$ and by repeating the proof of \cite[Lemma 5.16]{compositionpaper} we infer that $\mathcal{M}_{\omega}\{\approx\}\mathcal{M}_{\sigma}$ and $\mathcal{M}_{\omega}(\approx)\mathcal{M}_{\sigma}$: For more details see the comments in \cite[Sect. 6 $(III)$]{modgrowthstrangeII}, in particular formula \cite[$(6.6)$]{modgrowthstrangeII}, and also Lemma \ref{om1om6lemma1}. Here, the sequences $\mathbf{S}^{(\ell)}\in\mathcal{M}_{\sigma}$ are defined analogously as in $(i)$ above when involving $\varphi_{\sigma}^{*}$. Each $\mathbf{S}^{(\ell)}$ satisfies all properties of the set \hyperlink{LCset}{$\mathcal{LC}$} except normalization in general. Therefore, $\mathcal{M}_{\sigma}$ is formally not standard log-convex but this is only a minor technical issue.

\item[$(c)$] Now, in the situation of comment $(b)$, when $\omega$ satisfies in addition \hyperlink{om1}{$(\omega_1)$}, then $\sigma$ too and (as l.c.v.s.) one has $$\mathfrak{F}_{[\omega]}=\mathfrak{F}_{[\mathcal{M}_{\omega}]}=\mathfrak{F}_{[\mathcal{M}_{\sigma}]}=\mathfrak{F}_{[\sigma]},$$
    where $\mathfrak{F}\in\{\mathcal{E}, \mathcal{B}, \mathcal{A}, \mathcal{S}, \Lambda, \mathcal{F}\}$ and the last space above is defined analogously by involving $\varphi_{\sigma}^{*}$. (Note that when $\sigma$ is not normalized, then the natural definition for $\varphi^{*}_{\sigma}$ is to consider the supremum over all $y\in\RR$ in \eqref{legendreconjugate}; see again \cite[Sect. 6 $(III)$]{modgrowthstrangeII}.)
\end{itemize}

\section{Generalized lower Legendre conjugate}\label{lowersection}
We are investigating now the effects of the conjugate $\check{\star}$ for BMT-weight functions and their corresponding associated weight matrices.\vspace{6pt}

First, let us illustrate how the growth indices $\gamma(\cdot)$ and $\overline{\gamma}(\cdot)$ are preserved within the associated weight matrix. Of course, the main results \cite[Thm. 3.4, Cor. 3.5]{genLegendreconj} and the comments from \cite[Rem. 3.6]{genLegendreconj} can be applied to any $\sigma,\tau\in\hyperlink{omset0}{\mathcal{W}_0}$ and to any of the corresponding associated weight functions $\omega_{\mathbf{S}^{(\ell)}}$, $\omega_{\mathbf{T}^{(j)}}$ separately. However, the next result shows the flexibility w.r.t. the matrix parameter(s).

\begin{theorem}\label{lowertransformindexthm1}
Let $\sigma,\tau\in\hyperlink{omset0}{\mathcal{W}_0}$ be given with associated weight matrices $\mathcal{M}_{\sigma}:=\{\mathbf{S}^{(\ell)}: \ell>0\}$, $\mathcal{M}_{\tau}:=\{\mathbf{T}^{(\ell)}: \ell>0\}$.
\begin{itemize}
\item[$(i)$] Assume that $\gamma(\sigma),\gamma(\tau)>0$, then
$$\forall\;\ell,\ell_1,j,j_1>0:\;\;\;\gamma(\omega_{\mathbf{S}^{(\ell)}})+\gamma(\omega_{\mathbf{T}^{(j)}})=\gamma(\sigma)+\gamma(\tau)\le\gamma(\sigma\check{\star}\tau)=\gamma(\omega_{\mathbf{S}^{(\ell_1)}}\check{\star}\omega_{\mathbf{T}^{(j_1)}}).$$

\item[$(ii)$] Assume that $\overline{\gamma}(\sigma),\overline{\gamma}(\tau)<+\infty$, then
$$\forall\;\ell,\ell_1,j,j_1>0:\;\;\;\overline{\gamma}(\omega_{\mathbf{S}^{(\ell)}}\check{\star}\omega_{\mathbf{T}^{(j)}})=\overline{\gamma}(\sigma\check{\star}\tau)\le\overline{\gamma}(\sigma)+\overline{\gamma}(\tau)=\overline{\gamma}(\omega_{\mathbf{S}^{(\ell_1)}})+\overline{\gamma}(\omega_{\mathbf{T}^{(j_1)}}).$$

\item[$(iii)$] When $\gamma(\sigma)>0$, then
$$\forall\;\ell,\ell_1,\alpha>0:\;\;\;\gamma(\omega_{\mathbf{S}^{(\ell)}})+\alpha\le\gamma((((\omega_{\mathbf{S}^{(\ell_1)}}^{\iota})^{\alpha})_{\star})^{1/\alpha}),$$
and $\overline{\gamma}(\sigma)<+\infty$ implies
$$\forall\;\ell,\ell_1,\alpha>0:\;\;\;\overline{\gamma}((((\omega_{\mathbf{S}^{(\ell)}}^{\iota})^{\alpha})_{\star})^{1/\alpha})\le\overline{\gamma}(\omega_{\mathbf{S}^{(\ell_1)}})+\alpha.$$

\item[$(iv)$] When $0<\gamma(\sigma)=\overline{\gamma}(\sigma)<+\infty$, $0<\gamma(\tau)=\overline{\gamma}(\tau)<+\infty$, then
$$\forall\;\ell,\ell_1,\ell_2,j,j_1,j_2>0:\;\;\;\gamma(\omega_{\mathbf{S}^{(\ell)}})+\gamma(\omega_{\mathbf{T}^{(j)}})=\gamma(\omega_{\mathbf{S}^{(\ell_1)}}\check{\star}\omega_{\mathbf{T}^{(j_1)}})=\overline{\gamma}(\omega_{\mathbf{S}^{(\ell_2)}}\check{\star}\omega_{\mathbf{T}^{(j_2)}}).$$
\end{itemize}
\end{theorem}

\demo{Proof}
The conclusions follow by taking into account \eqref{goodequivalenceclassic}, the fact that the indices $\gamma(\cdot)$ and $\overline{\gamma}(\cdot)$ are preserved under relation \hyperlink{sim}{$\sim$}, and by using \cite[Lemma 3.3, Thm. 3.4, Cor. 3.5]{genLegendreconj}. Finally, recall formula \cite[$(3.4)$]{genLegendreconj}; see also \eqref{lowerLegendregeneral}.
\qed\enddemo

\emph{Note:}

\begin{itemize}
\item[$(*)$] In view of Remark \ref{om1om6indexrem}, in $(i)$ we assume that both $\sigma$ and $\tau$ have \hyperlink{om1}{$(\omega_1)$} and \hyperlink{om6}{$(\omega_6)$} in $(ii)$.

\item[$(*)$] However, even without these assumptions on the indices the corresponding implications are valid by taking into account the comments from \cite[Rem. 3.6]{genLegendreconj} and the fact that \eqref{goodequivalenceclassic} preserves both indices and that for this property neither \hyperlink{om1}{$(\omega_1)$} nor \hyperlink{om6}{$(\omega_6)$} is involved; the defining growth conditions for the set \hyperlink{omset0}{$\mathcal{W}_0$} are sufficient to conclude.
\end{itemize}

Let $\sigma,\tau\in\hyperlink{omset0}{\mathcal{W}_0}$ be given with associated matrices $\mathcal{M}_{\sigma}$, $\mathcal{M}_{\tau}$ and define
\begin{equation}\label{productmatrix}
\mathcal{M}_{\sigma}\cdot\mathcal{M}_{\tau}:=\{\mathbf{S}^{(\ell)}\cdot\mathbf{T}^{(\ell)}: \ell>0\}.
\end{equation}
Moreover, as a special case for given $\mathbf{N}\in\RR_{>0}^{\NN}$ set
\begin{equation}\label{productmatrixspecial}
\mathcal{M}_{\sigma}\cdot\mathbf{N}:=\{\mathbf{S}^{(\ell)}\cdot\mathbf{N}: \ell>0\}.
\end{equation}

Using this new notation we are going to show the main result of this section:

\begin{theorem}\label{mainweighfctthm}
Let $\sigma,\tau\in\hyperlink{omset0}{\mathcal{W}_0}$ be given with associated matrices $\mathcal{M}_{\sigma}$, $\mathcal{M}_{\tau}$.

\begin{itemize}
\item[$(i)$] It holds that
     \begin{equation}\label{mainweighfctthmequ0}
    \forall\;\ell,\ell_1>0:\;\;\;\omega_{\mathbf{S}^{(\ell)}\cdot\mathbf{T}^{(\ell_1)}}=\omega_{\mathbf{S}^{(\ell)}}\check{\star}\omega_{\mathbf{T}^{(\ell_1)}}\hyperlink{sim}{\sim}\sigma\check{\star}\tau,
    \end{equation}
    which implies, in particular, that
    \begin{equation}\label{mainweighfctthmequ}
    \forall\;\ell,j>0:\;\;\;\omega_{\mathbf{S}^{(\ell)}\cdot\mathbf{T}^{(\ell)}}\hyperlink{sim}{\sim}\omega_{\mathbf{S}^{(j)}\cdot\mathbf{T}^{(j)}}.
\end{equation}

\item[$(ii)$] $\mathcal{M}_{\sigma}\cdot\mathcal{M}_{\tau}$ satisfies \eqref{newmoderategrowth}.

\item[$(iii)$] If either $\sigma$ or $\tau$ satisfies in addition \hyperlink{om1}{$(\omega_1)$}, then $\mathcal{M}_{\sigma}\cdot\mathcal{M}_{\tau}$ has \eqref{newexpabsorb}, too.

\item[$(iv)$] If either $\sigma$ or $\tau$ satisfies in addition \hyperlink{om1}{$(\omega_1)$}, then as l.c.v.s.
\begin{equation}\label{mainweighfctthmequ1}
\forall\;\ell,j>0:\;\;\;\mathfrak{F}_{[\mathcal{M}_{\sigma}\cdot\mathcal{M}_{\tau}]}=\mathfrak{F}_{[\omega_{\mathbf{S}^{(\ell)}\cdot\mathbf{T}^{(j)}}]}=\mathfrak{F}_{[\sigma\check{\star}\tau]},
\end{equation}
with $\mathfrak{F}\in\{\mathcal{E}, \mathcal{B}, \mathcal{A}, \mathcal{S}, \Lambda, \mathcal{F}\}$.

\item[$(v)$] Let $\sigma$, $\tau$ be as in assertion $(iv)$. When $\tau$ satisfies in addition \hyperlink{om6}{$(\omega_6)$}, then
$\mathcal{M}_{\sigma}\cdot\mathcal{M}_{\tau}[\approx]\mathcal{M}_{\sigma}\cdot\mathbf{T}^{(j)}$ for any $j>0$ and if $\sigma$ satisfies in addition \hyperlink{om6}{$(\omega_6)$}, then $\mathcal{M}_{\sigma}\cdot\mathcal{M}_{\tau}[\approx]\mathbf{S}^{(\ell)}\cdot\mathcal{M}_{\tau}$ for any $\ell>0$. In the first case, the equalities in \eqref{mainweighfctthmequ1} can be extended by using the matrix $\mathcal{M}_{\sigma}\cdot\mathbf{T}^{(j)}$, $j>0$ arbitrary, and in the second one by using $\mathbf{S}^{(\ell)}\cdot\mathcal{M}_{\tau}$, $\ell>0$ arbitrary.
\end{itemize}
\end{theorem}

\demo{Proof}
$(i)$ By \eqref{goodequivalenceclassic} we have that $\omega_{\mathbf{S}^{(\ell)}}\hyperlink{sim}{\sim}\sigma$ and $\omega_{\mathbf{T}^{(\ell_1)}}\hyperlink{sim}{\sim}\tau$ for all $\ell,\ell_1>0$ and so, by \cite[Lemma 3.3 $(i)$]{genLegendreconj}: $$\forall\;\ell,\ell_1>0:\;\;\;\omega_{\mathbf{S}^{(\ell)}}\check{\star}\omega_{\mathbf{T}^{(\ell_1)}}\hyperlink{sim}{\sim}\sigma\check{\star}\tau.$$
Indeed, \eqref{goodequivalenceclassic} together with the given arguments in the proof of \cite[Lemma 3.3 $(i)$]{genLegendreconj} yield
\begin{align*}
&\forall\;\ell,\ell_1>0\;\exists\;D_{\ell,\sigma},D_{\ell_1,\tau}\ge 1\;\forall\;t\ge 0:
\\&
\ell'\sigma_{\mathbf{S}^{(\ell)}}\check{\star}\omega_{\mathbf{T}^{(\ell_1)}}(t)\le\sigma\check{\star}\tau(t)\le 2\ell''\omega_{\mathbf{S}^{(\ell)}}\check{\star}\omega_{\mathbf{T}^{(\ell_1)}}(t)+D_{\ell,\sigma}+D_{\ell_1,\tau},
\end{align*}
with $\ell':=\min\{\ell,\ell_1\}$ and $\ell'':=\max\{\ell,\ell_1\}$.

Moreover, \cite[Thm. 5.4]{genLegendreconj} applied to $\mathbf{M}\equiv\mathbf{S}^{(\ell)}$ and $\mathbf{N}\equiv\mathbf{T}^{(\ell_1)}$ gives $\omega_{\mathbf{S}^{(\ell)}\cdot\mathbf{T}^{(\ell_1)}}(t)=\omega_{\mathbf{S}^{(\ell)}}\check{\star}\omega_{\mathbf{T}^{(\ell_1)}}(t)$ for all $t\ge 0$. Note that $\mathbf{S}^{(\ell)}_{\iota}=+\infty=\mathbf{T}^{(\ell_1)}_{\iota}$ since $\mathbf{S}^{(\ell)},\mathbf{T}^{(\ell_1)}\in\hyperlink{LCset}{\mathcal{LC}}$ for each $\ell,\ell_1>0$.\vspace{6pt}

$(ii)$ $\mathcal{M}_{\sigma}\cdot\mathcal{M}_{\tau}$ satisfies \eqref{newmoderategrowth} because
$$\forall\;\ell>0\;\forall\;p,q\in\NN:\;\;\;S^{(\ell)}_{p+q}\cdot T^{(\ell)}_{p+q}\le S^{(2\ell)}_pS^{(2\ell)}_q\cdot T^{(2\ell)}_pT^{(2\ell)}_q.$$

$(iii)$ If $\sigma$ satisfies \hyperlink{om1}{$(\omega_1)$} then $\mathcal{M}_{\sigma}$ has \eqref{newexpabsorb}; recall $(iv)$ in Section \ref{assomatrixsection}. Thus, by the point-wise order of the sequences in the matrices:
$$\forall\;h\ge 1\;\exists\;d\ge 1\;\forall\;\ell>0\;\exists\;D\ge 1\;\forall\;p\in\NN:\;\;\;
h^pS^{(\ell)}_p\cdot T^{(\ell)}_p\le DS^{(d\ell)}_p\cdot T^{(\ell)}_p\le DS^{(d\ell)}_p\cdot T^{(d\ell)}_p.$$
If $\tau$ satisfies \hyperlink{om1}{$(\omega_1)$}, then proceed analogously.\vspace{6pt}

$(iv)$ By combining \eqref{mainweighfctthmequ}, $(ii)$ and $(iii)$ we are in position to apply \cite[Thm. 3.2, Cor. 3.17]{testfunctioncharacterization}; see also \cite[Lemma 3.1, Thm. 3.2 \& 5.4]{equalitymixedOregular}. Therefore, we get (as l.c.v.s.)
$$\forall\;\ell>0:\;\;\;\mathfrak{F}_{[\mathcal{M}_{\sigma}\cdot\mathcal{M}_{\tau}]}=\mathfrak{F}_{[\omega_{\mathbf{S}^{(\ell)}\cdot\mathbf{T}^{(\ell)}}]}=\mathfrak{F}_{[\omega\check{\star}\sigma]};$$
and by \eqref{mainweighfctthmequ0} also $\mathfrak{F}_{[\omega_{\mathbf{S}^{(\ell)}\cdot\mathbf{T}^{(\ell)}}]}=\mathfrak{F}_{[\omega_{\mathbf{S}^{(\ell_1)}\cdot\mathbf{T}^{(\ell_2)}}]}$ for all $\ell,\ell_1,\ell_2>0$.

Indeed, the cited results in \cite{testfunctioncharacterization} and \cite{equalitymixedOregular} have been shown for $\mathfrak{F}=\mathcal{E}$ but the other weighted classes follow since they are defined in terms of analogous growth requirements and so the verified growth relations of weights immediately transfer to the other weighed structures as well.\vspace{6pt}

$(v)$ This is immediate by definition and $(iii)$ in Section \ref{BMTdefsection}.
\qed\enddemo

A special case of this main statement is the following Corollary; here we are multiplying the associated weight matrix with a \emph{fixed Gevrey sequence:}

\begin{corollary}\label{mainweighfctcor}
Let $\sigma\in\hyperlink{omset0}{\mathcal{W}_0}$ be given with associated matrix $\mathcal{M}_{\sigma}:=\{\mathbf{S}^{(\ell)}: \ell>0\}$ and let $\alpha>0$. According to \eqref{productmatrixspecial} set
\begin{equation}\label{Gevreymult}
\mathcal{M}_{\sigma}\cdot\mathbf{G}^{\alpha}:=\{\mathbf{S}^{(\ell)}\cdot\mathbf{G}^{\alpha}: \ell>0\}.
\end{equation}
Then, as l.c.v.s., we get
\begin{equation}\label{mainweighfctcorequ}
\forall\;\ell>0:\;\;\;\mathfrak{F}_{[\mathcal{M}_{\sigma}\cdot\mathbf{G}^{\alpha}]}=\mathfrak{F}_{[\omega_{\mathbf{S}^{(\ell)}\cdot\mathbf{G}^{\alpha}}]}=\mathfrak{F}_{[\omega_{\mathbf{S}^{(\ell)}}\check{\star}\omega_{\mathbf{G}^{\alpha}}]}=\mathfrak{F}_{[\sigma\check{\star}\omega_{\mathbf{G}^{\alpha}}]}=\mathfrak{F}_{[\sigma\check{\star}\id^{1/\alpha}]}=\mathfrak{F}_{[(((\sigma^{\iota})^{\alpha})_{\star})^{1/\alpha}]},
\end{equation}
with $\mathfrak{F}\in\{\mathcal{E}, \mathcal{B}, \mathcal{A}, \mathcal{S}, \Lambda, \mathcal{F}\}$.
\end{corollary}

\demo{Proof}
First remark that $\id^{1/\alpha}$ satisfies both \hyperlink{om1}{$(\omega_1)$} and \hyperlink{om6}{$(\omega_6)$} but it is not normalized. However, it is equivalent to a weight $\widetilde{\id}^{1/\alpha}\in\hyperlink{omset0}{\mathcal{W}_0}$; e.g. set $\widetilde{\id}^{1/\alpha}(t):=0$ for $0\le t\le 1$ and $\widetilde{\id}^{1/\alpha}(t):=(t-1)^{1/\alpha}$ for $t>1$.   $\omega_{\mathbf{G}^{\alpha}}\hyperlink{sim}{\sim}\id^{1/\alpha}$ holds, see e.g. \cite[Example 2.9]{genLegendreconj}, and recall that $\omega_{\mathbf{G}^{\alpha}}\in\hyperlink{omset0}{\mathcal{W}_0}$, see Lemma \ref{assoweightomega0}. The equivalences yield \hyperlink{om1}{$(\omega_1)$} and \hyperlink{om6}{$(\omega_6)$} for $\widetilde{\id}^{1/\alpha}$ and $\omega_{\mathbf{G}^{\alpha}}$, too.

Then apply Theorem \ref{mainweighfctthm} to $\sigma$ and $\tau:=\omega_{\mathbf{G}^{\alpha}}$ (see \eqref{mainweighfctthmequ1}) and note that $\mathbf{G}^{\alpha}=(\mathbf{G}^{\alpha})^{(1)}$ with formally $(\mathbf{G}^{\alpha})^{(\ell)}$ denoting the $\ell$-th sequence associated with $\omega_{\mathbf{G}^{\alpha}}$; see e.g. \cite[p. 407, $(2.13)$]{subaddlike}. Indeed, by having \hyperlink{om6}{$(\omega_6)$} all sequences in the matrix $\mathcal{M}_{\omega_{\mathbf{G}^{\alpha}}}$ are equivalent to $\mathbf{G}^{\alpha}$.

Involve again \cite[Lemma 3.3 $(i)$, Thm. 5.4]{genLegendreconj} and, finally, for the last equality in \eqref{mainweighfctcorequ} recall \eqref{lowerLegendregeneral}.
\qed\enddemo

\section{An application of the generalized lower Legendre conjugate to dynamics}\label{applicationlowersection}
The aim is to give a concrete application of the main result of the previous section and we focus on a recent result dealing with the dynamics of the composition operator on weighted global spaces of ultradifferentiable functions; i.e. \emph{spaces of Gelfand-Shilov type.} We recall some basic notation and refer to \cite{ArizaFernandezGalbis25} for more details and citations:

First, let $\mathcal{E}(\RR)$ be the space of smooth functions on whole $\RR$ and for any $\omega\in\hyperlink{omset0}{\mathcal{W}_0}$ let us define the Fr\'{e}chet space
$$\mathcal{S}_{(\omega)}(\RR):=\{f\in\mathcal{E}(\RR):\;\sup_{x\in\RR,j,k\in\NN}\frac{(1+|x|)^k|f^{(j)}(x)|}{W^{(\ell)}_{j+k}}<+\infty,\;\forall\;\ell>0\}.$$
In \cite{ArizaFernandezGalbis25} for this space the notation $\mathcal{S}_{\omega}(\RR)$ is used, see \cite[Def. 2]{ArizaFernandezGalbis25} and, analogously, weighted spaces of sequences of complex numbers are introduced. According to our previously used notation we write $\Lambda_{(\omega)}$ whereas in \cite{ArizaFernandezGalbis25} the symbol $\mathcal{E}_{\omega}(\{0\})$ is used; see again \cite[Def. 2]{ArizaFernandezGalbis25}. Let us recall the space defined by the associated matrix $\mathcal{M}_{\omega}$
$$\mathcal{S}_{(\mathcal{M}_{\omega})}(\RR):=\{f\in\mathcal{E}(\RR):\;\sup_{x\in\RR,j,k\in\NN}\frac{(1+|x|)^k|f^{(j)}(x)|}{h^{j+k}W^{(\ell)}_{j+k}}<+\infty,\;\forall\;\ell,h>0\};$$
for spaces defined in terms of (abstractly given) weight matrices we refer to \cite{nuclearglobal2} and \cite{GelfandShilovincl}.
When $d>1$, then for the weight $\id^{1/d}$ in \cite{ArizaFernandezGalbis25} the notation $\Sigma_d(\RR)=\mathcal{S}_{(\id^{1/d})}(\RR)$ is used and in this case as l.c.v.s. $$\mathcal{S}_{(\id^{1/d})}(\RR)=\mathcal{S}_{(\omega_{\mathbf{G}^d})}(\RR)=\mathcal{S}_{(\mathcal{M}_{\omega_{\mathbf{G}^d}})}(\RR)=\mathcal{S}_{(\mathbf{G}^d)}(\RR),$$
since $\mathcal{M}_{\omega_{\mathbf{G}^d}}$ is constant: indeed, all sequences in this matrix are equivalent to $\mathbf{G}^d$; see the proof of Corollary \ref{mainweighfctcor}.\vspace{6pt}

In \cite{ArizaFernandezGalbis25} the authors study the composition operator $C_{\psi}: f\mapsto f\circ\psi$ acting on these spaces of Gelfand-Shilov type (in the Beurling setting) and $\psi$ is a polynomial. $\psi_n$ denotes the $n$-th iteration of $\psi$; i.e. $\psi_n:=\underbrace{\psi\circ\dots\circ\psi}_{\text{n-times}}$.

When comparing the notation from \cite[Def. 1]{ArizaFernandezGalbis25} with the conditions in Section \ref{BMTdefsection} we see that \hyperlink{om1}{$(\omega_1)$} corresponds to $(\alpha)$, \hyperlink{om3}{$(\omega_3)$} is equivalent to $(\gamma)$ and \hyperlink{om4}{$(\omega_4)$} is $(\delta)$. In addition in \cite{ArizaFernandezGalbis25} a weight is called \emph{strong} if \eqref{assostrongnq} is valid which precisely corresponds to \cite[$(1)$]{ArizaFernandezGalbis25} and occasionally this is also denoted by $(\epsilon)$ in the literature. Formally, normalization for $\omega$ is not assumed in \cite{ArizaFernandezGalbis25} but which does not restrict the generality of our considerations and can be assumed w.l.o.g. (i.e. without changing the class $\mathcal{S}_{(\omega)}(\RR)$).\vspace{6pt}

The next result generalizes \cite[Thm. 4.6 $2.$]{ArizaFernandezGalbis25}; we are involving the new information established in the previous section and follow the ideas given in the proof of \cite[Thm. 4.6 $2.$]{ArizaFernandezGalbis25}.

\begin{theorem}\label{ArizaFernandezGalbisthm}
Let $\sigma,\tau\in\hyperlink{omset0}{\mathcal{W}_0}$ be given with associated weight matrices $\mathcal{M}_{\sigma}$, $\mathcal{M}_{\tau}$. Let $0<\alpha<2$, $\psi(x):=x^2+\frac{1}{4}$, $x\in\RR$, $\mu\in\CC$ with $|\mu|>1$ and set $R_{\mu}:=\sum_{m=0}^{+\infty}\frac{C_{\psi_m}}{\mu^{m+1}}$ (i.e. $R_{\mu}=(\mu-C_{\psi})^{-1}$). We assume that
\begin{itemize}
\item[$(a)$] $\gamma(\sigma)>1$,

\item[$(b)$] the matrices $\mathcal{M}_{\sigma}$, $\mathcal{M}_{\tau}$ are related by
$$\exists\;k_0>0\;\forall\;\ell>0:\;\;\;\mathbf{T}^{(k_0)}\hyperlink{preceq}{\preceq}\mathbf{S}^{(\ell)}.$$
\end{itemize}
Then there exists $f\in\mathcal{S}_{(\sigma)}(\RR)$ such that $R_{\mu}(f)\notin\mathcal{S}_{((((\tau^{\iota})^{\alpha})_{\star})^{1/\alpha})}(\RR)=\mathcal{S}_{(\tau\check{\star}\id^{1/\alpha})}(\RR)=\mathcal{S}_{(\mathcal{M}_{\tau}\cdot\mathbf{G}^{\alpha})}(\RR)$.
\end{theorem}

Note that in \cite[Sect. 4]{ArizaFernandezGalbis25} the authors have been concerned with the study of the range of the \emph{resolvent operator} $R_{\mu}=(\mu-C_{\psi})^{-1}$, $|\mu|>1$, when $R_{\mu}$ is restricted to spaces $\mathcal{S}_{(\omega)}(\RR)$.

\demo{Proof}
\emph{Step I} By \cite[Cor. 2.13]{index} it follows that $\gamma(\sigma)>1$ is equivalent to the fact that \eqref{assostrongnq} is valid; indeed, note that \cite[Cor. 2.13]{index} even holds for weight functions in the sense of Definition \ref{weightfctdef} by \cite[Sect. 2.3, Thm. 2.11]{index}.

Thus $\sigma$ is a \emph{strong weight} in the notion of \cite{ArizaFernandezGalbis25} and this implies the fact that the \emph{Borel map} $B: \mathcal{S}_{(\sigma)}(\RR)\rightarrow\Lambda_{(\sigma)}$, $f\mapsto(f^{(j)}(0))_{j\in\NN}$, is surjective; for this see the comments given in \cite[Sect. 1, p. 3]{ArizaFernandezGalbis25}. Moreover, for $B$ being surjective the point $0$ can be replaced by any $x_0\in\RR$ (translation).

The equalities as l.c.v.s. $\mathcal{S}_{(\sigma)}(\RR)=\mathcal{S}_{(\mathcal{M}_{\sigma})}(\RR)$ and $\Lambda_{(\sigma)}=\Lambda_{(\mathcal{M}_{\sigma})}$ are valid since $\gamma(\sigma)>1$ and so $\sigma$ also satisfies \hyperlink{om1}{$(\omega_1)$}; recall Remark \ref{om1om6indexrem} and $(iv)$ in Section \ref{assomatrixsection}.\vspace{6pt}

\emph{Step II} For any $n\in\NN$ let $\delta_n:=(\delta_{n,j})_{j\in\NN}$ with $\delta_{n,j}$ denoting the ``Kronecker delta''. Now take a strictly increasing sequence of positive reals (or even integers) $(\ell_j)_{j\in\NN_{>0}}$ such that $\lim_{j\rightarrow+\infty}\ell_j=+\infty$ and $\ell_1:=1$; in the proof given in \cite{ArizaFernandezGalbis25} the authors have chosen $(\lambda_j=)\ell_j:=\log(j)$ as a strictly increasing sequence tending to infinity but one rather should put there $\lambda_j:=\log(j+1)$ since otherwise $\lambda_1=0$. By assumption $(b)$ we have $\mathbf{T}^{(k_0)}\hyperlink{preceq}{\preceq}\mathbf{S}^{(\ell)}$ for any $\ell>0$ and so (since $T^{(k_0)}_0=1=S^{(\ell)}_0$) the following holds:
$$\forall\;j\in\NN_{>0}\;\exists\;h_j\ge 1\;\forall\;p\in\NN:\;\;\;T^{(k_0)}_p\le h_j^pS^{(1/\ell_j)}_p.$$
The mapping $j\mapsto h_j$ is non-decreasing by the point-wise order of the sequences in the matrix $\mathcal{M}_{\sigma}$ and one can assume w.l.o.g. resp. one will have in general that $h_j\ge 1$ for all $j$ and that $\lim_{j\rightarrow+\infty}h_j=+\infty$. However, for any $j\in\NN_{>0}$ there exists $p_j\in\NN_{>0}$ such that $h_j\le\log(p)$ for all $p\ge p_j$ and w.l.o.g. we can assume that $j\mapsto p_j$ is strictly increasing and $\lim_{j\rightarrow+\infty}p_j=+\infty$. Summarizing, one has
\begin{equation}\label{ArizaFernandezGalbisthmStepIIequ}
\forall\;j\in\NN_{>0}\;\exists\;p_j\in\NN_{>0}\;\forall\;p\ge p_j:\;\;\;T^{(k_0)}_p\le\log(p)^pS^{(1/\ell_j)}_p.
\end{equation}
Introduce sequences $\mathbf{a}:=(a_n)_{n\in\NN}$ and $(\mathbf{b}_n)_{n\in\NN}$ with $\mathbf{b}_n=(b_{n,j})_{j\in\NN}$ given by
\begin{equation}\label{ArizaFernandezGalbisthmStepIIequ1}
a_k:=S^{(1)}_k,\;\;\;0\le k<p_1,\hspace{15pt}a_k:=S^{(1/\ell_j)}_k,\;\;\;p_j\le k<p_{j+1},\;j\in\NN_{>0},
\end{equation}
and
$$\mathbf{b}_n:=\delta_n\cdot\mathbf{a}.$$
This product between two sequences is (again) understood to be point-wise and so $b_{n,j}=0$ if $n\neq j$ and $b_{n,n}=a_n$ for any $n\in\NN$. Note that $(\mathbf{b}_n)_{n\in\NN}$ is a bounded sequence in $\Lambda_{(\sigma)}(=\Lambda_{(\mathcal{M}_{\sigma})})$: The estimate
\begin{equation}\label{ArizaFernandezGalbisthmStepIIequ2}
\forall\;\ell>0\;\exists\;C\ge 1\;\forall\;n\in\NN\;\forall\;j\in\NN:\;\;\;|b_{n,j}|=b_{n,j}\le CS^{(\ell)}_j
\end{equation}
clearly holds since the only non-trivial value is given by $b_{n,n}=a_n$, by \eqref{ArizaFernandezGalbisthmStepIIequ1}, $\lim_{j\rightarrow+\infty}\ell_j=+\infty$ and by the point-wise order of the sequences in the matrix $\mathcal{M}_{\sigma}$.\vspace{6pt}

\emph{Step III} Let $x_0\ge 2$ be arbitrary but fixed and put recursively $x^2_{n+1}:=x_n-\frac{1}{4}$ which yields a decreasing sequence $(x_n)_{n\in\NN}$ in $\RR$ tending to $\frac{1}{2}$. Since the Borel map $B: \mathcal{S}_{(\sigma)}(\RR)\rightarrow\Lambda_{(\sigma)}$ is surjective there exists a bounded sequence of functions $(f_n)_{n\in\NN}$ in $\mathcal{S}_{(\sigma)}(\RR)$ such that $B(f_n)=\mathbf{b}_n$ for all $n$ and so $f^{(j)}_n(x_0)=b_{n,j}$ for all $j,n\in\NN$. The fact that $(f_n)_n$ is bounded holds, as mentioned in the proof of \cite[Thm. 4.6 $1.$]{ArizaFernandezGalbis25}, because $(\mathbf{b}_n)_{n\in\NN}$ is bounded (via \eqref{ArizaFernandezGalbisthmStepIIequ2}) and since $\Lambda_{(\sigma)}$ is a \emph{Fr\'{e}chet nuclear space:} This property follows by \cite[Lemma 3.4]{ArizaFernandezGalbis24}. Note that this result holds for any $\sigma\in\hyperlink{omset0}{\mathcal{W}_0}$ satisfying \hyperlink{om1}{$(\omega_1)$} and hence can be applied to our situation. Indeed, \cite[Lemma 3.3]{ArizaFernandezGalbis24} precisely corresponds to the proof of\eqref{newexpabsorb}. W.l.o.g. one can assume that each $f_n$ has compact support contained in $(x_1,\psi(x_0))$ and so $f_n^{(k)}(x_j)=0=f_n^{(k)}(\psi(x_0))$ for all $k,n\in\NN$ and $j\in\NN_{>0}$.

Then, by the \emph{closed graph theorem,} when assuming that $R_{\mu}(\mathcal{S}_{(\sigma)}(\RR))\subseteq\mathcal{S}_{((((\tau^{\iota})^{\alpha})_{\star})^{1/\alpha})}(\RR)$ (even as sets), we have that the map $R_{\mu}:\mathcal{S}_{(\sigma)}(\RR)\rightarrow\mathcal{S}_{((((\tau^{\iota})^{\alpha})_{\star})^{1/\alpha})}(\RR)$ is continuous and hence $(R_{\mu}(f_n))_{n\in\NN}$ is a bounded sequence in $\mathcal{S}_{((((\tau^{\iota})^{\alpha})_{\star})^{1/\alpha})}(\RR)$. By Corollary \ref{mainweighfctcor}, see \eqref{mainweighfctcorequ}, it follows that $(R_{\mu}(f_n))_{n\in\NN}$ is a bounded sequence in $\mathcal{S}_{(\mathcal{M}_{\tau}\cdot\mathbf{G}^{\alpha})}(\RR)$ and consequently, when choosing the index $\ell=k_0$ and $h=1$, then
\begin{equation}\label{ArizaFernandezGalbisthmStepIIIequ}
\exists\;C>0\;\forall\;n\in\NN:\;\;\;|(R_{\mu}(f_n))^{(n)}(x_n)|\le C T^{(k_0)}_nn!^{\alpha}.
\end{equation}
On the other hand, the proof given in \cite[Thm. 4.6 $1.$]{ArizaFernandezGalbis25} (involving \emph{Fa\`{a} di Bruno's formula}) implies that
\begin{equation}\label{ArizaFernandezGalbisthmStepIIIequ1}
\forall\;n\in\NN_{>0}:\;\;\;|(R_{\mu}f_n)^{(n)}(x_n)|>\frac{(n+2)^{2n}}{(4|\mu|)^{n+1}}b_{n,n}=\frac{(n+2)^{2n}}{(4|\mu|)^{n+1}}a_n.
\end{equation}
Let now $n\in\NN_{>0}$ be such that $p_j\le n<p_{j+1}$ for some $j\in\NN_{>0}$. Then \eqref{ArizaFernandezGalbisthmStepIIequ1} gives $a_n=S^{(1/\ell_j)}_n$ and \eqref{ArizaFernandezGalbisthmStepIIequ} implies $T^{(k_0)}_n\le\log(n)^nS^{(1/\ell_j)}_n$. Combining this information with \eqref{ArizaFernandezGalbisthmStepIIIequ} and \eqref{ArizaFernandezGalbisthmStepIIIequ1} yields
$$\exists\;C>0\;\forall\;n,j\in\NN_{>0},\;p_j\le n<p_{j+1}:\;\;\;S^{(1/\ell_j)}_n(n+2)^{2n}\le C(4|\mu|)^{n+1}n!^{\alpha}\log(n)^nS^{(1/\ell_j)}_n.$$
By the estimates $(n+2)^{2n}\ge n^{2n}$ and $n!^{\alpha}\le n^{n\alpha}$ for all $n\in\NN_{>0}$ we infer
$$\exists\;C>0\;\forall\;n\in\NN_{>0},\;n\ge p_1:\;\;\;n^{n(2-\alpha)}\le C(4|\mu|)^{n+1}\log(n)^n.$$
Since $2-\alpha>0$ this estimate is impossible for all (large) $n$ which gives the contradiction.
\qed\enddemo

\emph{Note:}
\begin{itemize}
\item[$(*)$] Assumption $(b)$ in Theorem \ref{ArizaFernandezGalbisthm} follows when $\mathcal{M}_{\tau}\vartriangleleft\mathcal{M}_{\sigma}$ and $(b)$ should be compared with the relations studied in Lemma \ref{om1om6lemma1}.

\item[$(*)$] By assumption $(a)$ in Theorem \ref{ArizaFernandezGalbisthm} one infers \hyperlink{om1}{$(\omega_1)$} for $\sigma$, equivalently $\gamma(\sigma)>0$, and hence via $(i)$ and $(iii)$ in Lemma \ref{om1om6lemma1} the relation $\mathcal{M}_{\tau}\vartriangleleft\mathcal{M}_{\sigma}$ is equivalent to requiring $\tau\hyperlink{omvartriangle}{\vartriangleleft}\sigma$.

\item[$(*)$] Assumption $(b)$ simplifies if $\mathcal{M}_{\sigma}$ is constant, equivalently if $\sigma$ satisfies \hyperlink{om6}{$(\omega_6)$}, and so, in particular, when both matrices are constant.

\item[$(*)$] On the other hand, when $\sigma=\tau$ then assumption $(b)$ holds if and only if the matrix $\mathcal{M}_{\sigma}$ is constant; confirm Corollary \ref{strongrelationlemmaauxprimesremnew}.
\end{itemize}

We close this section by showing how Theorem \ref{ArizaFernandezGalbisthm} generalizes \cite[Thm. 4.6 $2.$]{ArizaFernandezGalbis25}: Let $d,d'>1$ be given with $d<d'<d+2$, set $\sigma=\tau:=\id^{1/d}$ and so $\gamma(\sigma)=d$ which implies $(a)$ in Theorem \ref{ArizaFernandezGalbisthm}. Put $\alpha:=d'-d$, hence $\alpha\in(0,2)$ and $$(((\sigma^{\iota})^{\alpha})_{\star})^{1/\alpha}=\id^{1/d}\check{\star}\id^{1/\alpha}\hyperlink{sim}{\sim}\omega_{\mathbf{G}^{d}}\check{\star}\omega_{\mathbf{G}^{\alpha}}=\omega_{\mathbf{G}^{d+\alpha}}=\omega_{\mathbf{G}^{d'}};$$ for this recall again \cite[Example 2.9, $(3.4)$, Lemma 3.3, Thm. 5.4]{genLegendreconj}. Then recall that the matrices $\mathcal{M}_{\tau}\cdot\mathbf{G}^{\alpha}$ and $\mathcal{M}_{\tau}$ both are constant: Indeed, all sequences in the first matrix are equivalent to $\mathbf{G}^{d+\alpha}=\mathbf{G}^{d'}$ and in the second one to $\mathbf{G}^d$; see the proof of Corollary \ref{mainweighfctcor} and $(iii)$ in Section \ref{assomatrixsection}. In particular, assumption $(b)$ in Theorem \ref{ArizaFernandezGalbisthm} is clear.

Note that in the proof of \cite[Thm. 4.6 $2.$]{ArizaFernandezGalbis25} the assumption $d<d'$ has not been stated and required but which is natural in our context to ensure $\alpha>0$. On the other hand, for this result one is interested in all (large) $d'<d+2$ and so $d<d'$ is not giving a crucial restriction on the generality of the above statement.

\section{Generalized upper Legendre conjugate}\label{uppersection}
\subsection{Preliminaries}\label{upperpreliminarysection}
In this section we proceed with the detailed study of the generalized upper Legendre conjugate $\sigma\widehat{\star}\tau$ for  $\sigma,\tau\in\hyperlink{omset0}{\mathcal{W}_0}$ and investigate the effects of this operation on the corresponding associated weight matrices $\mathcal{M}_{\sigma}$, $\mathcal{M}_{\tau}$. We summarize some immediate consequences:\vspace{6pt}

\begin{itemize}
\item[$(*)$] Note that $\tau(0)=0$ and so one has automatically $(A)$ in \cite[Sect. 4.2]{genLegendreconj}; see \cite[Lemma 4.2]{genLegendreconj}. The same holds for all (appearing) associated weight functions.

\item[$(*)$] However, in order to ensure that $\sigma\widehat{\star}\tau$ is well-defined, and hence to be at least a weight function in the sense of Definition \ref{weightfctdef}, one has to ensure \eqref{equ39}.

\item[$(*)$] Recall that $\omega_{\mathbf{M}}\widehat{\star}\omega_{\mathbf{N}}$ is well-defined if and only if $\mathbf{N}\hyperlink{mtriangle}{\vartriangleleft}\mathbf{M}$ is satisfied; see \cite[Prop. 5.6]{genLegendreconj} and also \cite[Def. 5.9]{genLegendreconj}. The detailed study of the well-definedness is provided in the next Section \ref{welldefinedsection} and there frequently this equivalence from \cite[Prop. 5.6]{genLegendreconj} is mentioned in the statements. For convenience in the proofs we are then not citing \cite[Prop. 5.6]{genLegendreconj} explicitly again.

\item[$(*)$] All crucial sequences $\mathbf{M}$ under consideration are belonging to associated weight matrices and hence are elements of the set \hyperlink{LCset}{$\mathcal{LC}$}. Therefore, with the notation $\mathbf{M}_{\iota}:=\liminf_{p\rightarrow+\infty}\left(\frac{M_p}{M_0}\right)^{1/p}$, we have that $\mathbf{M}_{\iota}=+\infty$ holds and in view of this observation note that \cite[$(5.6)$]{genLegendreconj} is then valid for all $u\ge 0$.
\end{itemize}

Moreover, concerning the special case $\sigma=\tau$ let us recall the comments from \cite[Rem. 4.4 \& 4.13]{genLegendreconj}:

\begin{itemize}
\item[$(*)$] Given $\sigma\in\hyperlink{omset0}{\mathcal{W}_0}$ then, as stated in \cite[Rem. 4.13]{genLegendreconj}, the fact that $\sigma\widehat{\star}\sigma$ is well-defined implies $\overline{\gamma}(\sigma)=+\infty$. Then note that $\overline{\gamma}(\sigma)=\overline{\gamma}(\omega_{\mathbf{S}^{(\ell)}})$ for some/any $\ell>0$ since this index is preserved under equivalence and by taking into account \eqref{goodequivalenceclassic}.

\item[$(*)$] And $\overline{\gamma}(\sigma)<+\infty$ holds if and only if $\sigma$ satisfies \hyperlink{om6}{$(\omega_6)$}, see Remark \ref{om1om6indexrem}, equivalently if and only if $\mathcal{M}_{\sigma}$ is constant; recall $(iii)$ in Section \ref{assomatrixsection}. Thus the assumption that $\sigma\widehat{\star}\sigma$ is well-defined implies the fact that $\mathcal{M}_{\sigma}$ is non-constant.

\item[$(*)$] In order to treat $\widehat{\star}$ for associated weight functions expressed by $\mathbf{S}^{(\ell)}\in\mathcal{M}_{\sigma}$ one requires naturally that $\omega_{\mathbf{S}^{(\ell_1)}}\widehat{\star}\omega_{\mathbf{S}^{(\ell)}}$ is well-defined for (at least) some indices $\ell_1,\ell$. And this is equivalent to the fact that $\mathbf{S}^{(\ell)}\hyperlink{mtriangle}{\vartriangleleft}\mathbf{S}^{(\ell_1)}$ should be valid for (at least) some indices $\ell_1>\ell$. This situation can occur but then necessarily $\mathcal{M}_{\sigma}$ has to be again non-constant.
\end{itemize}

\subsection{On the well-definedness of the generalized upper Legendre conjugate}\label{welldefinedsection}
We start with the first immediate consequence of \eqref{goodequivalenceclassic}:

\begin{lemma}\label{strongrelationlemmaaux}
Let $\sigma,\tau\in\hyperlink{omset0}{\mathcal{W}_0}$ be given with associated matrices $\mathcal{M}_{\sigma}:=\{\mathbf{S}^{(\ell)}: \ell>0\}$, $\mathcal{M}_{\tau}:=\{\mathbf{T}^{(\ell)}: \ell>0\}$. Consider the following assertions:

\begin{itemize}
\item[$(i)$] $\omega_{\mathbf{S}^{(\ell)}}\widehat{\star}\omega_{\mathbf{T}^{(2\ell)}}$ is well-defined for some/all $\ell>0$, equivalently
\item[$(i)'$] $\mathbf{T}^{(2\ell)}\hyperlink{mtriangle}{\vartriangleleft}\mathbf{S}^{(\ell)}$ for some/all $\ell>0$.

\item[$(ii)$] $\sigma\widehat{\star}\tau$ is well-defined.

\item[$(iii)$] $\omega_{\mathbf{S}^{(2\ell)}}\widehat{\star}\omega_{\mathbf{T}^{(\ell)}}$ is well-defined for some/all $\ell>0$,
equivalently
\item[$(iii)'$] $\mathbf{T}^{(\ell)}\hyperlink{mtriangle}{\vartriangleleft}\mathbf{S}^{(2\ell)}$ for some/all $\ell>0$.
\end{itemize}
Then $(i)\Rightarrow(ii)\Rightarrow(iii)$ holds.
\end{lemma}

\emph{Note:} The proof shows that when $(i)$ holds for some $\ell_0$, then $(iii)$ too. $(i)\Rightarrow(ii)$ does not use \hyperlink{om4}{$(\omega_4)$} for $\tau$, whereas this property is not required for $\sigma$ for showing $(ii)\Rightarrow(iii)$. Of course, $(i)\Rightarrow(iii)$ follows directly by involving the point-wise order of the sequences in the matrices and the definition of $\widehat{\star}$.

\demo{Proof}
We exploit \eqref{goodequivalenceclassic} for both matrices; first we get the following estimate for all $t\in(0,+\infty)$:
\begin{align*}
&\sigma\widehat{\star}\tau(t)=\sup_{s\ge 0}\{\sigma(s)-\tau(s/t)\}\le\sup_{s\ge 0}\{2\ell\omega_{\mathbf{S}^{(\ell)}}(s)-2\ell\omega_{\mathbf{T}^{(2\ell)}}(s/t)\}+D_{\ell}
\\&
=2\ell\sup_{s\ge 0}\{\omega_{\mathbf{S}^{(\ell)}}(s)-\omega_{\mathbf{T}^{(2\ell)}}(s/t)\}+D_{\ell}=2\ell\omega_{\mathbf{S}^{(\ell)}}\widehat{\star}\omega_{\mathbf{T}^{(2\ell)}}(t)+D_{\ell}.
\end{align*}
Similarly,
\begin{align*}
&\sigma\widehat{\star}\tau(t)\ge\sup_{s\ge 0}\{2\ell\omega_{\mathbf{S}^{(2\ell)}}(s)-2\ell\omega_{\mathbf{T}^{(\ell)}}(s/t)\}-D'_{\ell}
\\&
=2\ell\sup_{s\ge 0}\{\omega_{\mathbf{S}^{(2\ell)}}(s)-\omega_{\mathbf{T}^{(\ell)}}(s/t)\}-D'_{\ell}=2\ell\omega_{\mathbf{S}^{(2\ell)}}\widehat{\star}\omega_{\mathbf{T}^{(\ell)}}(t)-D'_{\ell}.
\end{align*}
Finally, for $t=0$ we recall that $\sigma\widehat{\star}\tau(0)=\sigma(0)-\tau(0)=0-0=0$ holds by \cite[Lemma 4.1 $(b)$]{genLegendreconj}, similarly this applies to the associated weight functions and hence even equality for $t=0$ follows. Summarizing,
$$\forall\;\ell>0\;\exists\;D_{\ell},D'_{\ell}\ge 1\;\forall\;t\ge 0:\;\;\;2\ell\omega_{\mathbf{S}^{(2\ell)}}\widehat{\star}\omega_{\mathbf{T}^{(\ell)}}(t)-D'_{\ell}\le\sigma\widehat{\star}\tau(t)\le 2\ell\omega_{\mathbf{S}^{(\ell)}}\widehat{\star}\omega_{\mathbf{T}^{(2\ell)}}(t)+D_{\ell},$$
which gives $(i)\Rightarrow(ii)\Rightarrow(iii)$.
\qed\enddemo

This result illustrates that relation \hyperlink{mtriangle}{$\vartriangleleft$} between sequences in the associated matrices is crucial in order to detect whether $\widehat{\star}$ is well-defined for given (associated) weight functions. On the other hand, relation $\mathcal{M}_{\tau}\vartriangleleft\mathcal{M}_{\sigma}$ seems to be too restrictive resp. is not required necessarily, since in $(i)'$ and $(iii)'$ only the relation for a certain pair of sequences is required. In view of this observation we show the next technical result:

\begin{lemma}\label{strongrelationlemmaauxprimes}
Let $\sigma,\tau\in\hyperlink{omset0}{\mathcal{W}_0}$ be given with associated matrices $\mathcal{M}_{\sigma}:=\{\mathbf{S}^{(\ell)}: \ell>0\}$, $\mathcal{M}_{\tau}:=\{\mathbf{T}^{(\ell)}: \ell>0\}$. Then the following are equivalent:
\begin{itemize}
\item[$(i)$] $\sigma$, $\tau$ satisfy
\begin{equation}\label{strongrelationlemmaauxprimesequ}
\exists\;C>0\;\forall\;H\ge 1\;\exists\;D_H\ge 1\;\forall\;t\ge 0:\;\;\;\sigma(Ht)\le C\tau(t)+D_H.
\end{equation}

\item[$(ii)$] $\mathcal{M}_{\sigma}$, $\mathcal{M}_{\tau}$ satisfy
$$\exists\;d>0\;\forall\;\ell>0:\;\;\;\mathbf{T}^{(\ell)}\hyperlink{mtriangle}{\vartriangleleft}\mathbf{S}^{(d\ell)}.$$

\item[$(iii)$] $\mathcal{M}_{\sigma}$, $\mathcal{M}_{\tau}$ satisfy
$$\exists\;\ell,\ell_1>0:\;\;\;\mathbf{T}^{(\ell)}\hyperlink{mtriangle}{\vartriangleleft}\mathbf{S}^{(\ell_1)}.$$
\end{itemize}
The proof shows the following relations: If $(i)$ holds with $C$, then $(ii)$ with $d=C$ and when $(ii)$ holds with some $d$, then it suffices to put $C=2d$ in $(i)$.
\end{lemma}

\emph{Note:} $(i)'$ in Lemma \ref{strongrelationlemmaaux} precisely corresponds to $d=\frac{1}{2}$ in assertion $(ii)$ and $(iii)'$ to $d=2$ there. By the point-wise order of the sequences, if $(ii)$ holds then $\mathbf{T}^{(\ell)}\hyperlink{mtriangle}{\vartriangleleft}\mathbf{S}^{(\ell_1)}$ for any $\ell_1\ge d\ell$ and the relation in $(iii)$ holds for all $\ell_2,\ell_3>0$ satisfying $\ell_2<\ell$ and $\ell_3>\ell_1$, too. Moreover, the proof does not require necessarily \hyperlink{om4}{$(\omega_4)$} for $\tau$.

\demo{Proof}
$(ii)\Rightarrow(iii)$ is trivial and we continue with $(iii)\Rightarrow(i)$: The assumption implies
\begin{align*}
\exists\;\ell,\ell_1>0\;\forall\;h>0\;\exists\;C_h\ge 1\;\forall\;p\in\NN\;\forall\;t\ge 0:\;\;\;\frac{t^p}{S^{(\ell_1)}_p}\le C_h\frac{(ht)^p}{T^{(\ell)}_p},
\end{align*}
and so by definition $\omega_{\mathbf{S}^{(\ell_1)}}(t)\le\omega_{\mathbf{T}^{(\ell)}}(ht)+\log(C_h)$ for all $t\ge 0$. Then \eqref{goodequivalenceclassic} gives
$$\sigma(t)\le 2\ell_1\omega_{\mathbf{S}^{(\ell_1)}}(t)+D_{\ell_1}\le 2\ell_1\omega_{\mathbf{T}^{(\ell)}}(ht)+2\ell_1\log(C_h)+D_{\ell_1}\le\frac{2\ell_1}{\ell}\tau(ht)+2\ell_1\log(C_h)+D_{\ell_1},$$
hence \eqref{strongrelationlemmaauxprimesequ} is verified with $C:=\frac{2\ell_1}{\ell}$, $H:=h^{-1}$, $D_H:=2\ell_1\log(C_h)+D_{\ell_1}$. (So, if $\ell_1=d\ell$ for some $d>0$ like in assertion $(ii)$, then $C=2d$.)\vspace{6pt}

$(i)\Rightarrow(ii)$: In \eqref{strongrelationlemmaauxprimesequ} set $t=e^s$, $H=e^h$ and so
$$\exists\;C>0\;\forall\;h\ge 0\;\exists\;D_h\ge 1\;\forall\;s\in\RR:\;\;\;\varphi_{\sigma}(h+s)=\sigma(e^he^s)\le C\tau(e^s)+D_h=C\varphi_{\tau}(s)+D_h.$$
Using this, we estimate as follows for all $x\ge 0$:
\begin{align*}
\varphi^{*}_{\tau}(x)&=\sup_{s\ge 0}\{xs-\varphi_{\tau}(s)\}=\sup_{s\in\RR}\{xs-\varphi_{\tau}(s)\}\le\sup_{s\in\RR}\{xs-\frac{1}{C}\varphi_{\sigma}(h+s)\}+\frac{D_h}{C}
\\&
=\frac{1}{C}\sup_{s\in\RR}\{(Cx)s-\varphi_{\sigma}(h+s)\}+\frac{D_h}{C}=\frac{1}{C}\sup_{u\in\RR}\{(Cx)(u-h)-\varphi_{\sigma}(u)\}+\frac{D_h}{C}
\\&
=\frac{1}{C}\sup_{u\in\RR}\{(Cx)u-\varphi_{\sigma}(u)\}-xh+\frac{D_h}{C}=\frac{1}{C}\sup_{u\ge 0}\{(Cx)u-\varphi_{\sigma}(u)\}-xh+\frac{D_h}{C}
\\&
=\frac{1}{C}\varphi^{*}_{\sigma}(Cx)-xh+\frac{D_h}{C}.
\end{align*}
For this computation also recall that both weights are normalized; however the same conclusion holds when one or both are not by naturally considering $\sup_{s\in\RR}$ then. We apply this estimate to $x=\ell p$, $\ell>0$ and $p\in\NN$ arbitrary, divide both sides by $\ell$ and apply $\exp$ in order to get:
$$\exists\;C>0\;\forall\;H\ge 1\;\exists\;D_H\ge 1\;\forall\;\ell>0\;\forall\;p\in\NN:\;\;\;T^{(\ell)}_p\le\exp(D_H/(\ell C))H^{-p}S^{(C\ell)}_p.$$
Therefore, $(ii)$ is verified with $d:=C$ and note that $d$ is neither depending on $\ell$ nor on $H$ and $H$ is not depending on $\ell$.
\qed\enddemo

Using the previous result we immediately obtain:

\begin{corollary}\label{strongrelationlemmaauxprimescor}
Let $\sigma,\tau\in\hyperlink{omset0}{\mathcal{W}_0}$ be given with associated matrices $\mathcal{M}_{\sigma}:=\{\mathbf{S}^{(\ell)}: \ell>0\}$, $\mathcal{M}_{\tau}:=\{\mathbf{T}^{(\ell)}: \ell>0\}$. Then the following are equivalent:
\begin{itemize}
\item[$(i)$] \eqref{strongrelationlemmaauxprimesequ} is valid for some $C\in(0,+\infty)$.

\item[$(ii)$] There exists $C\in(0,+\infty)$ such that $\omega_{\mathbf{S}^{(C\ell)}}\widehat{\star}\omega_{\mathbf{T}^{(\ell)}}$ is well-defined for any $\ell>0$ (and so $\omega_{\mathbf{S}^{(\ell_1)}}\widehat{\star}\omega_{\mathbf{T}^{(\ell)}}$ is well-defined for any $\ell_1\ge C\ell$).

\item[$(iii)$] There exist $\ell_1,\ell>0$ such that $\omega_{\mathbf{S}^{(\ell_1)}}\widehat{\star}\omega_{\mathbf{T}^{(\ell)}}$ is well-defined.
\end{itemize}
\end{corollary}

\demo{Proof}
We apply and combine the characterizations obtained in \cite[Prop. 5.6]{genLegendreconj} and Lemma \ref{strongrelationlemmaauxprimes}.
\qed\enddemo

\begin{remark}\label{strongrelationlemmaauxprimesrem}
\emph{\eqref{strongrelationlemmaauxprimesequ} should be compared with \eqref{equ39}. Indeed, if $C\in(0,1]$ in \eqref{strongrelationlemmaauxprimesequ} then \eqref{equ39} follows and this condition implies \eqref{strongrelationlemmaauxprimesequ} for any $C\in[1,+\infty)$.}
\end{remark}

When considering in the above Lemma \ref{strongrelationlemmaauxprimes} the same weight function and following the ideas of the proof, then we immediately infer:

\begin{corollary}\label{strongrelationlemmaauxprimesremnew}
Let $\omega\in\hyperlink{omset0}{\mathcal{W}_0}$ be given, then the following are equivalent:
\begin{itemize}
\item[$(i)$] $\omega$ satisfies \hyperlink{om6}{$(\omega_6)$}.

\item[$(ii)$] $\mathcal{M}_{\omega}$ is constant.

\item[$(iii)$] $\mathcal{M}_{\omega}$ satisfies
$$\exists\;\ell>2\ell_1>0:\;\;\;\mathbf{W}^{(\ell)}\hyperlink{preceq}{\preceq}\mathbf{W}^{(\ell_1)}.$$
\end{itemize}
\end{corollary}

\emph{Note:} This result slightly improves the proof of \cite[Lemma 5.9]{compositionpaper} where we have chosen $\ell:=4\ell_1$ for proving the analogous implication $(iii)\Rightarrow(i)$. In \cite[Prop. 4.10]{modgrowthstrangeII} the implication $(iii)\Rightarrow(i)$ has already been mentioned; there inspired by similar ideas and techniques involving the convolved sequence appearing in \cite[Sect. 3]{modgrowthstrangeII}.

\demo{Proof}
Indeed, $(i)\Rightarrow(ii)$ is shown in \cite[Lemma 5.9]{compositionpaper}, $(ii)\Rightarrow(iii)$ is trivial and for $(iii)\Rightarrow(i)$ follow the proof of $(iii)\Rightarrow(i)$ in Lemma \ref{strongrelationlemmaauxprimes}: Therefore, $(iii)$ implies $\frac{\ell}{2\ell_1}\omega(t)\le\omega(Ht)+H$ for some $H\ge 1$ and all $t\ge 0$ and then \hyperlink{om6}{$(\omega_6)$} follows by iteration; see also \cite[Rem. 2.2]{modgrowthstrangeII}.
\qed\enddemo

Let us now comment on more statements concerning the special case $\sigma=\tau$ in Lemmas \ref{strongrelationlemmaaux} and \ref{strongrelationlemmaauxprimes}:

\begin{itemize}
\item[$(*)$] In this case $(i)'$ in Lemma \ref{strongrelationlemmaaux} is violated since this assertion gives $\mathbf{S}^{(2\ell)}\hyperlink{mtriangle}{\vartriangleleft}\mathbf{S}^{(\ell)}\le\mathbf{S}^{(2\ell)}$ which fails for any $\ell>0$.

\item[$(*)$] By analogous reasons one infers that in \eqref{strongrelationlemmaauxprimesequ} one has to choose $C>1$ resp. $d>1$ in $(ii)$ and $\ell_1>\ell$ in $(iii)$ in Lemma \ref{strongrelationlemmaauxprimes}. Therefore, \eqref{strongrelationlemmaauxprimesequ} for $\sigma$ does not imply necessarily \eqref{equ39}; see also Remark \ref{strongrelationlemmaauxprimesrem} concerning the role of the parameter $C$ in this condition.

\item[$(*)$] $(iii)'$ in Lemma \ref{strongrelationlemmaaux} turns into $\mathbf{S}^{(\ell)}\hyperlink{mtriangle}{\vartriangleleft}\mathbf{S}^{(2\ell)}$ for some/all $\ell>0$. Hence, for $(iii)'$ being valid, \hyperlink{om6}{$(\omega_6)$} for $\sigma$ is impossible since in this case $\mathcal{M}_{\sigma}$ has to be constant and so all sequences are equivalent; see $(iii)$ in Section \ref{assomatrixsection}. Note that this observation and $(ii)\Rightarrow(iii)\Leftrightarrow(iii)'$ in Lemma \ref{strongrelationlemmaaux} is consistent with \cite[Rem. 4.4 \& 4.8, Cor. 5.7]{genLegendreconj}.

\item[$(*)$] By the characterization in Lemma \ref{strongrelationlemmaauxprimes}, when \eqref{strongrelationlemmaauxprimesequ} holds for $\sigma$ (with $C>1$), then \hyperlink{om6}{$(\omega_6)$} is impossible. Thus, by taking into account Remark \ref{om1om6indexrem}, one infers $\overline{\gamma}(\sigma)=+\infty$.

    The fact that \hyperlink{om6}{$(\omega_6)$} fails can also be seen directly for any weight function in the sense of Definition \ref{weightfctdef} by an iterated application of \hyperlink{om6}{$(\omega_6)$}: Let $n\in\NN_{>0}$ be arbitrary (but fixed), then by \hyperlink{om6}{$(\omega_6)$} one can find $H\ge 1$ such that $2^n\sigma(t)\le\sigma(Ht)+H$ for all $t\ge 0$ and \eqref{strongrelationlemmaauxprimesequ} applied to this $H$ implies $2^n\sigma(t)\le C\sigma(t)+D_H+H$ for some $C\ge 1$ not depending on $H$ and all $t\ge 0$. But this estimate gives a contradiction if $2^n>C$ since $\lim_{t\rightarrow+\infty}\sigma(t)=+\infty$.

\item[$(*)$] Let $\sigma$ by a weight function according to Definition \ref{weightfctdef}. Then \eqref{strongrelationlemmaauxprimesequ} for $\sigma=\tau$ should be compared with the following known growth condition
\begin{equation}\label{om7}
\exists\;H,C>0\;\forall\;t\ge 0:\;\;\;\sigma(t^2)\le C\sigma(Ht)+C,
\end{equation}
denoted by $(\omega_7)$ e.g. in \cite{sectorialextensions1}. Since it implies \hyperlink{om1}{$(\omega_1)$}, see \cite[Appendix A]{sectorialextensions1}, this property is equivalent to having $\sigma(t^2)\le C\sigma(t)+C$ for some $C\ge 1$ and all $t\ge 0$. It is known that $(\omega_7)$ implies $\gamma(\sigma)=+\infty$, see \cite[Lemma A.1]{sectorialextensions1}, and let us verify
\begin{equation}\label{crucialchain}
\eqref{om7}\Longrightarrow\eqref{strongrelationlemmaauxprimesequ}\Longrightarrow\gamma(\sigma)=+\infty(\Longrightarrow\overline{\gamma}(\sigma)=+\infty).
\end{equation}
The first implication is clear since $\sigma$ is non-decreasing and so $\sigma(Ht)\le\sigma(t^2)$ for all $t\ge t_H:=H$. For the second implication even a uniform choice for $K$ in $(P_{\sigma,\gamma})$, $\gamma>0$ arbitrary, is sufficient: Consider a fixed $\epsilon>0$, set $K:=(1+\epsilon)C$ with $C$ denoting the uniform constant from \eqref{strongrelationlemmaauxprimesequ}. Let $\gamma>0$ be arbitrary (large), then apply \eqref{strongrelationlemmaauxprimesequ} to $H:=K^{\gamma}$ and so $(P_{\sigma,\gamma})$ holds. The last implication is trivial and is consistent with the previous comment above.

Summarizing, in view of Remark \ref{strongrelationlemmaauxprimesrem} the implications verified in \eqref{crucialchain} are again consistent with \cite[Rem. 4.4 \& 4.8]{genLegendreconj}. For this recall that, as stated in these remarks, the assumption that $\sigma\widehat{\star}\sigma$ is a weight function implies that $\sigma$ has to be \emph{slowly varying} and hence $\gamma(\sigma)=+\infty$.

\item[$(*)$] The equivalent assertions stated in Lemma \ref{strongrelationlemmaauxprimes} can be extended by the following properties introduced in \cite[Sect. 4.1]{compositionpaper} for abstractly given weight matrices:
$$(\mathcal{M}_{\{\on{BR}\}}):\;\;\;\forall\;\ell>0\;\exists\;\ell_1>0:\;\;\;\mathbf{S}^{(\ell)}\hyperlink{mtriangle}{\vartriangleleft}\mathbf{S}^{(\ell_1)},$$
$$(\mathcal{M}_{(\on{BR})}):\;\;\;\forall\;\ell_1>0\;\exists\;\ell>0:\;\;\;\mathbf{S}^{(\ell)}\hyperlink{mtriangle}{\vartriangleleft}\mathbf{S}^{(\ell_1)}.$$
Obviously, $(ii)$ in Lemma \ref{strongrelationlemmaauxprimes} implies both $(\mathcal{M}_{\{\on{BR}\}})$ and $(\mathcal{M}_{(\on{BR})})$, whereas any of these conditions yields $(iii)$. Note that in the above conditions $\ell_1>\ell$ is necessary by the point-wise order of the sequences and any (abstractly given) weight matrix satisfying one of these assertions cannot be constant; independently this verifies that \hyperlink{om6}{$(\omega_6)$} has to be violated.
\end{itemize}

In the next statement we are using a precise relation for the indices $\gamma(\sigma)$ and $\overline{\gamma}(\tau)$ of the involved weights. Recall that $\gamma(\sigma)>0$ if and only if \hyperlink{om1}{$(\omega_1)$} and $\overline{\gamma}(\tau)<+\infty$ if and only if \hyperlink{om6}{$(\omega_6)$}.

\begin{lemma}\label{strongrelationlemmaauxvar}
Let $\sigma,\tau\in\hyperlink{omset0}{\mathcal{W}_0}$ be given with associated matrices $\mathcal{M}_{\sigma}:=\{\mathbf{S}^{(\ell)}: \ell>0\}$, $\mathcal{M}_{\tau}:=\{\mathbf{T}^{(\ell)}: \ell>0\}$. Assume that
\begin{equation}\label{strongrelationlemmaauxvarequ}
0\le\overline{\gamma}(\tau)<\gamma(\sigma)\le+\infty
\end{equation}
is satisfied. Then the following are equivalent:
\begin{itemize}
\item[$(i)$] $\sigma\widehat{\star}\tau$ is well-defined.

\item[$(ii)$] $\omega_{\mathbf{S}^{(\ell)}}\widehat{\star}\tau$ is well-defined for some/all $\ell>0$.

\item[$(iii)$]  $\sigma\widehat{\star}\omega_{\mathbf{T}^{(\ell)}}$ is well-defined for some/all $\ell>0$.

\item[$(iv)$] $\omega_{\mathbf{S}^{(\ell)}}\widehat{\star}\omega_{\mathbf{T}^{(\ell_1)}}$ is well-defined for all $\ell,\ell_1>0$ equivalently $\mathbf{T}^{(\ell_1)}\hyperlink{mtriangle}{\vartriangleleft}\mathbf{S}^{(\ell)}$ for all $\ell,\ell_1>0$; i.e. $\mathcal{M}_{\tau}\vartriangleleft\mathcal{M}_{\sigma}$ holds.
\end{itemize}
Moreover, if any of these stated equivalent assertions holds, then $$\forall\;\ell,\ell_1,\ell_2,\ell_3>0:\;\;\;\sigma\widehat{\star}\tau\hyperlink{sim}{\sim}\omega_{\mathbf{S}^{(\ell)}}\widehat{\star}\tau\hyperlink{sim}{\sim}\sigma\widehat{\star}\omega_{\mathbf{T}^{(\ell_1)}}\hyperlink{sim}{\sim}\omega_{\mathbf{S}^{(\ell_2)}}\widehat{\star}\omega_{\mathbf{T}^{(\ell_3)}}.$$
\end{lemma}

\demo{Proof}
Since the growth indices are preserved under equivalence, by assumption and using \eqref{goodequivalenceclassic} one has $$\forall\;\ell,\ell_1>0:\;\;\;0\le\overline{\gamma}(\omega_{\mathbf{T}^{(\ell)}})=\overline{\gamma}(\tau)<\gamma(\sigma)=\gamma(\omega_{\mathbf{S}^{(\ell_1)}})\le+\infty.$$

Then the desired equivalences follow by applying \cite[Thm. 4.14]{genLegendreconj}.
\qed\enddemo

We comment on the failure of Lemma \ref{strongrelationlemmaauxvar} when $\sigma=\tau$, especially on \eqref{strongrelationlemmaauxvarequ}:
 
\begin{remark}\label{strongrelationlemmaauxvarrem}
\emph{When $\sigma=\tau$, then \eqref{strongrelationlemmaauxvarequ} fails since $\overline{\gamma}(\sigma)<\gamma(\sigma)$ and also $\mathcal{M}_{\sigma}\vartriangleleft\mathcal{M}_{\sigma}$ is impossible. Moreover, for having $(i)$ Lemma \ref{strongrelationlemmaauxvar} necessarily $\sigma$ has to be \emph{slowly varying} but this contradicts \hyperlink{om6}{$(\omega_6)$} and so $\overline{\gamma}(\sigma)=+\infty$ follows; see \cite[Rem. 4.4 \& 4.8]{genLegendreconj} and the citations there.}

\emph{Similarly, this result cannot be applied to (different) associated weight functions $\omega_{\mathbf{W}^{(\ell)}}$, $\omega_{\mathbf{W}^{(\ell_1)}}$ when $\omega$ satisfies $\overline{\gamma}(\omega)=+\infty$: This value is preserved under equivalence and hence by \eqref{goodequivalenceclassic} each associated weight function shares the same value for $\overline{\gamma}(\cdot)$.}

\emph{Since any slowly varying weight function $\omega$ satisfies $\gamma(\omega)=\overline{\gamma}(\omega)=+\infty$ we see that the estimates shown in \cite[Thm. 4.10]{genLegendreconj} formally remain valid when treating $\sigma=\tau$: In this case \cite[$(4.14)\&(4.15)$]{genLegendreconj} turn into the ``equalities'' $+\infty=+\infty+(+\infty)$; see \cite[Rem. 4.13]{genLegendreconj}.}
\end{remark}

Inspired by assertion $(iv)$ in the previous lemma we study in more detail the strong growth relation $\mathcal{M}_{\tau}\vartriangleleft\mathcal{M}_{\sigma}$ but without making the assumption on the growth indices in Lemma \ref{strongrelationlemmaauxvar}. In this case, $\omega_{\mathbf{S}^{(\ell_1)}}\widehat{\star}\omega_{\mathbf{T}^{(\ell)}}$ is well-defined \emph{for any indices $\ell_1,\ell>0$;} see the comments in Section \ref{upperpreliminarysection}. Note that $\mathcal{M}_{\sigma}\vartriangleleft\mathcal{M}_{\sigma}$ can never be valid. Now we establish a precise connection to the known natural growth relation $\tau\hyperlink{omvartriangle}{\vartriangleleft}\sigma$ which is crucial for having resp. characterizing the inclusion $\mathfrak{F}_{\{\tau\}}\subseteq\mathfrak{F}_{(\sigma)}$; we also refer to \cite[Lemma 2.2]{GelfandShilovincl} and \cite[Lemma 5.16, Cor. 5.17]{compositionpaper} (and the characterizations shown in these works), but there formally more restrictive assumptions on the weights have been considered.

\begin{lemma}\label{om1om6lemma1}
Let $\sigma,\tau\in\hyperlink{omset0}{\mathcal{W}_0}$ and let $\mathcal{M}_{\sigma}$, $\mathcal{M}_{\tau}$ be the associated weight matrices.
\begin{itemize}
\item[$(i)$] If either $\sigma$ or $\tau$ satisfies \hyperlink{om1}{$(\omega_1)$}, then $\tau\hyperlink{omvartriangle}{\vartriangleleft}\sigma$ implies $\mathcal{M}_{\tau}\vartriangleleft\mathcal{M}_{\sigma}$.

\item[$(ii)$] $\tau\hyperlink{ompreceq}{\preceq}\sigma$ implies $\mathcal{M}_{\tau}\{\preceq\}\mathcal{M}_{\sigma}$ and $\mathcal{M}_{\tau}(\preceq)\mathcal{M}_{\sigma}$.

\item[$(iii)$] $\mathcal{M}_{\tau}\vartriangleleft\mathcal{M}_{\sigma}$ implies $\tau\hyperlink{omvartriangle}{\vartriangleleft}\sigma$.

\item[$(iv)$] If either $\sigma$ or $\tau$ satisfies \hyperlink{om1}{$(\omega_1)$}, then
\begin{equation}\label{om1om6lemma1equ}
\exists\;j,\ell>0:\;\;\;\mathbf{T}^{(j)}\hyperlink{preceq}{\preceq}\mathbf{S}^{(\ell)}
\end{equation}
implies $\tau\hyperlink{ompreceq}{\preceq}\sigma$.
\end{itemize}
\end{lemma}

\emph{Note:} In $(iii)$ and $(iv)$ we can allow $\tau\in\hyperlink{omset0}{\mathcal{W}_0}\backslash\{\hyperlink{om4}{(\omega_4)}\}$ and in $(i)$ and $(ii)$ even $\sigma,\tau\in\hyperlink{omset0}{\mathcal{W}_0}\backslash\{\hyperlink{om4}{(\omega_4)}\}$. Relation \eqref{om1om6lemma1equ} is precisely $\mathcal{M}_{\tau}(\preceq\}\mathcal{M}_{\sigma}$ introduced in \cite[Sect. 4.2, p. 111]{compositionpaper} and it follows from both $\mathcal{M}_{\tau}\{\preceq\}\mathcal{M}_{\sigma}$ and $\mathcal{M}_{\tau}(\preceq)\mathcal{M}_{\sigma}$.

\demo{Proof}
$(i)$ In \cite[Lemma 5.16 $(2)$]{compositionpaper} it was shown that $\tau\hyperlink{omvartriangle}{\vartriangleleft}\sigma$ implies
\begin{equation}\label{Lemma516}
\forall\;H>0\;\forall\;\ell>0\;\exists\;A\ge 1\;\forall\;p\in\NN:\;\:\:T^{(\ell)}_p\le AS^{(H\ell)}_p.
\end{equation}
Let $h\ge 1$, $\ell,j>0$ be given, arbitrary but fixed. If $\tau$ satisfies \hyperlink{om1}{$(\omega_1)$}, then involve \eqref{newexpabsorb} and get some $d,D\ge 1$, $d$ depending on $h$ and $D$ on $h$ and $\ell$, such that $h^pT^{(\ell)}_p\le DT^{(d\ell)}_p$ for all $p\in\NN$. Then \eqref{Lemma516} applied to $H:=\frac{j}{d\ell}$ yields for all $p\in\NN$:
$$h^pT^{(\ell)}_p\le DT^{(d\ell)}_p\le ADS^{(j)}_p.$$
The constant $A$ is depending on $h$, $\ell$ and $j$ and this shows $\mathbf{T}^{(j)}\hyperlink{mtriangle}{\vartriangleleft}\mathbf{S}^{(\ell)}$ for any $\ell,j>0$; i.e. $\mathcal{M}_{\tau}\vartriangleleft\mathcal{M}_{\sigma}$.

If $\sigma$ satisfies \hyperlink{om1}{$(\omega_1)$}, then the conclusion follows similarly by
$$h^pT^{(\ell)}_p\le Ah^pS^{(H\ell)}_p\le ADS^{dH\ell}_p=ADS^{(j)}_p.$$

$(ii)$ In \cite[Lemma 5.16 $(1)$]{compositionpaper} it was shown that $\tau\hyperlink{ompreceq}{\preceq}\sigma$ implies
\begin{equation}\label{Lemma516weak}
\exists\;H>0\;\forall\;\ell>0\;\exists\;A\ge 1\;\forall\;p\in\NN:\;\:\:T^{(\ell)}_p\le AS^{(H\ell)}_p,
\end{equation}
which clearly yields both $\mathcal{M}_{\tau}\{\preceq\}\mathcal{M}_{\sigma}$ and $\mathcal{M}_{\tau}(\preceq)\mathcal{M}_{\sigma}$.

A careful inspection of the proofs in \cite[Lemma 5.9 $(5.10)$, Lemma 5.16]{compositionpaper} gives that \hyperlink{om4}{$(\omega_4)$} is not required necessarily for these steps; see also \cite[Sect. 3.1]{dissertation} and especially comments $(iii)-(v)$ there, and more recently \cite[Sect. 6 $(III)$]{modgrowthstrangeII}.\vspace{6pt}

$(iii)$ $\mathcal{M}_{\tau}\vartriangleleft\mathcal{M}_{\sigma}$ implies
$$\forall\;\ell,j,h>0\;\exists\;C=C_{\ell,j,h}\ge 1\;\forall\;p\in\NN\;\forall\;t\ge 0:\;\;\;\frac{t^p}{S^{(j)}_p}\le C\frac{(th)^p}{T^{(\ell)}_p},$$
hence $\omega_{\mathbf{S}^{(j)}}(t)\le\omega_{\mathbf{T}^{(\ell)}}(ht)+\log(C)$. And so \eqref{goodequivalenceclassic} yields for all $t\ge 0$:
$$\sigma(t)\le 2j\omega_{\mathbf{S}^{(j)}}(t)+D_{j}\le 2j\omega_{\mathbf{T}^{(\ell)}}(ht)+2j\log(C)+D_{j}\le\frac{2j}{\ell}\tau(ht)+2j\log(C)+D_{j}.$$
Then choose $h:=1$ and $\ell:=2$ and let $j\rightarrow0$.\vspace{6pt}

$(iv)$ Analogously, \eqref{om1om6lemma1equ} implies
$$\exists\;j,\ell>0\;\exists\;h\ge 1\;\exists\;C\ge 1\;\forall\;p\in\NN\;\forall\;t\ge 0:\;\;\;\frac{t^p}{S^{(j)}_p}\le C\frac{(th)^p}{T^{(\ell)}_p},$$
hence $\omega_{\mathbf{S}^{(j)}}(t)\le\omega_{\mathbf{T}^{(\ell)}}(ht)+\log(C)$ and via \eqref{goodequivalenceclassic} similarly as in $(iii)$ one has $\sigma(t)\le\frac{2j}{\ell}\tau(ht)+2j\log(C)+D_{j}$ for all $t\ge 0$. Let now $n\in\NN_{>0}$ be minimal such that $h\le 2^n$.

If $\tau$ satisfies \hyperlink{om1}{$(\omega_1)$}, then iterating this property $n$-times gives $\sigma(t)\le\frac{2jL}{\ell}\tau(t)+\frac{2jL}{\ell}+2j\log(C)+D_{j}$ for some $L\ge 1$ and all $t\ge 0$ and so $\tau\hyperlink{ompreceq}{\preceq}\sigma$ is verified.

If $\sigma$ satisfies \hyperlink{om1}{$(\omega_1)$}, then similarly $\sigma(t)\le L\sigma(h^{-1}t)+L\le\frac{2jL}{\ell}\tau(t)+2jL\log(C)+LD_{j}+L$ for all $t\ge 0$ and so $\tau\hyperlink{ompreceq}{\preceq}\sigma$ holds again.
\qed\enddemo

Exploiting this information, the next result relates the property for $\sigma\widehat{\star}\tau$ being well-defined to relation $\tau\hyperlink{omvartriangle}{\vartriangleleft}\sigma$; i.e. to a condition directly involving the weight functions and not their associated matrices.

\begin{proposition}\label{strongrelationlemma}
Let $\sigma,\tau\in\hyperlink{omset0}{\mathcal{W}_0}$ be given with associated matrices $\mathcal{M}_{\sigma}:=\{\mathbf{S}^{(\ell)}: \ell>0\}$, $\mathcal{M}_{\tau}:=\{\mathbf{T}^{(\ell)}: \ell>0\}$. Consider the following assertions:
\begin{itemize}
\item[$(i)$] $\tau\hyperlink{omvartriangle}{\vartriangleleft}\sigma$ holds,

\item[$(ii)$] $\omega_{\mathbf{T}^{(\ell_0)}}\hyperlink{omvartriangle}{\vartriangleleft}\omega_{\mathbf{S}^{(j_0)}}$ holds for some $\ell_0,j_0>0$.

\item[$(iii)$] $\omega_{\mathbf{T}^{(\ell)}}\hyperlink{omvartriangle}{\vartriangleleft}\omega_{\mathbf{S}^{(j)}}$ holds for any $\ell,j>0$.

\item[$(iv)$] $\sigma\widehat{\star}\tau$ is well-defined (i.e. \eqref{equ39} holds).

\item[$(v)$] $\omega_{\mathbf{S}^{(\ell_0)}}\widehat{\star}\omega_{\mathbf{T}^{(j_0)}}$ is well-defined for some $\ell_0,j_0>0$, equivalently $\mathbf{T}^{(j_0)}\hyperlink{mtriangle}{\vartriangleleft}\mathbf{S}^{(\ell_0)}$ for some $\ell_0,j_0>0$.

\item[$(vi)$] $\omega_{\mathbf{S}^{(\ell)}}\widehat{\star}\omega_{\mathbf{T}^{(j)}}$ is well-defined for any $\ell,j>0$, equivalently $\mathcal{M}_{\tau}\vartriangleleft\mathcal{M}_{\sigma}$.
\end{itemize}

Then $(i)\Leftrightarrow(ii)\Leftrightarrow(iii)$ and $(vi)\Rightarrow(i),(v)$; if either $\sigma$ or $\tau$ satisfies \hyperlink{om1}{$(\omega_1)$}, then $(i)\Rightarrow(iv)$, $(ii)\Rightarrow(v)$, $(iii)\Rightarrow(vi)$; if either $\sigma$ or $\tau$ satisfies \hyperlink{om6}{$(\omega_6)$}, then $(vi)\Rightarrow(iii)$, $(v)\Rightarrow(ii)$, $(iv)\Rightarrow(i)$.
\end{proposition}

Summarizing, if either $\sigma$ or $\tau$ satisfies \hyperlink{om1}{$(\omega_1)$} and either $\sigma$ or $\tau$ satisfies \hyperlink{om6}{$(\omega_6)$}, then all listed assertions are equivalent.

\demo{Proof}
First, $(i)\Leftrightarrow(ii)\Leftrightarrow(iii)$ follows by \eqref{goodequivalenceclassic} and the fact that relation \hyperlink{omvartriangle}{$\vartriangleleft$} is preserved under \hyperlink{sim}{$\sim$}.\vspace{6pt}

The stated equivalences in $(v)$ and $(vi)$ are valid by recalling \cite[Prop. 5.6]{genLegendreconj}. Consequently, $(vi)\Rightarrow(i)$ follows by $(iii)$ in Lemma \ref{om1om6lemma1} and the implication $(vi)\Rightarrow(v)$ is trivial.

Moreover, \hyperlink{om1}{$(\omega_1)$} and \hyperlink{om6}{$(\omega_6)$} are preserved under \hyperlink{sim}{$\sim$}, hence $\sigma$ resp. $\tau$ satisfies one of these conditions if and only if some/any $\omega_{\mathbf{S}^{(\ell)}}$ resp. $\omega_{\mathbf{T}^{(\ell)}}$ does so.

Then $(i)\Rightarrow(iv)$, $(ii)\Rightarrow(v)$, $(iii)\Rightarrow(vi)$ follows by \cite[Lemma 4.7 $(ii)\Rightarrow(iii)$]{genLegendreconj}.

Finally, for $(vi)\Rightarrow(iii)$, $(v)\Rightarrow(ii)$, $(iv)\Rightarrow(i)$ we apply \cite[Lemma 4.3 $(ii)\Rightarrow(iii)$]{genLegendreconj} and \cite[Lemma 4.7 $(i)\Rightarrow(ii)$]{genLegendreconj}.
\qed\enddemo

We apply this result to the special but important case $\tau=\id^{1/\alpha}$.

\begin{corollary}\label{strongrelationlemmacor}
Let $\sigma\in\hyperlink{omset0}{\mathcal{W}_0}$ be given with associated matrix $\mathcal{M}_{\sigma}:=\{\mathbf{S}^{(\ell)}: \ell>0\}$ and let $\alpha>0$. Then the following assertions are equivalent:
\begin{itemize}
\item[$(i)$] $\sigma(t)=o(t^{1/\alpha})$ as $t\rightarrow+\infty$.

\item[$(ii)$] $\sigma(t)=o(\omega_{\mathbf{G}^{\alpha}}(t))$ as $t\rightarrow+\infty$.

\item[$(iii)$] $\omega_{\mathbf{G}^{\alpha}}\hyperlink{omvartriangle}{\vartriangleleft}\omega_{\mathbf{S}^{(j)}}$ holds for some/any $j>0$.

\item[$(iv)$] $\sigma\widehat{\star}\id^{1/\alpha}$ is well-defined (i.e. \eqref{equ39} holds).

\item[$(v)$] $\sigma\widehat{\star}\omega_{\mathbf{G}^{\alpha}}$ is well-defined.

\item[$(vi)$] $\omega_{\mathbf{S}^{(\ell)}}\widehat{\star}\omega_{\mathbf{G}^{\alpha}}$ is well-defined for some/any $\ell>0$.

\item[$(vii)$] $\mathbf{G}^{\alpha}\hyperlink{mtriangle}{\vartriangleleft}\mathbf{S}^{(\ell)}$ for some/any $\ell>0$.
\end{itemize}
\end{corollary}

\demo{Proof}
Like in the proof of Corollary \ref{mainweighfctcor} remark that $\id^{1/\alpha}$ is not normalized but $\widetilde{\id}^{1/\alpha}\in\hyperlink{omset0}{\mathcal{W}_0}$, $\omega_{\mathbf{G}^{\alpha}}\hyperlink{sim}{\sim}\id^{1/\alpha}\hyperlink{sim}{\sim}\widetilde{\id}^{1/\alpha}$ and $\omega_{\mathbf{G}^{\alpha}}\in\hyperlink{omset0}{\mathcal{W}_0}$ satisfies both \hyperlink{om1}{$(\omega_1)$} and \hyperlink{om6}{$(\omega_6)$}.\vspace{6pt}

This and \eqref{goodequivalenceclassic} imply $(i)\Leftrightarrow(ii)\Leftrightarrow(iii)$; note that $\sigma(t)=o(\omega_{\mathbf{G}^{\alpha}}(t))$ as $t\rightarrow+\infty$ precisely means $\omega_{\mathbf{G}^{\alpha}}\hyperlink{omvartriangle}{\vartriangleleft}\sigma$ and equivalence preserves this relation.

Moreover, one can apply Proposition \ref{strongrelationlemma} to $\sigma$ and $\tau:=\omega_{\mathbf{G}^{\alpha}}$ and so the equivalences follow by taking into account the facts that the matrix $\mathcal{M}_{\omega_{\mathbf{G}^{\alpha}}}$ is constant (in view of \hyperlink{om6}{$(\omega_6)$}) and that $\mathbf{G}^{\alpha}=(\mathbf{G}^{\alpha})^{(1)}$ with $(\mathbf{G}^{\alpha})^{(\ell)}$ denoting the $\ell$-th sequence associated with $\omega_{\mathbf{G}^{\alpha}}$; see again the proof of Corollary \ref{mainweighfctcor} and \cite[p. 407, $(2.13)$]{subaddlike}.
\qed\enddemo

We gather the content of this section in the following remark:

\begin{remark}\label{closingremark}
\emph{Relation \hyperlink{omvartriangle}{$\vartriangleleft$} between weight functions is known and natural in the setting of Braun-Meise-Taylor weight functions, see \cite[Lemma 5.16, Cor. 5.17]{compositionpaper} for the characterization of inclusion relations in the ultradifferentiable setting, and it is straight-forward to check. However, it is not reflexive and so we cannot treat $\sigma=\tau$ in Proposition \ref{strongrelationlemma} and $\sigma=\id^{1/\alpha}$ in Corollary \ref{strongrelationlemmacor}.}

\emph{On the other hand, due to the appearance of the weight matrix parameter one can consider weights $\omega$ such that the (non-constant) associated matrix $\mathcal{M}_{\omega}$ allows for having ``enough space'' when modifying this parameter. More precisely, it is possible to get relation \hyperlink{mtriangle}{$\vartriangleleft$} between sequences $\mathbf{W}^{(\ell)},\mathbf{W}^{(\ell_1)}\in\mathcal{M}_{\omega}$ and so to investigate the generalized upper Legendre conjugate transform between associate weight functions when the matrix parameter is chosen properly. Corollary \ref{strongrelationlemmaauxprimescor} together with Remark \ref{strongrelationlemmaauxprimesrem} illustrate this situation.}

\emph{However, for any given weight $\omega$ one can never expect that $\omega_{\mathbf{W}^{(\ell)}}\widehat{\star}\omega_{\mathbf{W}^{(j)}}$ is a weight function for all arbitrary $\ell,j>0$: Recall that this requirement is equivalent to $\mathbf{W}^{(\ell)}\hyperlink{mtriangle}{\vartriangleleft}\mathbf{W}^{(j)}$ for all $\ell,j>0$ which is clearly impossible. Therefore, in this situation one loses uniformity on the indices.}
\end{remark}

\subsection{On the growth indices}
The aim is to focus now on the effects of the generalized upper Legendre conjugate for weight functions in the sense of Braun-Meise-Taylor concerning the growth indices $\gamma(\cdot)$ and $\overline{\gamma}(\cdot)$.

First, when being well-defined, the conclusions from \cite[Thm. 4.10, Cor. 4.11, Rem. 4.12 \& 4.13]{genLegendreconj} directly apply to all $\sigma,\tau\in\hyperlink{omset0}{\mathcal{W}_0}$ and to all the corresponding associated weight functions $\omega_{\mathbf{S}^{(\ell)}}$, $\omega_{\mathbf{T}^{(\ell_1)}}$, separately. Moreover, as already seen before, \eqref{goodequivalenceclassic} implies
$$\forall\;\ell>0:\;\;\;\gamma(\sigma)=\gamma(\omega_{\mathbf{S}^{(\ell)}}),\hspace{10pt}\overline{\gamma}(\sigma)=\overline{\gamma}(\omega_{\mathbf{S}^{(\ell)}}),$$
and similarly for $\tau$.

In view of Corollary \ref{strongrelationlemmaauxprimescor} and Remarks \ref{strongrelationlemmaauxprimesrem} and \ref{closingremark} it follows that $\omega_{\mathbf{S}^{(\ell_1)}}\widehat{\star}\omega_{\mathbf{T}^{(\ell)}}$ can be well-defined for some $\ell_1,\ell>0$, and \cite[Thm. 4.10, Cor. 4.11]{genLegendreconj} can be applied to this weight, even in situations for which this is not clear for $\sigma\widehat{\star}\tau$. This is due to the fact that the matrix parameter allows for having additional ``flexibility''.\vspace{6pt}

Now, by taking into account this comment we want to treat mixed settings and investigate the role of the matrix parameter. First we state and prove a generalization of \cite[Cor. 4.11]{genLegendreconj}:

\begin{proposition}\label{indexprop}
Let $\sigma,\tau\in\hyperlink{omset0}{\mathcal{W}_0}$ be given with associated matrix $\mathcal{M}_{\sigma}:=\{\mathbf{S}^{(\ell)}: \ell>0\}$, $\mathcal{M}_{\tau}:=\{\mathbf{T}^{(\ell)}: \ell>0\}$.

\begin{itemize}
\item[$(i)$] Let $\alpha>0$ and assume that $\sigma(s)=o(s^{1/\alpha})$ as $s\rightarrow+\infty$. When $\gamma(\sigma)>\alpha$, then $$\forall\;\ell,\ell_1>0:\;\;\;\gamma(\omega_{\mathbf{S}^{(\ell_1)}})\le\gamma((((\omega_{\mathbf{S}^{(\ell)}}^{\alpha})^{\star})^{\iota})^{1/\alpha})+\alpha,$$
    and the additional assumption $(\alpha\le)\overline{\gamma}(\sigma)<+\infty$ implies
$$\forall\;\ell,\ell_1>0:\;\;\;\overline{\gamma}((((\omega_{\mathbf{S}^{(\ell)}}^{\alpha})^{\star})^{\iota})^{1/\alpha})+\alpha\le\overline{\gamma}(\omega_{\mathbf{S}^{(\ell_1)}}).$$

\item[$(ii)$] Assume that either $\omega_{\mathbf{S}^{(x)}}\widehat{\star}\omega_{\mathbf{T}^{(y)}}$ is well-defined for some $x,y>0$ or that $\sigma\widehat{\star}\tau$ is well-defined. If then $0<\gamma(\tau)=\overline{\gamma}(\tau)<\gamma(\sigma)\le\overline{\gamma}(\sigma)<+\infty$ is valid, then for all $\ell_j>0$, $j=1,\dots,8$ we get
$$\gamma(\omega_{\mathbf{S}^{(\ell_1)}})\le\gamma(\omega_{\mathbf{S}^{(\ell_2)}}\widehat{\star}\omega_{\mathbf{T}^{(\ell_3)}})+\overline{\gamma}(\omega_{\mathbf{T}^{(\ell_4)}})\le\overline{\gamma}(\omega_{\mathbf{S}^{(\ell_5)}}\widehat{\star}\omega_{\mathbf{T}^{(\ell_6)}})+\gamma(\omega_{\mathbf{T}^{(\ell_7)}})\le\overline{\gamma}(\omega_{\mathbf{S}^{(\ell_8)}}).$$
\end{itemize}
\end{proposition}

\demo{Proof}
$(i)$ We apply $(i)$ in \cite[Cor. 4.11]{genLegendreconj} to the (associated) weights and involve the following facts: The indices are preserved under equivalence and, again by \eqref{goodequivalenceclassic}, relation $\omega(s)=o(s^{1/\alpha})$ as $s\rightarrow+\infty$ holds if and only if $\omega_{\mathbf{S}^{(\ell)}}(s)=o(s^{1/\alpha})$ for some/any $\ell>0$. Moreover, recall \eqref{upperLegendregeneral} and then \cite[Thm. 4.10 $(i)$]{genLegendreconj} together with $\gamma(\sigma)>\alpha$ and $\overline{\gamma}(\id^{1/\alpha})=\alpha<+\infty$ gives $\gamma((((\sigma^{\alpha})^{\star})^{\iota})^{1/\alpha})>0$. Thus \cite[Thm. 4.14 $(i)$]{genLegendreconj} implies $(((\omega_{\mathbf{S}^{(\ell)}}^{\alpha})^{\star})^{\iota})^{1/\alpha}\hyperlink{sim}{\sim}(((\sigma^{\alpha})^{\star})^{\iota})^{1/\alpha}$ for all $\ell>0$.\vspace{6pt}

$(ii)$ By assumption on the indices and \cite[Thm. 4.10 $(i)$]{genLegendreconj} we get $\gamma(\omega_{\mathbf{S}^{(x)}}\widehat{\star}\omega_{\mathbf{T}^{(y)}})>0$ resp. $\gamma(\sigma\widehat{\star}\tau)>0$. Then \cite[Cor. 4.11 $(ii)$]{genLegendreconj}, \cite[Thm. 4.14 $(i)$ \& $(ii)$]{genLegendreconj} and again the fact that both indices are preserved under equivalence imply the conclusion.
\qed\enddemo

Inspired by $\sigma(s)=o(s^{1/\alpha})$ in the previous result we finish this section by involving the non-reflexive relation $\tau\hyperlink{omvartriangle}{\vartriangleleft}\sigma$ and so giving an application of Proposition \ref{strongrelationlemma} for the growth indices $\gamma(\cdot)$ and $\overline{\gamma}(\cdot)$, see also Lemma \ref{strongrelationlemmaauxvar}.

\begin{proposition}\label{indexprop1}
Let $\sigma,\tau\in\hyperlink{omset0}{\mathcal{W}_0}$ be given with associated matrices $\mathcal{M}_{\sigma}:=\{\mathbf{S}^{(\ell)}: \ell>0\}$, $\mathcal{M}_{\tau}:=\{\mathbf{T}^{(\ell)}: \ell>0\}$. Assume that $\tau\hyperlink{omvartriangle}{\vartriangleleft}\sigma$ holds.
\begin{itemize}
\item[$(i)$] If $0\le\overline{\gamma}(\tau)<\gamma(\sigma)\le+\infty$, then for all $\ell,\ell_1,j,j_1>0$:
$$\gamma(\omega_{\mathbf{S}^{(\ell_1)}})\le\gamma(\omega_{\mathbf{S}^{(\ell)}}\widehat{\star}\omega_{\mathbf{T}^{(j)}})+\overline{\gamma}(\omega_{\mathbf{T}^{(j_1)}}).$$

\item[$(ii)$] If $0<\gamma(\tau)\le\overline{\gamma}(\tau)<\gamma(\sigma)\le\overline{\gamma}(\sigma)<+\infty$, then for all $\ell,\ell_1,j,j_1>0$:
$$\overline{\gamma}(\omega_{\mathbf{S}^{(\ell)}}\widehat{\star}\omega_{\mathbf{T}^{(j)}})+\gamma(\omega_{\mathbf{T}^{(j_1)}})\le\overline{\gamma}(\omega_{\mathbf{S}^{(\ell_1)}}).$$
\end{itemize}
\end{proposition}

\demo{Proof}
First, by $(i)\Rightarrow(iv),(vi)$ in Proposition \ref{strongrelationlemma} one has that $\sigma\widehat{\star}\tau$ and $\omega_{\mathbf{S}^{(\ell)}}\widehat{\star}\omega_{\mathbf{T}^{(j)}}$ is well-defined for any $\ell,j>0$.

Then use \cite[Thm. 4.10 $(i)$]{genLegendreconj} and \cite[Thm. 4.14 $(i) \& (ii)$]{genLegendreconj} for part $(i)$, \cite[Thm. 4.10 $(i) \& (ii)$]{genLegendreconj} and \cite[Thm. 4.14 $(i) \& (ii)$]{genLegendreconj} for part $(ii)$. Again recall that $\gamma(\cdot)$, $\overline{\gamma}(\cdot)$ are preserved under equivalence.
\qed\enddemo

\emph{Note:} The analogous comments as the ones given in \cite[Rem. 4.12]{genLegendreconj} also apply to Propositions \ref{indexprop}, \ref{indexprop1}. On the other hand, we cannot apply these propositions to $\sigma=\tau$ since then the assumptions $\sigma\hyperlink{omvartriangle}{\vartriangleleft}\sigma$ and $\overline{\gamma}(\sigma)<\gamma(\sigma)$ fail and $\overline{\gamma}(\sigma)=+\infty$. This observation is consistent with the comments in Remark \ref{closingremark}: For any given weight $\sigma$ one can never expect that $\omega_{\mathbf{S}^{(\ell)}}\widehat{\star}\omega_{\mathbf{S}^{(j)}}$ is a weight function for all arbitrary $\ell,j>0$.

\subsection{On the division of weight matrices}\label{divisionsection}
Analogously as in \eqref{productmatrix} resp. \eqref{productmatrixspecial} for given $\sigma,\tau\in\hyperlink{omset0}{\mathcal{W}_0}$ with associated matrices $\mathcal{M}_{\sigma}:=\{\mathbf{S}^{(\ell)}: \ell>0\}$ and $\mathcal{M}_{\tau}:=\{\mathbf{T}^{(\ell)}: \ell>0\}$ and any $\mathbf{N}\in\RR_{>0}^{\NN}$we set
\begin{equation}\label{quotientmatrix}
\frac{\mathcal{M}_{\sigma}}{\mathcal{M}_{\tau}}:=\left\{\frac{\mathbf{S}^{(\ell)}}{\mathbf{T}^{(\ell)}}: \ell>0\right\},\hspace{15pt}\frac{\mathcal{M}_{\sigma}}{\mathbf{N}}:=\left\{\frac{\mathbf{S}^{(\ell)}}{\mathbf{N}}: \ell>0\right\}.
\end{equation}
A special but useful situation is given when considering $\mathbf{N}\equiv\mathbf{G}^{\alpha}$ (for fixed $\alpha>0$) and hence the analogue of \eqref{Gevreymult}:
\begin{equation}\label{Gevreydiv}
\frac{\mathcal{M}_{\sigma}}{\mathbf{G}^{\alpha}}:=\left\{\frac{\mathbf{S}^{(\ell)}}{\mathbf{G}^{\alpha}}: \ell>0\right\}.
\end{equation}
However, in order to get more information it is natural and desirable to assume that:
\begin{itemize}
\item[$(a)$] Each $\frac{\mathbf{S}^{(\ell)}}{\mathbf{T}^{(\ell)}}$ is at least equivalent to a log-convex sequence, say $\mathbf{Q}^{(\ell)}$.

\item[$(b)$] $\mathcal{M}_{\tau}\vartriangleleft\mathcal{M}_{\sigma}$ holds since then $\lim_{p\rightarrow+\infty}(Q^{(\ell)}_p)^{1/p}=\lim_{p\rightarrow+\infty}\left(\frac{S^{(\ell)}_p}{T^{(\ell)}_p}\right)^{1/p}=+\infty$ is ensured for each $\ell>0$.
\end{itemize}
But $(a)$ is \emph{not clear in general} which is a completely different situation compared with the point-wise product considered in \eqref{productmatrix}, \eqref{productmatrixspecial}.

\subsection{Technical preliminaries on the division of weight matrices and the growth index $\gamma(\mathbf{M})$}\label{technicaldivsection}
When dealing with the matrix from \eqref{Gevreydiv}, then more information concerning requirement $(a)$ from the previous section can be obtained by involving the \emph{growth index for weight sequences $\gamma(\mathbf{M})$} introduced by V. Thilliez in \cite[Sect. 1.3]{Thilliezdivision} and studied in detail in \cite{index}. On the one hand, in the weight function setting one should expect to involve the index $\gamma(\omega)$ for weight functions but, on the other hand, due to the definition of the matrix and the appearance of the Gevrey sequence (with index $\alpha$), the following facts and explanations from \cite{index} seem to be crucial and natural; see also the citations in \cite{index}:

\begin{itemize}
\item[$(a)$] For any $\mathbf{M}\in\RR_{>0}^{\NN}$ and $\beta\ge 0$ the following condition was appearing in \cite[Thm. 3.11 $(v)$]{index}:
\begin{equation}\label{beta3cond}
\exists\;Q\in\NN_{\ge 2}:\;\;\;\liminf_{p\rightarrow+\infty}\frac{\mu_{Qp}}{\mu_p}>Q^{\beta}.
\end{equation}
Indeed, in \cite{index} a slightly different notation involving an index shift for the quotient sequences is used, see \cite[Sect. 3.1]{index}, but this does not effect our considerations: we refer to \cite[Rem. 3.8]{index} and the comments after \cite[Cor. 3.12]{index}.

\item[$(b)$] Let $\mathbf{M}$ be log-convex and such that $\lim_{p\rightarrow+\infty}M_p^{1/p}=+\infty$ and $M_0=1$; i.e. $\mathbf{M}$ is \emph{weight sequence} in the terminology of \cite{index}; see \cite[Sect. 3.1]{index}. Then by \cite[Thm. 3.11 $(v)\Leftrightarrow(viii)$]{index} we have that \eqref{beta3cond} holds if and only if $\gamma(\mathbf{M})>\beta$ and $\gamma(\mathbf{M})$ is denoting the aforementioned growth index by V. Thilliez; we refer to the comments and citations in \cite[Sect. 3.1]{index}.

\item[$(c)$] Indeed, \cite[Thm. 3.11]{index} gathers several equivalent reformulations which can be used for introducing $\gamma(\mathbf{M})$. The original definition by Thilliez in \cite[Sect. 1.3]{Thilliezdivision} is as follows:

    For any $\beta<\gamma(\mathbf{M})$ it is required that $\frac{\mathbf{M}}{\mathbf{G}^{\beta}}$ is equivalent to a log-convex weight sequence $\mathbf{L}$. Indeed, the equivalence is even required w.r.t. the corresponding quotient sequences which is stronger than $\frac{\mathbf{M}}{\mathbf{G}^{\beta}}\hyperlink{approx}{\approx}\mathbf{L}$ in general; see \cite[Thm. 3.11 $(ii)$ \& $(iii)$]{index}.

\item[$(d)$] On the other hand, let $\omega\in\hyperlink{omset0}{\mathcal{W}_0}$ be given with associated matrix $\mathcal{M}_{\omega}:=\{\mathbf{W}^{(\ell)}: \ell>0\}$, then
    \begin{equation}\label{commentdequ}
    \forall\;\ell>0:\;\;\;\gamma(\omega)=\gamma(\omega_{\mathbf{W}^{(\ell)}})\ge\gamma(\mathbf{W}^{(\ell)}).
    \end{equation}
    The equality holds by \eqref{goodequivalenceclassic} and the inequality by \cite[Cor. 4.6 $(i)$]{index}. If $\mathbf{W}^{(\ell)}$ satisfies \hyperlink{mg}{$(\on{mg})$} then all quantities in \eqref{commentdequ} coincide; see \cite[Cor. 4.6 $(iii)$]{index}.

    But in this case the matrix is already constant by $(iii)$ in Section \ref{assomatrixsection} and, in general, the difference between $\gamma(\omega_{\mathbf{M}})$ and $\gamma(\mathbf{M})$ can become (very) large; we refer to \cite[Sect. 5]{index} for an example.
\end{itemize}

Concerning the matrix $\frac{\mathcal{M}_{\sigma}}{\mathbf{G}^{\alpha}}$ from \eqref{Gevreydiv}, in order to ensure $(a)$ in Section \ref{divisionsection} and in view of these explanations, it is crucial to study \eqref{beta3cond} for all $\mathbf{S}^{(\ell)}\in\mathcal{M}_{\sigma}$ and to investigate the role of the matrix parameter $\ell$. This has already been the goal in \cite[Prop. 5.1]{modgrowthstrange}; indeed, we reformulate and complete this statement and improve $(II)$ there.

\begin{proposition}\label{modProp51new}
Let $\omega\in\hyperlink{omset0}{\mathcal{W}_0}$ be given with associated matrix $\mathcal{M}_{\omega}:=\{\mathbf{W}^{(\ell)}: \ell>0\}$. Assume that there exists $\ell>0$ such that $\mathbf{W}^{(\ell)}$ satisfies \eqref{beta3cond} for some $\beta\ge 0$ and $Q\in\NN_{\ge 2}$. Then the following holds:

\begin{itemize}
\item[$(i)$] For all $c\in\NN_{>0}$ the sequences $\mathbf{W}^{(c\ell)}$ satisfy \eqref{beta3cond} for $\beta$ and with the same uniform choice $Q$.

\item[$(ii)$] When $\beta=0$, then for all $c\in\NN_{>0}$ the sequences $\mathbf{W}^{(\ell/c)}$ satisfy \eqref{beta3cond} with $\beta$ and with $Q':=4Q$.

    If $\beta>0$, then for all $c\in\NN_{>0}$ the sequences $\mathbf{W}^{(\ell/c)}$ satisfy \eqref{beta3cond} for all $0<\beta'<\beta$ with $Q':=4Q^n$ and $n\in\NN_{>0}$ is chosen (minimal) such that $n>\frac{2}{\beta/\beta'-1}$.

\item[$(iii)$] For all $a>\ell$ such that $a\neq c\ell$, $c\in\NN_{>0}$, we have that $\mathbf{W}^{(a)}$ satisfies \eqref{beta3cond} for $\beta$ and with the uniform choice $2Q$.

\item[$(iv)$] Let $a<\ell$ such that $a\neq\frac{\ell}{c}$, $c\in\NN_{>0}$. Then $\mathbf{W}^{(a)}$ satisfies \eqref{beta3cond} for $\beta=0$ and with the uniform choice $8Q$ and, if $\beta>0$, then $\mathbf{W}^{(a)}$ satisfies \eqref{beta3cond} for all $0<\beta'<\beta$ with $Q':=8Q^n$ and $n\in\NN_{>0}$ is chosen as in $(ii)$.
\end{itemize}
Therefore, one obtains
\begin{equation}\label{modProp51newequ}
\forall\;\ell_1>0:\;\;\;\gamma(\mathbf{W}^{(\ell_1)})=\gamma(\mathbf{W}^{(\ell)}).
\end{equation}
\end{proposition}

\demo{Proof}
$(i)$ For any $\ell>0$ recall that $\vartheta^{(\ell)}=(\vartheta^{(\ell)}_p)_p$ denotes the corresponding quotient sequence; i.e. $\vartheta^{(\ell)}_p:=\frac{W^{(\ell)}_p}{W^{(\ell)}_{p-1}}$, $p\in\NN_{>0}$, and $\vartheta^{(\ell)}_0:=1$. Then for all $c,p\in\NN_{>0}$ we obtain, see also \cite[$(4.9)$]{modgrowthstrange},
\begin{equation}\label{modstrange49}
\vartheta^{(c\ell)}_p=\frac{W^{(c\ell)}_p}{W^{(c\ell)}_{p-1}}=\left(\frac{W^{(\ell)}_{cp}}{W^{(\ell)}_{c(p-1)}}\right)^{1/c}=(\vartheta^{(\ell)}_{c(p-1)+1}\cdots\vartheta^{(\ell)}_{cp})^{1/c},
\end{equation}
and for $p=0$ we have clearly the equality $1=1$.\vspace{6pt}

Let $\ell>0$ be such that $\mathbf{W}^{(\ell)}$ has \eqref{beta3cond} with $\beta\ge 0$ and $Q\in\NN_{\ge 2}$; i.e. $\liminf_{p\rightarrow\infty}\frac{\vartheta^{(\ell)}_{Qp}}{\vartheta^{(\ell)}_p}>Q^{\beta}$ holds. Then, by \eqref{modstrange49}, one infers
$$\frac{\vartheta^{(c\ell)}_{Qp}}{\vartheta^{(c\ell)}_p}=\left(\frac{\vartheta^{(\ell)}_{Qcp-c+1}\cdots\vartheta^{(\ell)}_{Qcp}}{\vartheta^{(\ell)}_{cp-c+1}\cdots\vartheta^{(\ell)}_{cp}}\right)^{1/c},$$
and for all $1\le i\le c$ we have
\begin{equation*}\label{strongnonquasilemmaequ}
Qcp-c+i\ge Q(cp-c+i)\Leftrightarrow Qcp-c+i\ge Qcp-Qc+Qi\Leftrightarrow Q(c-i)\ge(c-i).
\end{equation*}
By this and log-convexity for $\mathbf{W}^{(\ell)}$ we obtain
$$\liminf_{p\rightarrow\infty}\frac{\vartheta^{(c\ell)}_{Qp}}{\vartheta^{(c\ell)}_p}\ge\prod_{i=1}^c\left(\liminf_{p\rightarrow\infty}\frac{\vartheta^{(\ell)}_{Qcp-c+i}}{\vartheta^{(\ell)}_{cp-c+i}}\right)^{1/c}\ge\prod_{i=1}^c\left(\liminf_{p\rightarrow\infty}\frac{\vartheta^{(\ell)}_{Q(cp-c+i)}}{\vartheta^{(\ell)}_{cp-c+i}}\right)^{1/c}>Q^{\beta},$$ verifying for any $\mathbf{W}^{(c\ell)}$ property \eqref{beta3cond} with $\beta\ge 0$ and the same $Q$.\vspace{6pt}

$(ii)$ By inspecting the proof of \cite[$(4.10)$]{modgrowthstrange} we have the following equality \cite[$(5.2)$]{modgrowthstrange}:
$$\vartheta^{(\ell/c)}_{cp}=\frac{W^{(\ell/c)}_{cp}}{W^{(\ell/c)}_{cp-1}}=\frac{W^{(\ell/c)}_{cp}}{W^{(\ell/c)}_{cp-c}}\frac{1}{\vartheta^{(\ell/c)}_{cp-1}}\frac{1}{\vartheta^{(\ell/c)}_{cp-2}}\cdots\frac{1}{\vartheta^{(\ell/c)}_{cp-c+1}}=\left(\frac{W^{(\ell)}_{p}}{W^{(\ell)}_{p-1}}\right)^c\frac{1}{\vartheta^{(\ell/c)}_{cp-1}}\frac{1}{\vartheta^{(\ell/c)}_{cp-2}}\cdots\frac{1}{\vartheta^{(\ell/c)}_{cp-c+1}},$$
and so
\begin{equation}\label{modstrange52}
\forall\;\ell>0\;\forall\;c\in\NN_{>0}\;\forall\;p\in\NN_{>0}:\;\;\;\vartheta^{(\ell)}_p=(\vartheta^{(\ell/c)}_{cp}\vartheta^{(\ell/c)}_{cp-1}\cdots\vartheta^{(\ell/c)}_{cp-c+1})^{1/c}.
\end{equation}
Note that \eqref{modstrange52} is analogous to \eqref{modstrange49} above. Now let $Q_1\in\NN_{>0}$, then by involving log-convexity for $\mathbf{W}^{(\ell/c)}$ and since $Q_1cp\le 2Q_1c(p-1)\Leftrightarrow 2\le p$ we get for all $\ell>0$, for all $c,Q_1\in\NN_{>0}$ and $p\ge 2$:
\begin{equation}\label{strongnonquasilemmaequ1}
\frac{\vartheta^{(\ell)}_{Q_1p}}{\vartheta^{(\ell)}_p}=\left(\frac{\vartheta^{(\ell/c)}_{Q_1cp}\vartheta^{(\ell/c)}_{Q_1cp-1}\cdots\vartheta^{(\ell/c)}_{Q_1cp-c+1}}{\vartheta^{(\ell/c)}_{cp}\vartheta^{(\ell/c)}_{cp-1}\cdots\vartheta^{(\ell/c)}_{cp-c+1}}\right)^{1/c}\le\frac{\vartheta^{(\ell/c)}_{Q_1cp}}{\vartheta^{(\ell/c)}_{cp-c+1}}\le\frac{\vartheta^{(\ell/c)}_{2Q_1c(p-1)}}{\vartheta^{(\ell/c)}_{c(p-1)}}.
\end{equation}
By assumption $\liminf_{p\rightarrow+\infty}\frac{\vartheta^{(\ell)}_{Qp}}{\vartheta^{(\ell)}_p}>Q^{\beta}$ holds for $\beta\ge 0$ and some choice $Q\in\NN_{\ge 2}$. Then we distinguish:

\emph{Case 1 - $\beta=0$} Then \eqref{strongnonquasilemmaequ1} applied to $Q_1=Q$ gives the estimate in \eqref{beta3cond} for $\mathbf{W}^{(\ell/c)}$ when choosing $Q':=2Q$ and restricting in the $\liminf$ to all $q\in\NN$ satisfying $q=c(p-1)$ for some $p\ge 2$.

Assume now that $q\in\NN_{>0}$ is such that $c(p-1)<q<cp$ for some $p\ge 2$. By log-convexity for $\mathbf{W}^{(\ell/c)}$ and \eqref{strongnonquasilemmaequ1} applied to $p':=2(p-1)+1=2p-1$, thus $p'\ge 2\Leftrightarrow p\ge\frac{3}{2}$ is valid with $p'-1=2(p-1)$, it follows that
\begin{equation}\label{strongnonquasilemmaequ2}
\frac{\vartheta^{(\ell)}_{Q_1(2p-1)}}{\vartheta^{(\ell)}_{2p-1}}\le\frac{\vartheta^{(\ell/c)}_{4Q_1c(p-1)}}{\vartheta^{(\ell/c)}_{2c(p-1)}}\le\frac{\vartheta^{(\ell/c)}_{4Q_1c(p-1)}}{\vartheta^{(\ell/c)}_{cp}}\le\frac{\vartheta^{(\ell/c)}_{4Q_1q}}{\vartheta^{(\ell/c)}_{q}},
\end{equation}
since $2c(p-1)\ge cp\Leftrightarrow p\ge 2$ and note that \eqref{strongnonquasilemmaequ2} holds for any $Q_1\in\NN_{\ge 2}$. Thus we have shown
$$\liminf_{p\rightarrow+\infty}\frac{\vartheta^{(\ell/c)}_{4Qp}}{\vartheta^{(\ell/c)}_p}>1;$$
i.e. \eqref{beta3cond} for $\mathbf{W}^{(\ell/c)}$ with $\beta=0$ when taking $Q':=4Q$.\vspace{6pt}

\emph{Case 2 - $\beta>0$} Let $0<\beta'<\beta$ be fixed and iteration of \eqref{beta3cond} gives that there exists some $Q\in\NN_{\ge 2}$ such that $\liminf_{p\rightarrow+\infty}\frac{\vartheta^{(\ell)}_{Q^np}}{\vartheta^{(\ell)}_p}>Q^{n\beta}$ for all $n\in\NN_{>0}$. Take $n\in\NN_{>0}$ (minimal) such that $n>\frac{2\beta'}{\beta-\beta'}=\frac{2}{\beta/\beta'-1}$ and so $2^{2\beta'}\le Q^{n(\beta-\beta')}$ because $Q\in\NN_{\ge 2}$. Then this choice and \eqref{strongnonquasilemmaequ1} applied to $Q_1=Q^n$ give
\begin{equation}\label{strongnonquasilemmaequ3}
(2Q^n)^{\beta'}<(4Q^n)^{\beta'}\le Q^{n\beta}<\liminf_{p\rightarrow+\infty}\frac{\vartheta^{(\ell)}_{Q^np}}{\vartheta^{(\ell)}_p}\le\liminf_{p\rightarrow+\infty}\frac{\vartheta^{(\ell/c)}_{2Q^nc(p-1)}}{\vartheta^{(\ell/c)}_{c(p-1)}},
\end{equation}
hence the estimate in \eqref{beta3cond} is verified with $\beta'$ and $Q':=2Q^n$ and again when restricting in the $\liminf$ to all $q\in\NN$ satisfying $q=c(p-1)$ for some $p\ge 2$.

Assume now that $q\in\NN_{>0}$ is such that $c(p-1)<q<cp$ for some $p\ge 2$, then combining \eqref{strongnonquasilemmaequ2} (with $Q_1:=Q^n$) and \eqref{strongnonquasilemmaequ3} yields
$$(4Q^n)^{\beta'}\le Q^{n\beta}<\liminf_{p\rightarrow+\infty}\frac{\vartheta^{(\ell)}_{Q^n(2p-1)}}{\vartheta^{(\ell)}_{2p-1}}\le\frac{\vartheta^{(\ell/c)}_{4Q^nc(p-1)}}{\vartheta^{(\ell/c)}_{cp}}\le\frac{\vartheta^{(\ell/c)}_{4Q^nq}}{\vartheta^{(\ell/c)}_{q}}.$$
This verifies the desired estimate in \eqref{beta3cond} with $\beta'$ and $Q':=4Q^n$. Summarizing, $\mathbf{W}^{(\ell/c)}$ satisfies \eqref{beta3cond} with $\beta'$ and $Q':=4Q^n$. Note that $n\rightarrow+\infty$ as $\beta'\rightarrow\beta$.\vspace{6pt}

$(iii)$ First, \eqref{modstrange49} and the log-convexity for $\mathbf{W}^{(\ell)}$ give
\begin{equation}\label{modstrange49new}
\forall\;\ell>0\;\forall\;c\in\NN_{>0}\;\forall\;p\in\NN:\;\;\;\vartheta^{(c\ell)}_p\le\vartheta^{(\ell)}_{cp}.
\end{equation}
Let $a>\ell$ such that $a\neq\ell c$, $c\in\NN_{>0}$. Then $\ell c<a<\ell(c+1)$ for some $c\in\NN_{>0}$ holds and in this case $2a>2\ell c$ and $2\ell c\ge\ell(c+1)\Leftrightarrow c\ge 1$. By \eqref{modstrange49new} and \eqref{quotientorderequ} (i.e. the point-wise order of the sequences of quotients) one has for any $Q'\in\NN_{\ge 2}$ and $p\in\NN$:
$$\frac{\vartheta^{(a)}_{2Q'p}}{\vartheta^{(a)}_p}\ge\frac{\vartheta^{(2a)}_{Q'p}}{\vartheta^{(\ell(c+1))}_p}\ge\frac{\vartheta^{(\ell(c+1))}_{Q'p}}{\vartheta^{(\ell(c+1))}_p}.$$
In view of $(i)$ applied to $\mathbf{W}^{(\ell(c+1))}$ and when choosing $Q':=Q$ in this estimate we have verified that $\liminf_{p\rightarrow+\infty}\frac{\vartheta^{(a)}_{2Qp}}{\vartheta^{(a)}_p}\ge\liminf_{p\rightarrow+\infty}\frac{\vartheta^{(\ell(c+1))}_{Qp}}{\vartheta^{(\ell(c+1))}_p}$; consequently $\mathbf{W}^{(a)}$ satisfies \eqref{beta3cond} with $\beta$ and with the choice $2Q$.\vspace{6pt}

$(iv)$ Analogously, let $0<a<\ell$ such that $a\neq\frac{\ell}{c}$, $c\in\NN_{>0}$. So $\frac{\ell}{c+1}<a<\frac{\ell}{c}$ for some $c\in\NN_{>0}$ holds and $2a>\frac{2\ell}{c+1}$ and $\frac{2\ell}{c+1}\ge\frac{\ell}{c}\Leftrightarrow c\ge 1$. Involving again \eqref{modstrange49new} and the point-wise order of the sequences of quotients one has for any $Q'\in\NN_{\ge 2}$ and $p\in\NN$:
$$\frac{\vartheta^{(a)}_{2Q'p}}{\vartheta^{(a)}_p}\ge\frac{\vartheta^{(2a)}_{Q'p}}{\vartheta^{(\ell/c)}_p}\ge\frac{\vartheta^{(\ell/c)}_{Q'p}}{\vartheta^{(\ell/c)}_p}.$$
In view of $(ii)$ we have to distinguish: If $\beta=0$, then the above estimate applied to $\mathbf{W}^{(\ell/c)}$ and when choosing $Q':=4Q$ we have verified that $\liminf_{p\rightarrow+\infty}\frac{\vartheta^{(a)}_{8Qp}}{\vartheta^{(a)}_p}>1$ is true; i.e. $\mathbf{W}^{(a)}$ satisfies \eqref{beta3cond} with $\beta$ and with the choice $8Q$.

Finally, if $\beta>0$, then fix $0<\beta'<\beta$ and similarly the above verifies $\liminf_{p\rightarrow+\infty}\frac{\vartheta^{(a)}_{2Q'p}}{\vartheta^{(a)}_p}\ge\liminf_{p\rightarrow+\infty}\frac{\vartheta^{(\ell/c)}_{Q'p}}{\vartheta^{(\ell/c)}_p}$. Then, $(ii)$ applied to $\mathbf{W}^{(\ell/c)}$ yields that $\mathbf{W}^{(a)}$ satisfies \eqref{beta3cond} for $\beta'$ with the choice $8Q^n$ and $n\in\NN_{>0}$ is like in $(ii)$ chosen minimal to ensure $n>\frac{2}{\beta/\beta'-1}$.\vspace{6pt}

Concerning \eqref{modProp51newequ} note the following: $(i)-(iv)$ together verify that $\gamma(\mathbf{W}^{(\ell_1)})\ge\gamma(\mathbf{W}^{(\ell)})$ for all $\ell_1>0$. Then fix an arbitrary $\ell_1>0$ and repeat the above arguments when $\ell$ is replaced by $\ell_1$. Therefore, $\gamma(\mathbf{W}^{(\ell)})\ge\gamma(\mathbf{W}^{(\ell_1)})$ follows as well and hence equality. Since $\ell_1$ was chosen arbitrarily we are done.
\qed\enddemo

\begin{remark}
\emph{Proposition \ref{modProp51new}, in particular the identity \eqref{modProp51newequ}, answers the conjecture of the author stated in \cite[Rem. 5.2 $(ii)$]{modgrowthstrange}. On the other hand, as mentioned in this observation, in order to treat Roumieu-type spaces it suffices to consider the (sub-)matrix given by the sequences in $(i)$ in Proposition \ref{modProp51new} and for the Beurling-type spaces it suffices to treat the weights studied in $(ii)$ there (when $\ell>0$ is arbitrary but fixed): These weight (sub-)matrices are R- resp. B-equivalent to given $\mathcal{M}_{\omega}$ and therefore contain ``enough information'' to define and study the corresponding weighted function spaces.}

\emph{Using the complete information from this result we are able to give a more direct proof of the main result Theorem \ref{mainweighfctcorinv} below. If one would only use the (sufficient) information from $(i)$ and $(ii)$ then one has to restrict to a convenient (sub-)matrix and only consider $c,c^{-1}\in\NN_{>0}$ instead of $c>0$ arbitrary but which requires some more technical comments.}
\end{remark}

\subsection{Main result on the division of weight matrices}
Using the preparation from the previous sections we prove the following result which is inverse to Corollary \ref{mainweighfctcor}.

\begin{theorem}\label{mainweighfctcorinv}
Let $\omega\in\hyperlink{omset0}{\mathcal{W}_0}$ be given with associated matrix $\mathcal{M}_{\omega}:=\{\mathbf{W}^{(\ell)}: \ell>0\}$ and let $\alpha>0$. Assume the following:
\begin{itemize}
\item[$(i)$] There exists $\ell_0>0$ such that $\mathbf{W}^{(\ell_0)}$ satisfies $\gamma(\mathbf{W}^{(\ell_0)})>\alpha$.

\item[$(ii)$] $\omega(t)=o(t^{1/\alpha})$ as $t\rightarrow+\infty$ holds.
\end{itemize}
Then, for the matrix $\frac{\mathcal{M}_{\omega}}{\mathbf{G}^{\alpha}}$ from \eqref{Gevreydiv} it follows that as l.c.v.s.
\begin{equation}\label{mainweighfctcorinvequ}
\forall\;c,\ell>0:\;\;\;\mathfrak{F}_{\left[\frac{\mathcal{M}_{\omega}}{\mathbf{G}^{\alpha}}\right]}=\mathfrak{F}_{[\omega_{\mathbf{W}^{(c\ell_0)}/\mathbf{G}^{\alpha}}]}=\mathfrak{F}_{[\omega_{\mathbf{W}^{(\ell)}}\widehat{\star}\omega_{\mathbf{G}^{\alpha}}]}=\mathfrak{F}_{[\omega\widehat{\star}\omega_{\mathbf{G}^{\alpha}}]}=\mathfrak{F}_{[\omega\widehat{\star}\id^{1/\alpha}]}=\mathfrak{F}_{[(((\omega^{\alpha})^{\star})^{\iota})^{1/\alpha}]},
\end{equation}
with $\mathfrak{F}\in\{\mathcal{E}, \mathcal{B}, \mathcal{A}, \mathcal{S}, \Lambda, \mathcal{F}\}$.
\end{theorem}

\demo{Proof}
The proof is technical and we divide it into different steps. \eqref{mainweighfctcorinvequ} applies to any symbol/functor under consideration by the analogous definitions of the spaces involving the same weighted semi-norms and by the fact that all forthcoming arguments exclusively deal with weights and their growth relations.\vspace{6pt}

\emph{Step I - On the growth indices} By Proposition \ref{modProp51new} and assumption $(i)$ we see that
\begin{equation}\label{stepiequ0}
\forall\;c>0:\;\;\;\gamma(\mathbf{W}^{(c\ell_0)})=\gamma(\mathbf{W}^{(\ell_0)})>\alpha>0,
\end{equation}
and \eqref{commentdequ} in comment $(d)$ in Section \ref{technicaldivsection} implies
\begin{equation}\label{stepiequ}
\forall\;\ell>0:\;\;\;\gamma(\omega)=\gamma(\omega_{\mathbf{W}^{(\ell)}})=\gamma(\omega_{\mathbf{W}^{(\ell_0)}})\ge\gamma(\mathbf{W}^{(\ell_0)})>\alpha.
\end{equation}
In particular, $\gamma(\omega),\gamma(\omega_{\mathbf{W}^{(\ell)}})>0$ for all $\ell>0$ and so all these weight functions satisfy \hyperlink{om1}{$(\omega_1)$} by taking into account \cite[Thm. 2.11, Rem. 2.12, Cor. 2.14]{index}.\vspace{6pt}

\emph{Step II - On the crucial quotient matrix} By \eqref{stepiequ0} one has $\gamma(\mathbf{W}^{(c\ell)})=\gamma(\mathbf{W}^{(\ell)})>\alpha$ for all $c>0$, thus by recalling the original definition of this growth index by V. Thilliez (see comment $(c)$ in Section \ref{technicaldivsection}) it follows that
\begin{equation}\label{stepiiiequ}
\forall\;c>0\;\exists\;\mathbf{Q}^{(c)}\in\RR_{>0}^{\NN}:\;\;\;\frac{\mathbf{W}^{(c\ell_0)}}{\mathbf{G}^{\alpha}}\hyperlink{approx}{\approx}\mathbf{Q}^{(c)},
\end{equation}
and $\mathbf{Q}^{(c)}$ is log-convex. Consider the set
\begin{equation}\label{stepiiiequ1}
\mathcal{Q}:=\{\mathbf{Q}^{(c)}: c>0\}.
\end{equation}
Remark that Corollary \ref{strongrelationlemmacor} and assumption $(ii)$ yield $\mathbf{G}^{\alpha}\hyperlink{mtriangle}{\vartriangleleft}\mathbf{W}^{(\ell)}$ for any $\ell>0$ and so $$\forall\;c>0:\;\;\;\mathbf{Q}^{(c)}_{\iota}=\lim_{p\rightarrow+\infty}(Q^{(c)}_p)^{1/p}=\lim_{p\rightarrow+\infty}\left(\frac{\mathbf{W}^{(c\ell_0)}_p}{p!^{\alpha}}\right)^{1/p}=+\infty.$$
Moreover, by \eqref{stepiiiequ} it holds that $\frac{\mathcal{M}_{\omega}}{\mathbf{G}^{\alpha}}\{\approx\}\mathcal{Q}$ and $\frac{\mathcal{M}_{\omega}}{\mathbf{G}^{\alpha}}(\approx)\mathcal{Q}$. Therefore, both matrices define the same weighted (ultradifferentiable) classes of the particular type and $\mathcal{Q}$ satisfies \hyperlink{R-mg}{$(\mathcal{M}_{\{\on{mg}\}})$} and \hyperlink{B-mg}{$(\mathcal{M}_{(\on{mg})})$}: $\frac{\mathcal{M}_{\omega}}{\mathbf{G}^{\alpha}}$ clearly has both  \hyperlink{R-mg}{$(\mathcal{M}_{\{\on{mg}\}})$} and \hyperlink{B-mg}{$(\mathcal{M}_{(\on{mg})})$} by \eqref{newmoderategrowth} and by taking into account the fact that both properties are preserved when dividing the sequences of the matrix point-wise by a sequence $\mathbf{L}$ satisfying $L_pL_q\le A^{p+q}L_{p+q}$ for some $A\ge 1$ and all $p,q\in\NN$. Finally, $\mathcal{M}_{[\on{mg}]}$ is preserved under relation $[\approx]$.\vspace{6pt}

Concerning $\mathcal{M}_{[\on{L}]}$ note that in general this property might not be preserved under $[\approx]$ and so it is not clear for $\mathcal{Q}$. On the one hand, by the technical result \cite[Lemma 2.5]{optimalBeurling} and the comment just before this statement, $\mathcal{Q}$ can be replaced by an R-equivalent log-convex matrix $\mathcal{Q}^R$ satisfying \hyperlink{R-L}{$(\mathcal{M}_{\{\on{L}\}})$} and by an B-equivalent log-convex matrix $\mathcal{Q}^B$ satisfying \hyperlink{B-L}{$(\mathcal{M}_{(\on{L})})$}. But, on the other hand, $\mathcal{M}_{[\on{L}]}$ is not required necessarily for the arguments since we use directly the characterizations \cite[Lemma 3.1, Thm. 3.2 \& 5.4]{equalitymixedOregular} and not the results from \cite{testfunctioncharacterization}; see \emph{Step VI} below.\vspace{6pt}

Note that by \eqref{stepiiiequ} in general one only has $\mathbf{Q}^{(c)}\hyperlink{preceq}{\preceq}\mathbf{Q}^{(d)}$ for $c\le d$ and hence $\mathcal{Q}$ is formally not a weight matrix according to Section \ref{matrixsection} since the point-wise order might fail. Moreover, $Q^{(c)}_0\le Q^{(c)}_1$ (normalization) is not clear. However, these formal failures are not effecting the forthcoming arguments and, on the other hand, one can achieve $\widetilde{Q}^{(c)}_0\le\widetilde{Q}^{(c)}_1$ and $\widetilde{\mathbf{Q}}^{(c)}\le\widetilde{\mathbf{Q}}^{(d)}$ by switching to an equivalent matrix $\widetilde{\mathcal{Q}}:=\{\widetilde{\mathbf{Q}}^{(c)}: c>0\}$; see e.g. also \cite[Rem. 3.1 $(i)$]{mixedramisurj}. Therefore, w.l.o.g. one can assume that $\mathcal{Q}$ is \emph{standard log-convex} according to Section \ref{matrixsection} and in order to lighten notation we write again $\mathcal{Q}$ instead of $\widetilde{\mathcal{Q}}$.\vspace{6pt}

\emph{Step III - On the auxiliary matrix} Introduce the matrix
\begin{equation}\label{stepivequnew}
\widetilde{\mathcal{M}}_{\omega}:=\{\widetilde{\mathbf{W}}^{(c)}: c>0\},\hspace{15pt}
\widetilde{\mathbf{W}}^{(c)}:=\mathbf{G}^{\alpha}\cdot\mathbf{Q}^{(c)},\;\;\;c>0.
\end{equation}
First, by \eqref{stepiiiequ} in \emph{Step II} it follows that
\begin{equation}\label{stepivequ0}
\forall\;c>0:\;\;\;\widetilde{\mathbf{W}}^{(c)}\hyperlink{approx}{\approx}\mathbf{W}^{(c\ell_0)},
\end{equation}
and so, in particular, $\widetilde{\mathcal{M}}_{\omega}$ and $\mathcal{M}_{\omega}$ are both R- and B-equivalent. Consequently, $\widetilde{\mathcal{M}}_{\omega}$ satisfies \hyperlink{R-mg}{$(\mathcal{M}_{\{\on{mg}\}})$} and \hyperlink{B-mg}{$(\mathcal{M}_{(\on{mg})})$}.\vspace{6pt}

Moreover, let us show that each $\omega_{\widetilde{\mathbf{W}}^{(c)}}$ has \hyperlink{om1}{$(\omega_1)$} and
\begin{equation}\label{stepivequ}
\forall\;c>0\;\forall\;\ell>0:\;\;\;\omega_{\widetilde{\mathbf{W}}^{(c)}}\hyperlink{sim}{\sim}\omega_{\mathbf{W}^{(\ell)}}\hyperlink{sim}{\sim}\omega.
\end{equation}
Indeed, first note that \eqref{stepivequ} implies \hyperlink{om1}{$(\omega_1)$} for each $\omega_{\widetilde{\mathbf{W}}^{(c)}}$ since this property is preserved under equivalence and because \hyperlink{om1}{$(\omega_1)$} holds for $\omega$ and for each $\omega_{\mathbf{W}^{(\ell)}}$ as verified in \emph{Step I.} Then recall that the second equivalence in \eqref{stepivequ} is valid by \eqref{goodequivalenceclassic}. Therefore, it remains to show the first equivalence in \eqref{stepivequ} and for this let us now prove the following abstract\vspace{6pt}

\emph{Claim: Let $\mathbf{M},\mathbf{N}\in\RR_{>0}^{\NN}$ be given such that $\mathbf{M}_{\iota}=+\infty=\mathbf{N}_{\iota}$. Assume that $\mathbf{M}\hyperlink{preceq}{\preceq}\mathbf{N}$ and that either $\omega_{\mathbf{M}}$ or $\omega_{\mathbf{N}}$ satisfies \hyperlink{om1}{$(\omega_1)$}. Then $\omega_{\mathbf{N}}(t)=O(\omega_{\mathbf{M}}(t))$ as $t\rightarrow+\infty$ holds.}\vspace{6pt}

\emph{Proof of the claim:} By definition of relation $\mathbf{M}\hyperlink{preceq}{\preceq}\mathbf{N}$ one has
$$\exists\;C\ge 1\;\forall\;t\ge 0:\;\;\;\omega_{\mathbf{N}}(t)\le\omega_{\mathbf{M}}(Ct)+C,$$
more precisely
$$\forall\;t\ge 0:\;\;\;\omega_{\mathbf{N}}(t)\le\omega_{\mathbf{M}}(tC_{\mathbf{M}\preceq\mathbf{N}}),$$
with
$$C_{\mathbf{M}\preceq\mathbf{N}}:=\sup_{p\in\NN_{>0}}\left(\frac{M_p/M_0}{N_p/N_0}\right)^{1/p}\in(0,+\infty);$$
see \cite[Sect. 2.2, Sect. 2.5 $(2.20)$]{genLegendreconj}.

When $\omega_{\mathbf{M}}$ satisfies \hyperlink{om1}{$(\omega_1)$}, then choose $n\in\NN_{>0}$ minimal such that $2^n\ge C$ and so by iteration of this property, for some $L\ge 1$ and all $t\ge 0$:
$$\omega_{\mathbf{N}}(t)\le\omega_{\mathbf{M}}(Ct)+C\le\omega_{\mathbf{M}}(2^nt)+C\le L\omega_{\mathbf{M}}(t)+L+C.$$
Similarly, when $\omega_{\mathbf{N}}$ satisfies \hyperlink{om1}{$(\omega_1)$} we estimate as follows for all $t\ge 0$:
$$\omega_{\mathbf{N}}(t)\le L\omega_{\mathbf{N}}(2^{-n}t)+L\le L\omega_{\mathbf{N}}(C^{-1}t)+L\le L\omega_{\mathbf{M}}(t)+CL+L,$$
which ends the proof of the claim.\vspace{6pt}

Now, this claim together with \eqref{stepivequ0} and the fact that each $\omega_{\mathbf{W}^{(\ell)}}$ satisfies \hyperlink{om1}{$(\omega_1)$} as just mentioned before verifies the first equivalence in \eqref{stepivequ} and finishes this part.\vspace{6pt}

\emph{Step IV - Involving the generalized upper Legendre conjugate} First, each $\omega_{\mathbf{W}^{(\ell)}}\widehat{\star}\omega_{\mathbf{G}^{\alpha}}$ is a weight function because $\omega_{\mathbf{G}^{\alpha}}(0)=0$ and since it is well-defined by Corollary \ref{strongrelationlemmacor} and assumption $\omega(t)=o(t^{1/\alpha})$ as $t\rightarrow+\infty$. Moreover, by this result also $\omega\widehat{\star}\omega_{\mathbf{G}^{\alpha}}$ and $\omega\widehat{\star}\id^{1/\alpha}$ is a weight function.

Next, \eqref{stepiequ} gives $\gamma(\omega)=\gamma(\omega_{\mathbf{W}^{(\ell)}})>\alpha>0$ for all $\ell>0$ and $\alpha=\overline{\gamma}(\id^{1/\alpha})=\overline{\gamma}(\omega_{\mathbf{G}^{\alpha}})$ is valid; see e.g. \cite[Rem. 2.4, Ex. 2.9]{genLegendreconj}. Then, by \cite[Thm. 4.10 $(i)$]{genLegendreconj}, we get that
$$\forall\;\ell>0:\;\;\;0<\gamma(\omega_{\mathbf{W}^{(\ell)}})-\alpha\le\gamma(\omega_{\mathbf{W}^{(\ell)}}\widehat{\star}\omega_{\mathbf{G}^{\alpha}});$$
thus each $\omega_{\mathbf{W}^{(\ell)}}\widehat{\star}\omega_{\mathbf{G}^{\alpha}}$ satisfies \hyperlink{om1}{$(\omega_1)$}.

\cite[Thm. 4.14 $(i)$]{genLegendreconj} applied to $\sigma=\omega_{\mathbf{W}^{(\ell)}}$, $\sigma_1=\omega$, $\tau=\omega_{\mathbf{G}^{\alpha}}$ and \cite[Thm. 4.14 $(ii)$]{genLegendreconj} applied to $\sigma=\omega$, $\tau=\omega_{\mathbf{G}^{\alpha}}$ and $\tau_1=\id^{1/\alpha}$ give
\begin{equation}\label{stepvequ}
\forall\;\ell>0:\;\;\;\omega_{\mathbf{W}^{(\ell)}}\widehat{\star}\omega_{\mathbf{G}^{\alpha}}\hyperlink{sim}{\sim}\omega\widehat{\star}\omega_{\mathbf{G}^{\alpha}}\hyperlink{sim}{\sim}\omega\widehat{\star}\id^{1/\alpha}.
\end{equation}
Note that, as seen before $0<\alpha=\overline{\gamma}(\omega_{\mathbf{G}^{\alpha}})=\overline{\gamma}(\id^{1/\alpha})<\gamma(\omega)\le+\infty$, and hence one can use both parts in \cite[Thm. 4.14]{genLegendreconj}.

By \emph{Step III} the same conclusions are valid when considering $\omega_{\widetilde{\mathbf{W}}^{(c)}}$ instead of $\omega_{\mathbf{W}^{(\ell)}}$; more precisely in view of \eqref{stepivequ} one can extend \eqref{stepvequ} and obtain
\begin{equation}\label{stepvequ1}
\forall\;c>0\;\forall\;\ell>0:\;\;\;\omega_{\widetilde{\mathbf{W}}^{(c)}}\widehat{\star}\omega_{\mathbf{G}^{\alpha}}\hyperlink{sim}{\sim}\omega_{\mathbf{W}^{(\ell)}}\widehat{\star}\omega_{\mathbf{G}^{\alpha}}\hyperlink{sim}{\sim}\omega\widehat{\star}\omega_{\mathbf{G}^{\alpha}}\hyperlink{sim}{\sim}\omega\widehat{\star}\id^{1/\alpha}.
\end{equation}
In particular, this implies that each $\omega_{\widetilde{\mathbf{W}}^{(c)}}\widehat{\star}\omega_{\mathbf{G}^{\alpha}}$ satisfies  \hyperlink{om1}{$(\omega_1)$}, too.\vspace{6pt}

\emph{Step V - Quotient matrix $\mathcal{Q}$ and the generalized upper Legendre conjugate} Recall that by \eqref{stepivequnew} in \emph{Step III} the equation $\widetilde{\mathbf{W}}^{(c)}=\mathbf{G}^{\alpha}\cdot\mathbf{Q}^{(c)}$ holds for all $c>0$ and, as seen in \emph{Step II,} each $\mathbf{Q}^{(c)}$ is log-convex and $\mathbf{Q}^{(c)}_{\iota}=+\infty$. Therefore, we can apply the main result \cite[Thm. 5.11]{genLegendreconj} to $\mathbf{M}\equiv\widetilde{\mathbf{W}}^{(c)}$, $\mathbf{N}\equiv\mathbf{G}^{\alpha}$ and get:
\begin{equation}\label{stepviequ}
\forall\;c>0\;\forall\;t\in[0,+\infty):\;\;\;\omega_{\widetilde{\mathbf{W}}^{(c)}}\widehat{\star}\omega_{\mathbf{G}^{\alpha}}(t)=\omega_{\mathbf{Q}^{(c)}}(t).
\end{equation}
Using this identity, by \eqref{stepvequ1} in \emph{Step IV} it holds that $\omega_{\mathbf{Q}^{(c)}}$ satisfies \hyperlink{om1}{$(\omega_1)$}, too. Then recall that $\frac{\mathbf{W}^{(c\ell_0)}}{\mathbf{G}^{\alpha}}\hyperlink{approx}{\approx}\mathbf{Q}^{(c)}$, see \emph{Step II,} and so the claim shown in \emph{Step III} gives
\begin{equation}\label{stepviequ0}
\forall\;c>0:\;\;\;\omega_{\mathbf{Q}^{(c)}}\hyperlink{sim}{\sim}\omega_{\frac{\mathbf{W}^{(c\ell_0)}}{\mathbf{G}^{\alpha}}}.
\end{equation}
Combining \eqref{stepvequ1}, \eqref{stepviequ} and \eqref{stepviequ0} yields
\begin{equation}\label{stepviequ1}
\forall\;c>0\;\forall\;\ell>0:\;\;\;\omega_{\frac{\mathbf{W}^{(c\ell_0)}}{\mathbf{G}^{\alpha}}}\hyperlink{sim}{\sim}\omega_{\mathbf{Q}^{(c)}}\hyperlink{sim}{\sim}\omega_{\widetilde{\mathbf{W}}^{(c)}}\widehat{\star}\omega_{\mathbf{G}^{\alpha}}\hyperlink{sim}{\sim}\omega_{\mathbf{W}^{(\ell)}}\widehat{\star}\omega_{\mathbf{G}^{\alpha}}\hyperlink{sim}{\sim}\omega\widehat{\star}\omega_{\mathbf{G}^{\alpha}}\hyperlink{sim}{\sim}\omega\widehat{\star}\id^{1/\alpha},
\end{equation}
which implies, in particular, that
\begin{equation}\label{stepviequ2}
\forall\;c,c_1>0:\;\;\;\omega_{\mathbf{Q}^{(c)}}\hyperlink{sim}{\sim}\omega_{\mathbf{Q}^{(c_1)}}.
\end{equation}
\vspace{6pt}

\emph{Step VI - Final conclusions}
By \emph{Step V} we have that \hyperlink{om1}{$(\omega_1)$} is valid for $\omega_{\mathbf{Q}^{(c)}}$ for any $c>0$ and $\mathcal{Q}$ satisfies \hyperlink{R-mg}{$(\mathcal{M}_{\{\on{mg}\}})$} and \hyperlink{B-mg}{$(\mathcal{M}_{(\on{mg})})$} as seen in \emph{Step II.}

Thus we have verified the assumptions to apply (the proof of) the main result \cite[Thm. 5.4]{equalitymixedOregular} to the matrix $\mathcal{Q}$ and get that as l.c.v.s.
\begin{equation}\label{stepviiequ}
\mathfrak{F}_{[\omega_{\mathcal{Q}}]}=\mathfrak{F}_{[\mathcal{Q}]},
\end{equation}
where $\omega_{\mathcal{Q}}:=\{\omega_{\mathbf{Q}^{(c)}}: c>0\}$ denotes the corresponding \emph{weight function matrix}; see \cite[Sect. 2.6 \& 2.7]{equalitymixedOregular} for precise definitions concerning the weighted spaces. Note the following:
\begin{itemize}
\item[$(*)$] The characterization in \cite[Thm. 3.2]{equalitymixedOregular} holds with the same index $x=y$ since each $\omega_{\mathbf{Q}^{(c)}}$ satisfies \hyperlink{om1}{$(\omega_1)$}.

\item[$(*)$] Recall that w.l.o.g. $\mathcal{Q}$ can be assumed to be standard log-convex, see \emph{Step II}, and hence $\mathcal{Q}$ is $(\mathcal{M}_{\on{sc}})$ in the notion of \cite[Sect. 2.5]{equalitymixedOregular}.

\item[$(*)$] Indeed, in \cite[Thm. 5.4]{equalitymixedOregular} even a characterization has been established in terms of \eqref{stepviiequ} for $\mathfrak{F}=\mathcal{E}$. However, since for this proof we are only interested in inferring the equality \eqref{stepviiequ} from growth properties on the weight matrix $\mathcal{Q}$, i.e. in having one half of \cite[Thm. 5.4]{equalitymixedOregular}, the desired equality holds for any symbol/functor under consideration. Analogously, this comment applies to the equalities stated in the proof below.

\item[$(*)$] For the converse and so for the full characterization each weighted setting requires crucially appropriate (optimal/special/characteristic) elements in the weighted spaces under consideration and functional analytic techniques \emph{(Grothendieck's factorization theorem, closed graph theorem).} In particular, the appropriate choice for these optimal elements depends on the weighted setting and this is inspired by the characterization of the inclusion relations; for more details we refer to \cite[Sect. 5]{equalitymixedOregular}, \cite{testfunctioncharacterization} and \cite[Prop. 4.6]{compositionpaper} in case $\mathfrak{F}=\mathcal{E}$.
\end{itemize}

By taking into account the second part of \eqref{stepviequ1} we can extend \eqref{stepviiequ} to
\begin{equation}\label{stepviiequ1}
\forall\;c>0:\;\;\;\mathfrak{F}_{[\omega_{\mathbf{Q}^{(c)}}]}= \mathfrak{F}_{[\omega_{\mathcal{Q}}]}=\mathfrak{F}_{[\mathcal{Q}]}.
\end{equation}
On the other hand, recall that by \eqref{stepviequ0} and definitions \eqref{stepivequnew}, \eqref{stepivequ0} in \emph{Step III} as l.c.v.s.
\begin{equation}\label{stepviiequ2}
\forall\;c>0:\;\;\;\mathfrak{F}_{[\omega_{\mathbf{W}^{(c\ell_0)}/\mathbf{G}^{\alpha}}]}=\mathfrak{F}_{[\omega_{\mathbf{Q}^{(c)}}]}=\mathfrak{F}_{[\omega_{\mathcal{Q}}]}=\mathfrak{F}_{[\mathcal{Q}]}=\mathfrak{F}_{[\frac{\widetilde{\mathcal{M}}_{\omega}}{\mathbf{G}^{\alpha}}]}=\mathfrak{F}_{[\frac{\mathcal{M}_{\omega}}{\mathbf{G}^{\alpha}}]}.
\end{equation}
\eqref{stepviiequ2} implies the first equality in \eqref{mainweighfctcorinvequ}, whereas the second, third and fourth one there follow by \eqref{stepviequ1}. Finally, the last equality in \eqref{mainweighfctcorinvequ} holds by \eqref{upperLegendregeneral}.
\qed\enddemo

Concerning the analogous inverse of Theorem \ref{mainweighfctthm} we can show:

\begin{theorem}\label{mainweighfctinv}
Let $\sigma,\tau\in\hyperlink{omset0}{\mathcal{W}_0}$ be given with associated matrices $\mathcal{M}_{\sigma}:=\{\mathbf{S}^{(\ell)}: \ell>0\}$ and $\mathcal{M}_{\tau}:=\{\mathbf{T}^{(\ell)}: \ell>0\}$. Let us assume the following:
\begin{itemize}
\item[$(i)$] For all $\ell>0$ there exists $\mathbf{Q}^{(\ell)}$ which is normalized and log-convex and satisfies $\mathbf{Q}^{(\ell)}\hyperlink{approx}{\approx}\frac{\mathbf{S}^{(\ell)}}{\mathbf{T}^{(\ell)}}$.

\item[$(ii)$] $\tau\hyperlink{omvartriangle}{\vartriangleleft}\sigma$ holds.

\item[$(iii)$] $0\le\overline{\gamma}(\tau)<\gamma(\sigma)\le+\infty$.
\end{itemize}
Then, for the matrix $\frac{\mathcal{M}_{\sigma}}{\mathcal{M}_{\tau}}$ from \eqref{quotientmatrix} and with $\mathcal{Q}:=\{\mathbf{Q}^{(\ell)}: \ell>0\}$ it follows that as l.c.v.s.
\begin{equation}\label{mainweighfctinvequ}
\forall\;\ell>0:\;\;\;\mathfrak{F}_{[\frac{\mathcal{M}_{\sigma}}{\mathcal{M}_{\tau}}]}=\mathfrak{F}_{[\mathcal{Q}]}=\mathfrak{F}_{[\omega_{\mathbf{S}^{(\ell)}}\widehat{\star}\omega_{\mathbf{T}^{(\ell)}}]}=\mathfrak{F}_{[\sigma\widehat{\star}\tau]},
\end{equation}
with $\mathfrak{F}\in\{\mathcal{E}, \mathcal{B}, \mathcal{A}, \mathcal{S}, \Lambda, \mathcal{F}\}$.
\end{theorem}

\demo{Proof}
The first equality in \eqref{mainweighfctinvequ} is clear by assumption $(i)$. Then assumption $(iii)$ implies, in particular, that $\gamma(\sigma)>0$ which is equivalent to \hyperlink{om1}{$(\omega_1)$} for $\sigma$. Hence, by $(i)$ in Lemma \ref{om1om6lemma1} and assumption $(ii)$, one infers $\mathcal{M}_{\tau}\vartriangleleft\mathcal{M}_{\sigma}$. Consequently, $\mathbf{T}^{(\ell)}\hyperlink{mtriangle}{\vartriangleleft}\mathbf{S}^{(\ell)}$ for any $\ell>0$ holds and so $\lim_{p\rightarrow+\infty}(Q^{(\ell)}_p)^{1/p}=+\infty$ is valid. Therefore, $\mathcal{Q}$ is standard log-convex (again by taking into account the technical comments concerning the point-wise order of the sequences $\mathbf{Q}^{(\ell)}$ made in \emph{Step II} in the previous proof).

By $(iii)$ one can apply (both parts of) \cite[Thm. 4.14]{genLegendreconj}. However, by $\overline{\gamma}(\tau)<+\infty$ it follows that $\tau$ satisfies \hyperlink{om6}{$(\omega_6)$} and so $\mathcal{M}_{\tau}$ is constant. This fact implies $\frac{\mathcal{M}_{\sigma}}{\mathcal{M}_{\tau}}\{\approx\}\frac{\mathcal{M}_{\sigma}}{\mathbf{T}^{(1)}}$, $\frac{\mathcal{M}_{\sigma}}{\mathcal{M}_{\tau}}(\approx)\frac{\mathcal{M}_{\sigma}}{\mathbf{T}^{(1)}}$ and assumption $(i)$ simplifies as follows:
$$\forall\;\ell>0:\;\;\;\mathbf{Q}^{(\ell)}\hyperlink{approx}{\approx}\frac{\mathbf{S}^{(\ell)}}{\mathbf{T}^{(1)}}.$$
Then follow the proof in Theorem \ref{mainweighfctcorinv} with $\mathbf{G}^{\alpha}$ being replaced by $\mathbf{T}^{(1)}$.
\qed\enddemo

Summarizing, the only difference between Theorem \ref{mainweighfctinv} and Theorem \ref{mainweighfctcorinv} is that $\mathbf{G}^{\alpha}$ has to be replaced by $\mathbf{T}^{(1)}$.

On the other hand, assumption $(i)$ in Theorem \ref{mainweighfctinv} motivates to introduce and study a new growth index relating two (different) sequences and hence generalizing the known growth index by Thilliez which corresponds to the case when the sequence in the denominator coincides with a Gevrey-sequence. Related to this idea it seems to be natural to introduce the new notion of the growth index $\gamma(\mathcal{M})$ even for abstractly given matrices $\mathcal{M}$. This allows to formulate assumption $(i)$ in a concise and quantitative form and more applications involving $\gamma(\mathcal{M})$ are expected.

\section{An application of the generalized upper Legendre conjugate to dynamics}\label{applicationuppersection}
Analogously as in Section \ref{applicationlowersection} we apply the derived information to the resolvent operator on Gelfand-Shilov-type spaces; the next result is a second variant of the generalization of \cite[Thm. 4.6 $2.$]{ArizaFernandezGalbis25} and involves $\widehat{\star}$; it should be compared with Theorem \ref{ArizaFernandezGalbisthm}: Roughly speaking we divide the involved weight matrices point-wise by a Gevrey sequence $\mathbf{G}^{\alpha}$ (with $0<\alpha<2$).

\begin{theorem}\label{ArizaFernandezGalbisthminv}
Let $\sigma,\tau\in\hyperlink{omset0}{\mathcal{W}_0}$ be given with associated weight matrices $\mathcal{M}_{\sigma}$, $\mathcal{M}_{\tau}$. Let $0<\alpha<2$, $\mu\in\CC$ with $|\mu|>1$ and $\psi$, $R_{\mu}$ as in Theorem \ref{ArizaFernandezGalbisthm}. We assume the following properties:
\begin{itemize}
\item[$(a)$] There exists $\ell_0>0$ such that $\mathbf{S}^{(\ell_0)}$ satisfies $\gamma(\mathbf{S}^{(\ell_0)})>\alpha+1$,

\item[$(b)$] $\sigma(t)=o(t^{1/\alpha})$ as $s\rightarrow+\infty$ holds.

\item[$(c)$] The matrices $\mathcal{M}_{\sigma}$, $\mathcal{M}_{\tau}$ are related by
$$\exists\;k_0>0\;\forall\;\ell>0:\;\;\;\mathbf{T}^{(k_0)}\hyperlink{preceq}{\preceq}\mathbf{S}^{(\ell)}.$$
\end{itemize}
Then there exists $f\in\mathcal{S}_{((((\sigma^{\alpha})^{\star})^{\iota})^{1/\alpha})}(\RR)=\mathcal{S}_{(\sigma\widehat{\star}\id^{1/\alpha})}(\RR)=\mathcal{S}_{(\frac{\mathcal{M}_{\sigma}}{\mathbf{G}^{\alpha}})}(\RR)$ such that $R_{\mu}(f)\notin\mathcal{S}_{(\tau)}(\RR)$.
\end{theorem}

\demo{Proof}
We proceed like in the proof of Theorem \ref{ArizaFernandezGalbisthm}.\vspace{6pt}

\emph{Step I} By the assumptions $(a)$ and $(b)$ we can apply Theorem \ref{mainweighfctcorinv} to $\sigma$ and have that (as l.c.v.s.) $\mathcal{S}_{((((\sigma^{\alpha})^{\star})^{\iota})^{1/\alpha})}(\RR)=\mathcal{S}_{(\frac{\mathcal{M}_{\sigma}}{\mathbf{G}^{\alpha}})}(\RR)$; see \eqref{mainweighfctcorinvequ}.

Next, assumption $(a)$ and \cite[Cor. 4.6 $(i)$]{index} imply $\gamma(\omega_{\mathbf{S}^{(\ell_0)}})\ge\gamma(\mathbf{S}^{(\ell_0)})>\alpha+1$ and so $\gamma(\sigma)>\alpha+1$ too by taking into account \eqref{goodequivalenceclassic}.

Now we compute $\gamma((((\sigma^{\alpha})^{\star})^{\iota})^{1/\alpha})$: For this let us apply the first part of \cite[Prop. 2.21]{index} to the weight function $\sigma^{\alpha}$ and obtain $\gamma(((\sigma^{\alpha})^{\star})^{\iota})\ge\gamma(\sigma^{\alpha})-1$. Then \eqref{newindex5} applied to $\sigma$ gives $\gamma(\sigma^{\alpha})-1=\frac{1}{\alpha}\gamma(\sigma)-1$ and \eqref{newindex5} applied to $((\sigma^{\alpha})^{\star})^{\iota}$, which is a weight function as seen in \cite[Sect. 2.5]{index}, yields $\gamma((((\sigma^{\alpha})^{\star})^{\iota})^{1/\alpha})=\alpha\gamma(((\sigma^{\alpha})^{\star})^{\iota})$. Combining everything with assumption $(a)$ we obtain: $$\gamma((((\sigma^{\alpha})^{\star})^{\iota})^{1/\alpha})=\alpha\gamma(((\sigma^{\alpha})^{\star})^{\iota})\ge\alpha(\gamma(\sigma^{\alpha})-1)=\gamma(\sigma)-\alpha>1.$$

Therefore, $(((\sigma^{\alpha})^{\star})^{\iota})^{1/\alpha}$ satisfies \eqref{assostrongnq}, recall again \cite[Cor. 2.13]{index}. In view of \eqref{upperLegendregeneral} and \eqref{stepviequ1}, also $\gamma(\omega_{\mathbf{Q}^{(c)}})>1$ for some/each $c>0$ with $\mathbf{Q}^{(c)}$ the sequences being defined in \eqref{stepiiiequ} in the proof of Theorem \ref{mainweighfctcorinv}. Thus the Borel map $B: \mathcal{S}_{(\omega_{\mathbf{Q}^{(c)}})}(\RR)\rightarrow\Lambda_{(\omega_{\mathbf{Q}^{(c)}})}$ is surjective because some/each $\omega_{\mathbf{Q}^{(c)}}$ is a strong weight function.

Consequently, again by \eqref{mainweighfctcorinvequ}, $B: \mathcal{S}_{(\frac{\mathcal{M}_{\sigma}}{\mathbf{G}^{\alpha}})}(\RR)\rightarrow\Lambda_{(\frac{\mathcal{M}_{\sigma}}{\mathbf{G}^{\alpha}})}$ is also surjective and recall that the point $0$ can be replaced by any $x_0\in\RR$ (translation).\vspace{6pt}

\emph{Step II} Proceed as in \emph{Step II} in the proof of Theorem \ref{ArizaFernandezGalbisthm}: For any $n\in\NN$ let $\delta_n:=(\delta_{n,j})_{j\in\NN}$ and consider a strictly increasing sequence of positive integers $(\ell_j)_{j\in\NN_{>0}}$. Since by assumption $(c)$ we have $\mathbf{T}^{(k_0)}\hyperlink{preceq}{\preceq}\mathbf{S}^{(\ell)}$ for any $\ell>0$ (and $T^{(k_0)}_0=1=S^{(\ell)}_0$) the following holds:
$$\forall\;j\in\NN_{>0}\;\exists\;h_j\ge 1\;\forall\;p\in\NN:\;\;\;T^{(k_0)}_p\le h_j^pS^{(1/\ell_j)}_p.$$
The comments concerning the choices for $(h_j)_{j\in\NN_{>0}}$ and $(p_j)_{j\in\NN_{>0}}$ from \emph{Step II} in the proof of Theorem \ref{ArizaFernandezGalbisthm} apply as well and \eqref{ArizaFernandezGalbisthmStepIIequ} is valid. Next introduce sequences $\mathbf{a}:=(a_n)_{n\in\NN}$ and $(\mathbf{b}_n)_{n\in\NN}$ with $\mathbf{b}_n=(b_{n,j})_{j\in\NN}$ given by
\begin{equation}\label{ArizaFernandezGalbisthmStepinvIIequ1}
a_k:=S^{(1)}_k,\;\;\;0\le k<p_1,\hspace{15pt}a_k:=S^{(1/\ell_j)}_kk!^{-\alpha},\;\;\;p_j\le k<p_{j+1},\;j\in\NN_{>0},
\end{equation}
and the point-wise product
$$\mathbf{b}_n:=\delta_n\cdot\mathbf{a}.$$
Note that $(\mathbf{b}_n)_{n\in\NN}$ is a bounded sequence in $\Lambda_{((((\sigma^{\alpha})^{\star})^{\iota})^{1/\alpha})}=\Lambda_{(\frac{\mathcal{M}_{\sigma}}{\mathbf{G}^{\alpha}})}$: The estimate
$$\forall\;\ell>0\;\exists\;C\ge 1\;\forall\;n\in\NN\;\forall\;j\in\NN:\;\;\;|b_{n,j}|=b_{n,j}\le CS^{(\ell)}_jj!^{-\alpha}$$
clearly holds by the point-wise order of the sequences in the matrix $\mathcal{M}_{\sigma}$ and note that the exponential growth factor formally required for the definition of the space $\Lambda_{(\frac{\mathcal{M}_{\sigma}}{\mathbf{G}^{\alpha}})}$ can be omitted. This is due to the fact that the matrix $\frac{\mathcal{M}_{\sigma}}{\mathbf{G}^{\alpha}}$ satisfies \hyperlink{B-L}{$(\mathcal{M}_{(\on{L})})$}: $\mathcal{M}_{\sigma}$ shares this property because $\gamma(\sigma)>0$ and by $(iv)$ in Section \ref{assomatrixsection}. Obviously \hyperlink{B-L}{$(\mathcal{M}_{(\on{L})})$} is preserved when dividing each sequence in the matrix point-wise w.r.t. a fixed sequence.\vspace{6pt}

\emph{Step III} As in \emph{Step III} in the proof of Theorem \ref{ArizaFernandezGalbisthm} let $x_0\ge 2$ be arbitrary but fixed and put recursively $x^2_{n+1}:=x_n-\frac{1}{4}$. Since the Borel map $B: \mathcal{S}_{((((\sigma^{\alpha})^{\star})^{\iota})^{1/\alpha})}(\RR)\rightarrow\Lambda_{((((\sigma^{\alpha})^{\star})^{\iota})^{1/\alpha})}$ is surjective there exists a bounded sequence of functions $(f_n)_{n\in\NN}$ in $\mathcal{S}_{((((\sigma^{\alpha})^{\star})^{\iota})^{1/\alpha})}(\RR)$ such that $f^{(j)}_n(x_0)=b_{n,j}$ for all $j,n\in\NN$. As in \emph{Step III} in the proof of Theorem \ref{ArizaFernandezGalbisthm} for this conclusion we apply here the fact that $\Lambda_{((((\sigma^{\alpha})^{\star})^{\iota})^{1/\alpha})}$ is a \emph{Fr\'{e}chet nuclear space} via \cite[Thm. 4.6 $1.$]{ArizaFernandezGalbis25} and \cite[Lemma 3.4]{ArizaFernandezGalbis24} which can be seen as follows: By \eqref{mainweighfctcorinvequ} and the proof of Theorem \ref{mainweighfctcorinv} we have that as l.c.v.s. $\Lambda_{((((\sigma^{\alpha})^{\star})^{\iota})^{1/\alpha})}=\Lambda_{(\omega_{\mathbf{W}^{(c\ell_0)}/\mathbf{G}^{\alpha}})}=\Lambda_{(\omega_{\mathbf{Q}^{(c)}})}$ for any $c>0$ and $\ell_0>0$ being the index appearing in assumption $(a)$. Each of these (associated) weight functions is a weight as considered in these works and satisfies \hyperlink{om1}{$(\omega_1)$}.

Again, w.l.o.g. one can assume that each $f_n$ has compact support contained in $(x_1,\psi(x_0))$ and so $f_n^{(k)}(x_j)=0=f_n^{(k)}(\psi(x_0))$ for all $k\in\NN$ and $j\in\NN_{>0}$.

By the \emph{closed graph theorem,} when assuming that $R_{\mu}(\mathcal{S}_{((((\sigma^{\alpha})^{\star})^{\iota})^{1/\alpha})}(\RR))\subseteq\mathcal{S}_{(\tau)}(\RR)$ (even as sets), then the map $R_{\mu}:\mathcal{S}_{((((\sigma^{\alpha})^{\star})^{\iota})^{1/\alpha})}(\RR)\rightarrow\mathcal{S}_{(\tau)}(\RR)$ is continuous and hence $(R_{\mu}(f_n))_{n\in\NN}$ is a bounded sequence in $\mathcal{S}_{(\tau)}(\RR)$ and consequently, when choosing the index $\ell=k_0$, then
\begin{equation}\label{ArizaFernandezGalbisthminvStepIIIequ}
\exists\;C>0\;\forall\;n\in\NN:\;\;\;|(R_{\mu}(f_n))^{(n)}(x_n)|\le C T^{(k_0)}_n.
\end{equation}
Recall that the proof of \cite[Thm. 4.6 $1.$]{ArizaFernandezGalbis25} implies
$$\forall\;n\in\NN:\;\;\;|(R_{\mu}f_n)^{(n)}(x_n)|>\frac{(n+2)^{2n}}{(4|\mu|)^{n+1}}b_{n,n}=\frac{(n+2)^{2n}}{(4|\mu|)^{n+1}}a_n.$$
Let now $n\in\NN_{>0}$ be such that $p_j\le n<p_{j+1}$ for some $j\in\NN_{>0}$. Then \eqref{ArizaFernandezGalbisthmStepinvIIequ1} gives $a_n=S^{(1/\ell_j)}_nn!^{-\alpha}$ and \eqref{ArizaFernandezGalbisthmStepIIequ} implies $T^{(k_0)}_n\le\log(n)^nS^{(1/\ell_j)}_n$. Combining this information with \eqref{ArizaFernandezGalbisthminvStepIIIequ} yields
$$\exists\;C>0\;\forall\;n,j\in\NN_{>0},\;p_j\le n<p_{j+1}:\;\;\;S^{(1/\ell_j)}_nn!^{-\alpha}(n+2)^{2n}\le C(4|\mu|)^{n+1}\log(n)^nS^{(1/\ell_j)}_n;$$
and by the estimates $(n+2)^{2n}\ge n^{2n}$ and $n!^{\alpha}\le n^{n\alpha}$ for all $n\in\NN_{>0}$ we infer
$$\exists\;C>0\;\forall\;n\in\NN_{>0},\;n\ge p_1:\;\;\;n^{n(2-\alpha)}\le C(4|\mu|)^{n+1}\log(n)^n.$$
Since $2-\alpha>0$ this estimate is impossible for all (large) $n$ which gives the contradiction.
\qed\enddemo

Theorem \ref{ArizaFernandezGalbisthminv} is related to \cite[Thm. 4.6 $2.$]{ArizaFernandezGalbis25} as follows: Let $d,d'>1$ be given with $d<d'<d+2$, set $\sigma=\tau:=\id^{1/d'}$ and so $\gamma(\mathbf{S}^{(\ell)})=\gamma(\sigma)=\gamma(\omega_{\mathbf{G}^{d'}})=\gamma(\mathbf{G}^{d'})=d'$ for all $\ell>0$ since $\sigma$ satisfies \hyperlink{om6}{$(\omega_6)$}; confirm $(iii)$ in Section \ref{assomatrixsection} and \eqref{commentdequ} in $(d)$ in Section \ref{technicaldivsection}. Recall that here both $\mathcal{M}_{\sigma}$ and $\mathcal{M}_{\tau}$ are constant: all sequences belonging to these matrices are equivalent to $\mathbf{G}^{d'}$ and so assumption $(c)$ in the Theorem is trivial.

Then put $\alpha:=d'-d$, hence $\alpha\in(0,2)$, $\gamma(\sigma)=d'=\alpha+d>\alpha+1$ and $\sigma(t)=t^{1/d'}=t^{1/(d+\alpha)}=o(t^{1/\alpha})$ which verifies assumptions $(a)$ and $(b)$. Moreover, $$(((\sigma^{\alpha})^{\star})^{\iota})^{1/\alpha}=\id^{1/d'}\widehat{\star}\id^{1/\alpha}\hyperlink{sim}{\sim}\omega_{\mathbf{G}^{d'}}\widehat{\star}\omega_{\mathbf{G}^{\alpha}}=\omega_{\mathbf{G}^{d'-\alpha}}=\omega_{\mathbf{G}^{d}};$$ and for this identity recall \cite[Example 2.9, $(4.10)$, Thm. 4.14 \& 5.11]{genLegendreconj}. Finally, recall that the matrix $\frac{\mathcal{M}_{\sigma}}{\mathbf{G}^{\alpha}}$ is constant: indeed, all sequences are equivalent to $\frac{\mathbf{G}^{d'}}{\mathbf{G}^{\alpha}}=\mathbf{G}^{d'-\alpha}=\mathbf{G}^{d}$.

\bibliographystyle{plain}
\bibliography{Bibliography}

\end{document}